\documentclass[envcountsect,orivec,runningheads,11pt]{llncs}
\usepackage{a4wide}
\usepackage{amsmath}
\usepackage{amssymb}
\usepackage{xspace}
\usepackage{graphicx}
\DeclareFontFamily{U}{matha}{\hyphenchar\font45}
\DeclareFontShape{U}{matha}{m}{n}{
      <10> gen * matha
      <10.95> matha10 <12> <14.4> <17.28> <20.74> <24.88> matha12
      }{}
\DeclareSymbolFont{matha}{U}{matha}{m}{n}
\DeclareMathSymbol{\obot}{2}{matha}{"6B}

\makeatletter
\DeclareRobustCommand*{\bfseries}{%
  \not@math@alphabet\bfseries\mathbf
  \fontseries\bfdefault\selectfont
  \boldmath
}   
\def\longenum{%
    \leftmargin0pt
    \ifnum \@enumdepth >3
        \@toodeep
    \else
        \advance\@enumdepth \@ne
        \edef\@enumctr{enum\romannumeral\the\@enumdepth}%
        \list{\csname label\@enumctr\endcsname}{%
            \usecounter{\@enumctr}%
            \labelwidth\z@\leftmargin0pt
            \itemindent\parindent\advance\itemindent\labelsep
        }%
    \fi
}

\makeatother
\newcommand{\impl}{\supset}
\newcommand{\dagg}{\dagger}

\newcommand{\lab}[1]{\label{#1}\mynote{#1}} 
\renewcommand{\lab}[1]{\label{#1}}
\newcommand{\myurl}[1]{{\rm\texttt{#1}}\xspace}
\newcommand{\reduceq}{\preceq}

\newcommand{\IR}{\mathbb{R}}
\newcommand{\IA}{\mathbb{A}}

\newcommand{\IV}{\mathbb{V}}

\newcommand{\IQ}{\mathbb{Q}}

\newcommand{\IC}{\mathbb{C}}

\newcommand{\IE}{\mathbb{E}}
\newcommand{\IF}{\mathbb{F}}
\newcommand{\IZ}{\mathbb{Z}}

\newcommand{\IN}{\mathbb{N}}
\newcommand{\adjoint}{\dagger}
\newcommand{\BSS}{BSS\xspace}
\newcommand{\calBP}{\operatorname{BP}}
\DeclareMathOperator{\Bor}{\normalfont{\texttt{\large|\!|}}}
\DeclareMathOperator*{\Band}{\normalfont{\texttt{\large\&\&}}}
\DeclareMathOperator{\Bneg}{\normalfont{\texttt{\large!}}}

\newcommand{\range}{\operatorname{range}}
\newcommand{\Real}{\operatorname{Re}}
\newcommand{\Imag}{\operatorname{Im}}
\renewcommand{\Re}{\Real}
\renewcommand{\Im}{\Imag}
\newcommand{\bin}{\operatorname{bin}}
\newcommand{\rank}{\operatorname{rank}}

\newcommand{\true}{{\rm\texttt{true}}\xspace}
\newcommand{\false}{{\rm\texttt{false}}\xspace}

\newcommand{\vep}{\varepsilon}

\newcommand{\DIM}{\operatorname{DIM}}

\newcommand{\SAT}{\operatorname{SAT}}
\newcommand{\UNSAT}{\overline{\SAT}}
\newcommand{\FEAS}{\operatorname{FEAS}}
\newcommand{\FEAD}{\FEAS^{\dagger}}
\newcommand{\QUAD}[1]{\operatorname{QUAD}_{1,#1}}
\newcommand{\QUART}[1]{\operatorname{QUART}^{\dagger}_{n\mapsto2n,#1}}
\newcommand{\sat}{\operatorname{sat}}
\newcommand{\unsat}{\overline{\sat}}
\newcommand{\Quantifier}{\operatorname{Q}}
\newcommand{\qf}{\Quantifier}

\newcommand{\Gr}{\operatorname{L}}
\newcommand{\Grl}{\Lat}
\newcommand{\Lat}{\operatorname{L}_{\text{l}}}
\newcommand{\Grf}{\operatorname{L}_*}
\newcommand{\Neu}{\operatorname{P}}
\newcommand{\Zero}{\text{\bf 0}}
\newcommand{\One}{\text{\bf 1}}

\newcommand{\id}{\operatorname{id}}

\newcommand{\Interpret}{\Theta}
 
\newcommand{\LNRarrows}{\genfrac{}{}{0pt}{}{\raisebox{-0.2ex}[0pt]{\!$\Leftarrow$}}{\raisebox{0.2ex}[0pt]{\,$\not\Rightarrow$}}}

\newcommand{\trinom}[3]{\Big(\!\!\begin{array}{c}\scriptstyle#1\\[-1.4ex]\scriptstyle#2\\[-1.4ex]\scriptstyle#3\end{array}\!\!\Big)}

\newcommand{\calO}{\mathcal{O}}
\newcommand{\calA}{\mathcal{A}}
\newcommand{\calB}{\mathcal{B}}
\newcommand{\calC}{\mathcal{C}}
\newcommand{\calU}{\mathcal{U}}

\newcommand{\calM}{\mathcal{M}}
\newcommand{\calMO}{\ensuremath{\mathcal{MO}}}
\newcommand{\Ring}{R}
\newcommand{\IntRing}{\calR}
\newcommand{\calQ}{\mathcal{Q}}
\newcommand{\calR}{\mathcal{R}}
\newcommand{\calS}{\mathcal{S}}

\newcommand{\calH}{\mathcal{H}}
\newcommand{\calL}{\mathcal{L}}

\newcommand{\calP}{\ensuremath{\mathcal{P}}}
\newcommand{\calNP}{\ensuremath{\mathcal{NP}}}
\newcommand{\calcoRP}{{\rm co}\mathcal{RP}}

\newcommand{\PSPACE}{\text{\sf PSPACE}\xspace}
\newcommand{\PAR}{\text{\sf PAR}\xspace}
\newcommand{\PAT}{\text{\sf PAT}\xspace}
\newcommand{\QSAT}{\operatorname{QSAT}}
\newcommand{\qsat}{\operatorname{qsat}}

\newcommand{\co}{\textsc{co}}
\newcommand{\PolyS}[1]{\Sigma^{\calP}_{#1}}
\newcommand{\PolyP}[1]{\Pi^{\calP}_{#1}}

\DeclareFontFamily{U}{bbold}{}
\DeclareFontShape{U}{bbold}{m}{n}
   {  <5> <6> <7> <8> <9> gen * bbold
      <10> <10.95> bbold10
      <12> <14.4> bbold12
      <17.28> <20.74> <24.88> bbold17
   }{}
\DeclareSymbolFont{bbold}{U}{bbold}{m}{n}
\DeclareMathSymbol{\BssD}{3}{bbold}{"01}
\DeclareMathSymbol{\BssS}{3}{bbold}{"06}
\DeclareMathSymbol{\BssP}{3}{bbold}{"05}
\newcommand{\BssPolyS}[1]{\PolyS{#1}}
\newcommand{\BssPolyP}[1]{\PolyP{#1}}

\newcommand{\person}[1]{\textsc{#1}}

\newcommand{\mycite}[2]{{\rm\cite[\textsc{#1}]{#2}}}

\newcommand{\Card}{\operatorname{Card}}
\newcommand{\End}{\operatorname{End}}
\newcommand{\Hom}{\operatorname{Hom}}

\newcommand{\COMMENTED}[1]{}
\newtheorem{observation}[theorem]{Observation}
\newtheorem{fact}[theorem]{Fact}
\newtheorem{scholium}[theorem]{Scholium}
\newtheorem{mydefinition}[theorem]{Definition}
\newtheorem{myremark}[theorem]{Remark}
\newtheorem{myexample}[theorem]{Example}
\newtheorem{myquestion}[theorem]{Question}
\newtheorem{mylemma}[theorem]{Lemma}
\newtheorem{myproposition}[theorem]{Proposition}
\newtheorem{convention}[theorem]{Convention}
\newtheorem{digression}[theorem]{Digression}
\newtheorem{comment}[theorem]{Comment}
\newtheorem{results}[theorem]{Results}
\renewcommand{\thefootnote}{\fnsymbol{footnote}}
\newcommand{\cx}[1]{#1} 
\newcommand{\cxx}[1]{#1} 

\begin{document}

\title{Computational Complexity of Quantum Satisfiability\thanks{
This work was initially supported by the \emph{German Research Foundation}
\texttt{DFG} under grant \texttt{Zi\,1009/2} 
and hosted by \person{Karl Svozil} at the Vienna University of Technology.
The authors thank 
\person{Manuel Bodirsky}, 
\person{Arno Pauly}, 
and \person{Peter Scheiblechner}
as well as local colleagues
for helpful discussions.
Further gratitude is due to anonymous referees for 
suggestions on improving the earlier version \cite{LiCS} --
and to \person{Martin Grohe} for his support.}}

\author{Christian Herrmann and Martin Ziegler}
\institute{Technical University of Darmstadt}
\date{}

\makeatletter
\renewcommand\maketitle{\newpage
  \refstepcounter{chapter}%
  \stepcounter{section}%
  \setcounter{section}{0}%
  \setcounter{subsection}{0}%
  \setcounter{figure}{0}
  \setcounter{table}{0}
  \setcounter{equation}{0}
  \setcounter{footnote}{0}%
  \begingroup
    \parindent=\z@
    \renewcommand\thefootnote{\@fnsymbol\c@footnote}%
    \if@twocolumn
      \ifnum \col@number=\@ne
        \@maketitle
      \else
        \twocolumn[\@maketitle]%
      \fi
    \else
      \newpage
      \global\@topnum\z@   
      \@maketitle
    \fi
    \thispagestyle{empty}\@thanks
    \def\\{\unskip\ \ignorespaces}\def\inst##1{\unskip{}}%
    \def\thanks##1{\unskip{}}\def\fnmsep{\unskip}%
    \instindent=\hsize
    \advance\instindent by-\headlineindent
    \if@runhead
       \if!\the\titlerunning!\else
         \edef\@title{\the\titlerunning}%
       \fi
       \global\setbox\titrun=\hbox{\small\rm\unboldmath\ignorespaces\@title}%
       \ifdim\wd\titrun>\instindent
          \typeout{Title too long for running head. Please supply}%
          \typeout{a shorter form with \string\titlerunning\space prior to
                   \string\maketitle}%
          \global\setbox\titrun=\hbox{\small\rm
          Title Suppressed Due to Excessive Length}%
       \fi
       \xdef\@title{\copy\titrun}%
    \fi
    \if!\the\tocauthor!\relax
      {\def\and{\noexpand\protect\noexpand\and}%
      \protected@xdef\toc@uthor{\@author}}%
    \else
      \def\\{\noexpand\protect\noexpand\newline}%
      \protected@xdef\scratch{\the\tocauthor}%
      \protected@xdef\toc@uthor{\scratch}%
    \fi
    \if@runhead
       \if!\the\authorrunning!
         \value{@inst}=\value{@auth}%
         \setcounter{@auth}{1}%
       \else
         \edef\@author{\the\authorrunning}%
       \fi
       \global\setbox\authrun=\hbox{\small\unboldmath\@author\unskip}%
       \ifdim\wd\authrun>\instindent
          \typeout{Names of authors too long for running head. Please supply}%
          \typeout{a shorter form with \string\authorrunning\space prior to
                   \string\maketitle}%
          \global\setbox\authrun=\hbox{\small\rm
          Authors Suppressed Due to Excessive Length}%
       \fi
       \xdef\@author{\copy\authrun}%
       \markboth{\@author}{\@title}%
     \fi
  \endgroup
  \setcounter{footnote}{\fnnstart}%
  \clearheadinfo}
\makeatother

\addtocounter{page}{-2}
\thispagestyle{empty}
\maketitle

\begin{abstract}
Quantum logic was introduced in 1936 by Garrett Birkhoff and
John von Neumann as a framework for capturing the logical
peculiarities of quantum observables.
It generalizes, and on 1-dimensional Hilbert space coincides with,
Boolean propositional logic. 

We introduce the weak and strong 
satisfiability problem for quantum logic terms.
It turns out that in dimension two both are also $\calNP$--complete.

For higher-dimensional spaces $\IR^d$ and $\IC^d$
with $d\geq3$ fixed, on the other hand, 
we show both problems to be complete 
for the nondeterministic Blum-Shub-Smale model of real computation.
This provides a unified view on both Turing
and real \BSS complexity theory; and 
extends the still relatively scarce family of $\calNP_\IR$--complete problems
with one perhaps closest in spirit to the
classical Cook-Levin Theorem.

Our investigations on the dimensions a term is weakly/strongly
satisfiable in lead to satisfiability problems in 
indefinite finite and finally in infinite dimension.
Here, strong satisfiability turns out as polynomial-time
equivalent to the feasibility 
of noncommutative integer polynomial equations
over matrix rings.
\end{abstract}

\newpage

\setcounter{tocdepth}{3}
\begin{center}
\begin{minipage}[c]{0.95\textwidth}\vspace*{-8ex}%
\tableofcontents
\end{minipage}
\end{center}

\pagebreak

This contribution connects both discrete and algebraic complexity
theory with the satisfiability problem in certain non-Boolean lattices.
It thus aims at particularly broad an audience 
of pure mathematicians as well as of computer scientists.
We have therefore chosen to include some background on 
computational complexity theory
as well as on quantum logic---in particular 
to the  ``modular'' first  version of a ``logic of quantum
mechanics''  due to   \person{von Neumann} and \person{Birkhoff}
{\cite{Birkhoff}} based on the finite dimensional Hilbert spaces
considered in nowadays Quantum Computation.

This exposition employs, adapts, and tailors 
several general concepts from various areas of logic 
for our particular purpose; whereas few dedicated
\emph{digressions} offer additional background 
and deeper connections which may be skipped at
first reading.

\section{Introduction} \lab{s:Intro}
Quantum physics is famous for its seemingly paradoxical (yet very real) effects.
Shaping them into a mathematically sound physical theory 
was a big achievement of the last century \cite{Neumann}.
There, physical \emph{observables} correspond to 
linear self-adjoint operators on some Hilbert space $\calH$;
and \emph{properties} (i.e. observables attaining only (eigen)values 0 or 1)
to projection operators---which in turn can be identified 
with closed subspaces of $\calH$. 
Their logical features (reflecting  non-commutativity of operators) have been captured 
\cite{Mackey,Piron,Beltrametti}
as abstract properties  of the \emph{ortholattice} 
of closed  subspaces; cmp. also \cite{Rota}.  
Let us start with the following

\begin{mydefinition} \lab{d:qlogic}
Fix a  $\ast$-subfield $\IF$ of the complex number field $\IC$,
i.e. a subfield closed under and equipped with
complex conjugation 
 $(a+bi)^*=\overline{a+bi}=a-bi=\Re(a+ib)-i\Im(a+ib)$ as involution
(which may be ignored if $\IF \subseteq \IR$).
 Popular cases are e.g. $\IF=\IC$ itself, $\IF=\IA$ algebraic numbers,
$\IF=\IR$ real numbers, $\IF=\IQ$ rationals, and $\IF=\IA\cap\IR$
algebraic reals.
An $\IF$-\textsf{unitary} space is an    $\IF$-vector space $\calH$
equipped with a unitary scalar product, i.e. a map $(\vec x,\vec y)
\mapsto \langle \vec x \mid  \vec y \rangle \in \IF$,
which is $\IF$-linear for fixed $\vec x$ and
satisfies $\langle \vec y \mid  \vec x \rangle
=\overline{ \langle \vec x \mid  \vec y \rangle}$ 
and $ \langle \vec x \mid  \vec x \rangle \neq 0$ for
$\vec x \neq \vec 0$ (cf. \mycite{\S 59-65}{Halmos},\mycite{\S
2-3}{Gelfand}, \mycite{Ch.6}{Axler}, or \mycite{Ch.1}{Farenick}).
But observe that we follow the convention of Physics 
and require linearity  in the right hand argument.  
The basic example is  $\IF^d$ where
 $\langle \vec x \mid  \vec y \rangle =\sum_{j=1}^d
 \overline{x}_iy_i$.
The orthogonality  relation on $\cal H$ is then defined by
 $\vec x \perp \vec y$  if and only if $\langle \vec x \mid  \vec y \rangle =0$.
A subset $U$ of $\cal H$ is \textsf{closed} 
if $U=U^{\bot\bot}$ where
$U^\bot=\{\vec x \in \calH\mid \forall \vec u \in U.\ \vec x \perp
\vec u  \}$ is the  linear subspace orthogonal to $U$.
We write $U \perp V$ if $U \subseteq V^\bot$, 
equivalently $V \subseteq  U^\bot$.
The  system of quantum logical
properties of $\calH$, shortly the  \textsf{quantum logic} 
of $\calH$,  consists of the
set $\Gr(\calH)$ of all closed subsets 
 of $\calH$
equipped with the following connectives:
\begin{description}
\item[$\wedge$:] $\Gr(\calH)\times\Gr(\calH)\to\Gr(\calH)$, $(U,V)\mapsto U\cap V$
\item[$\neg$:] $\Gr(\calH)\to\Gr(\calH)$, 
  $U\mapsto U^\bot$,
\item[$\vee$:] $\Gr(\calH)\times\Gr(\calH)\to\Gr(\calH)$,
  $(U,V)\mapsto (U\cup V)^{\bot\bot}$ 
\end{description}
and, in addition, with the constants $\Zero:=\{\vec 0\}\in\Gr(\calH)$ and $\One:=\calH\in\Gr(\calH)$.
\end{mydefinition}
Note that any closed  subset $U$ is  a linear subspace of $\calH$ 
and that $U \vee V =(U+V)^{\bot\bot}$ where 
 $U+V=\{\vec u+\vec v:\vec u\in U,\vec v\in V\}$ is the
least linear subspace of $\calH$ containing $U\cup V$.
Moreover, if $\calH$ is of finite dimension $d$,
the case we are primarily interested in,
then any linear subspace $U$ is closed
 whence $U \vee V= U+V$. 
(Indeed, if $\dim U=k$ then  $U^\bot$ has dimension 
$d-k$ being the solution set of the
independent linear equations $
  \langle \vec u_i \mid  \vec x \rangle$,
$\vec u_1, \ldots, \vec u_k$ a basis of $U$. 
By the same token,  $\dim U^{\bot\bot}=d-(d-k)=k$ 
and $U=U^{\bot\bot}$ from $U \subseteq U^{\bot\bot}$.)
By $\Gr_k(\calH)$  we denote the set of $k$-dimensional
subspaces of $\calH$.

\cx{
\begin{figure}[htb]
\begin{center}
\includegraphics[width=0.85\columnwidth]{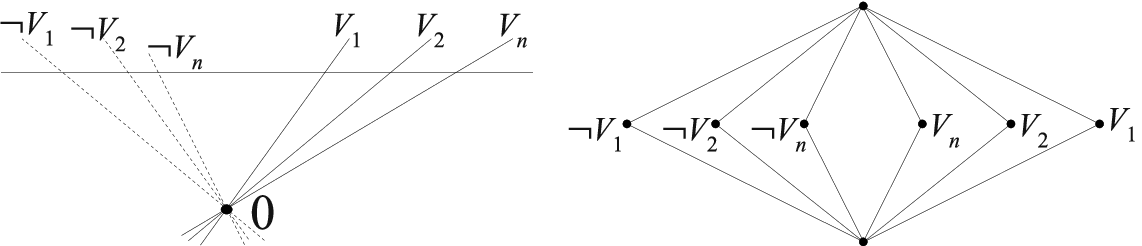}
\caption{\label{f:MO}Prototype $\calMO_{n}$ of a non-distributive
modular ortholattice}
\end{center}
\end{figure}}

\begin{myremark} \lab{r:Hilbert}
\begin{longenum}
\item[a)]
In case $d=1$, $\Gr(\calH)=\{\Zero,\One\}$ coincides with
the set of 
Boolean truth values. However, starting with dimension 2,
the distributive law ``$X\vee(Y\wedge Z)=(X\wedge Y)\vee(X\wedge Z)$''
generally fails: consider, e.g., $X:=\{(t,t):t\in\IF\}\in\Gr_1(\calH)$,
$Y:=\{(t,0):t\in\IF\}\in\Gr_1(\calH)$, 
$Z:=\{(0,t):t\in\IF\}\in\Gr_1(\calH)$; 
and compare also Figure~\ref{f:MO}.
This is generally seen as one cause underlying the
counter-intuitive effects of quantum physics.
\item[b)]
It is well-known that 
 $\big(\Gr(\calH),\wedge,\vee,\Zero,\One\big)$ constitutes 
an \textsf{lattice} 
a commutative, associative algebraic structure with respect to
idempotent operations ``$\wedge$'' (called \textsf{meet}) and ``$\vee$'' 
(called \textsf{join})
such that  $A\wedge(A\vee B)=A=A\vee(A\wedge B)$ and $0\wedge A=\Zero$,
$\One \wedge A=A$
(alternatively, $A \leq B \Leftrightarrow A=A \wedge B$ 
defines a   partial order  
 such that  $A \vee B$ is the supremum and
 $A\wedge B$  infimum  of $A,B$ and $\Zero$ and $\One$ are
the smallest and greatest element, respectively,
and then $A \leq B \Leftrightarrow B=A \vee B$). 
Moreover, $X \mapsto \neg X$ is an \textsf{involution} on 
$L$, i.e. $A \leq B$ if and only if $\neg B \leq \neg A$
and $\neg\neg A=A$. 
The \textsf{de Morgan laws} hold:
\begin{equation} \label{e:deMorgan}
\neg(A\vee B)=\neg A\wedge\neg B, \qquad
\neg(A\wedge B)=\neg A\vee\neg B. \end{equation}
In addition one has $A\wedge\neg A=\Zero$ and   $A\vee\neg A=\One$.
All this summarizes into saying that 
$\big(\Gr(\calH), \wedge,\neg,\vee,\Zero,\One\big)$ is
an \textsf{ortholattice}
and follows easily from the fact that $\perp$
is a symmetric  relation such that $\vec x \perp \vec x$ only for
$\vec x= \vec 0$.
\item[c)]
The crucial feature of finite dimensional $\calH$
is  the \textsf{modular law} 
\begin{equation} \label{e:Modular}
A \geq  C 
\quad\Rightarrow\quad  A\wedge (B \vee C) = (A \wedge B) \vee C  
\enspace .
\end{equation}
Indeed, if $\vec a \in   A\cap (B + C)$ then
$\vec a= \vec b +\vec c$ with $\vec b \in B$ and $\vec c \in C$
whence $\vec b =\vec a -\vec c \in A \cap B$,
proving $\leq$, while $\geq$ holds in any lattice.
An ortholattice obeying the modular law
is a \textsf{modular ortholattice}, MOL for short.
Actually, finite dimension is also necessary
for modularity of $\Gr(\calH)$ {\rm\cite{Birkhoff}}
and we restrict to this case
unless stated otherwise;
here,  $\Gr(\calH)$ amounts to the Grassmannian
of $\calH$ 
considered as an  ordered algebraic   structure  rather than as a topological space.
\end{longenum}
\end{myremark}

\cx{
\begin{digression}
In general, quantum mechanics lives in an \emph{in}finite-dimensional 
complex Hilbert space $\calH$; here,  only
the  special case $B=\neg C$ of the modular law,
the  \textsf{orthomodular law},
$A\geq C \Rightarrow  A=  (A\wedge\neg C)\vee C$
can be shown. Orthomodular lattices and their
generalizations  underly 
most of the work done in quantum logic.
Though, in {\rm\cite{Birkhoff,Neumann3}} \person{von Neumann}
stipulated certain modular sub-ortholattices  of $\Gr(\calH)$ 
(related to finite von Neumann algebras of operators) 
as the proper setting of a `logic of quantum mechanics'
cf. {\rm\cite{Duck,Redei}}. In the context of
the present investigation, these can be captured 
by the finite dimensional case---see Section~\ref{s:Neumann}.
These modular ortholattices, 
endowed with  a system of von Neumann  frames
corresponding to a continuous set
of matrix units {\rm\cite{Neumann4}}  respectively 
 tensor  products of
von Neumann  frames  {\rm\cite{Czedli}},
 might provide the  framework 
to take into account  
newer developments which  focus     on processes
and quantum  information and computation
\cite{OperationalQL,Coecke}.   
\end{digression}}
%
\subsection{Motivations from Physics, Logic, and Geometry}
We find quantum logic a promising object of study for various reasons
and to various communities \cite{Svozil,OperationalQL,Handbook}.

Physically,
it arises from quantum mechanics which our world seems to be
based on \cite{Weizsaecker} yet of which we regularly 
perceive only macroscopic approximations---which would
mean that Boolean logic is strictly speaking `wrong'
and to be replaced with quantum logic \cite{Putnam}:
perhaps even in the foundations of mathematics?

Indeed, $\Gr(\calH)$ constitutes a logic with semantics
$\false=\Zero$ and $\true=\One$ plus
(in dimensions$>1$) `intermediate' values like 
$V_n:=\IF\cdot\binom{1}{n}:=\{(t,n\cdot t):t\in\IF\}$, 
$n\in\IN=\{1,2,\ldots\}$. 
Such many-valuedness resembles probabilistic/fuzzy logic
and, more generally, G\"{o}del Logics
which `interpolate' between 0 (impossibility) and 1 (certainty) with 
values $r\in[0,1]$ (likelihood).
Note, though, that the above
$V_n$ are pairwise incomparable; cmp.
Figure~\ref{f:MO}a).

Geometry provides, at least in low dimensions, 
some important intuition to quantum logic. 
In fact the incidences 
among linear subspaces prescribed by some term relate 
to the field of Convex and Combinatorial Geometry 
\cite{Ziegler,Goodman}.

\subsection{Computational Complexity}
is at the core of theoretical computer science \cite{Papadimitriou}.
It classifies decision problems (i.e. countable families
of yes/no questions) according to the cost 
inherently incurred for their algorithmic solution:
asymptotically with respect to the binary length $n$ 
of the inputs and expressed in Landau's $\calO$-notation.
Both algorithms and costs pertain to some model of
computation; which usually amounts to the Turing machine
capable of reading, storing, processing, and printing
a constant number of bits in each step. This leads to
classical complexity classes $\calP$ (polynomial time),
$\PSPACE$ (polynomial space). 
The class $\calNP$ (nondeterministic $\calP$) 
on the other hand refers to polynomial-time Turing machines
that may `guess', but have to `verify' in polynomial time, bits. 

\begin{myexample} \lab{x:EvalSat}
\begin{enumerate}
\item[a)]
The following \emph{Boolean evaluation} problem 
is in $\calP$:
\begin{quote}
Given (the binary encoding $\langle t\rangle$ of)
a term $t$ over $\wedge,\vee,\neg,\Zero,\One$ in variables $x_1, \ldots ,x_n$ 
and an assignment within $\{\Zero, \One\}$, does it evaluate
to $\One$ ? 
\end{quote}
\item[b)]
The following \emph{Boolean satisfiability} problem,
$\SAT$, is in $\calNP$:
\begin{quote} 
Given (an encoding $\langle t\rangle$ of) a term $t$ over 
$\wedge,\vee,\neg$ and variables $x_1,\ldots,x_n$,
do they admit an assignment over $\{\Zero, \One\}$
for which $t$ evaluates to $\One$ ?
\end{quote}
\item[c)] 
One (out of many) reasonable way of encoding terms in a) and b) in binary
starts by assigning the function symbols $\Zero,\One,\neg,\vee,\wedge$ 
to bit strings $0,1,00,01,10$ and variable symbol $x_k$ ($k\in\IN$) 
to the number $k+2$ in binary. A term in polish notation thus gives
rise to a finite sequence of these  binary strings. And two such bit strings
$(a_1,\ldots,a_m)$ and $(b_1,\ldots,b_\ell)$ can be `concatenated'
reversibly into a single $(a_1,0,a_2,0,\ldots,a_m,1,b_1,0,b_2,0,\ldots,b_\ell,1)$.
\end{enumerate}
\end{myexample}
(\emph{Karp} or \emph{many-one}) reduction allows to compare
decision problems: $A\reduceq B$
('$A$ is at most as difficult as $B$')
means that instances for $A$ can efficiently be converted
into instances for $B$; formally: there exists a total 
function $f$ computable in polynomial time such that
$\forall x: x\in A\Leftrightarrow f(x)\in B$ holds.
Indeed, a (hypothetical) algorithm $\calM_B$ deciding membership
to $B$ gives rise to one at most polynomially slower for 
$A$: Given an instance ``$x\in A$?'', $\calM_A$
calculates $y:=f(x)$ and invokes $\calM_B$ to decide
whether $y\in B$ holds or not. A problem $B$ 
being $\calNP$-\emph{complete} means both that 
$B$ belongs to $\calNP$
(can be solved in nondeterministic polynomial time)
and is $\calNP$-\emph{hard}: every $A\in\calNP$ 
has $A\reduceq B$. In other words, $B$ is a
computationally most difficult problem within $\calNP$.
The \textsf{Cook-Levin Theorem} now states that
$\SAT$ is $\calNP$-complete---and thus not 
believed to admit a polynomial-time solution
by deterministic Turing machines. In fact a
vast variety of decision problems naturally
arising from various areas have all turned out
as polynomial-time equivalent to $\SAT$ \cite{Garey}.

Quantum computers have been suggested as an alternative
(and perhaps more powerful) model of computation.
They exploit linearity of quantum mechanics,
namely, that its evolution extends from states to superpositions.
Therefore a physical system realizing some `computation' on 
(say polynomially many) so-called qubits also works on 
linear combinations thereof---simultaneously: quantum parallelism.
The difficulty in exploiting this capability algorithmically 
consists in preparing the superposition and in extracting the 
output from the resulting state. 
Quantum logic, on the other hand, 
as describing operations on observables rather than states, 
has also been proposed as an approach 
to
computational purposes 
\cite{QLQC1,QLQC4,QLQC3}
 and   to   computational concepts \cite{Hagge1,QLQC2}.

\subsection{Blum-Shub-Smale Machines}
In algebraic complexity, the \BSS
machine (also known as \emph{real-RAM})
is a common model \cite{BCSS} 
capturing arithmetic on real numbers as entities 
with unit cost per arithmetic operation. 
More precisely it can 
read, store, operate, compare, and output a constant
number of reals in each step.

\begin{mydefinition} \lab{d:BSS}
Consider a commutative ring $\Ring$, possibly ordered.
A (deterministic) \textsf{\BSS machine $\calM$ over $\Ring$}
contains a finite number of constants\footnote{\label{f:BssConstant}%
Storing the non-recursive constant $\sum_{x\in H}2^{-n}$ may
be exploited by a \BSS-machine over $\IR$ to decide the
Halting problem $H\subseteq\IN$ for Turing (!) machines.
For machines without constants, refer to Fact~\ref{f:BSSBP}b).%
} $c\in\Ring$,
a register array, and three index registers.
It receives as input some finite tuple $\bar x\in\Ring^n$ 
together with its length $n=:|\bar x|\in\IN$.
$\calM$ can then apply arithmetic operations $+,-,\times$
(and $\div$ in case $\Ring$ is a field, but \emph{no}
conjugation $z\mapsto\bar z$ even in the complex case%
\footnote{This is the standard conception \mycite{\S 2.1}{BCSS},
noting that otherwise $\IC$ would become computationally
isomorphic to $\IR^2$ and violate the natural differences
between algebraic and semi-algebraic geometry \mycite{\S 2.3}{BCSS}.})
to these $x_j$, to its pre-stored constants,
or to some array elements accessed via index registers,
and store the result. It may furthermore branch 
based on the test for equality $=,\neq$ 
(in the ordered case also for $<,\leq,>,\geq$)
of two array elements. Each operation/branch
is counted for as one step.
On a fixed input $\bar x$, 
$\calM$ may accept, reject, or loop indefinitely.
\end{mydefinition}
Cf. \mycite{Definition~3.1}{BCSS} and compare also,
e.g., \mycite{\textsection 4.A}{Poizat}
or \mycite{\textsection 3}{Zucker}.
Note that operations and comparisons are presumed exact.

\begin{myexample} \lab{x:BSS}
\textsf{Gaussian Elimination} for $n\times n$--matrices over
a fixed field $\IF$ is a typical algorithm for \BSS-machines
over $\IF$ with polynomial running time $\calO(n^3)$.
Here, exact comparisons are employed during pivot search.

Furthermore this model commonly underlies algorithms 
devised, among others, in polynomial system solving {\rm\cite{Cox}}
and in computational geometry {\rm\cite{compGeom}}.
\end{myexample}
\BSS machines over the field $\IZ_2=\{0,1\}$ can be
seen equivalent to Turing machines. Definition~\ref{d:BSS}
thus extends the traditional, discrete theory of computation;
and has led to a rich structural complexity 
theory \cite{Meer}. In particular, 
a \emph{nondeterministic} \BSS machine over $\Ring$ 
may make and verify guesses from $\Ring$;
and doing so in polynomial time gives rise to
the complexity class $\calNP_\Ring$, thus 
naturally translating to this setting 
the classical question ``$\calP=\calNP$?''---which 
has turned out as equally inaccessible \cite{Fournier1,Fournier2}. 
As a matter of fact, both ``$\calP$ versus $\calNP$'' 
\emph{and} ``$\calP_{\IC}$ versus $\calNP_{\IC}$'' are propagated 
as \emph{Third Problem for the Next Century} \cite{Smale}. 

\begin{myexample} \lab{x:Feas}
Fix a commutative\footnote{Section~\ref{s:Indef} will 
naturally arrive at considering also noncommutative rings.}
ring $\Ring\supseteq\IZ$. The following problem $\FEAS_{\Ring,\Ring}$
can be decided by a nondeterministic polynomial-time \BSS machine
over $\Ring$, i.e. belongs to $\calNP_{\Ring}$:
\begin{quote}
Given 
($n\in\IN$ and the list of monomials and coefficients of each of)
finitely many polynomials $p_1,\ldots,p_k\in\Ring[X_1,\ldots,X_n]$, 
do they admit a common root in $\Ring$, i.e. some 
$\bar x\in\Ring^n$ such that $p_1(\bar x)=\cdots=p_k(\bar x)=0$ ?
\end{quote}
\end{myexample} 
Indeed, such a machine (does not need constants and) will simply `guess' 
an assignment $x_1,\ldots,x_n\in\Ring$ and `verify' 
by evaluating the polynomials:
which is clearly possible in a number steps
polynomial in (and noting that $n$ is bounded by)
the length of (the descriptions of) the polynomials.
$\FEAS_{\IC,\IC}$ is classically characterized by 
the famous \textsf{Hilbert's Nullstellensatz}
in algebraic geometry.

Generalizing the Cook-Levin Theorem, $\FEAS_{\Ring,\Ring}$ has been
established complete for $\calNP_{\Ring}$;
cf. \mycite{Theorem~5.1}{BCSS} and \cite{Goode}.
More precisely a Turing (!) machine can, given the description of a nondeterministic
polynomial-time \BSS machine $\calM$ (with symbolic references to its constants from $\Ring$)
and given an input $\vec y\in\Ring^m$,
output within time polynomial in $n$ some multivariate polynomials over $\Ring$
(with references to $\calM$'s constants as above)
such that the following holds:
$\calM$ accepts $\vec y$ ~iff~ these polynomials admit a common root in $\Ring$
\cite{Cucker,Meer}. The proof  employs 

\begin{fact} \lab{f:Constructible}
Finite Boolean combinations of polynomial in-/equalities
over a field can be expressed as (the feasibility of a system 
of) polynomial equations: 
\begin{enumerate}
\item[a)] $p(\bar x)=0\;\vee\;q(\bar x)=0 \quad\Leftrightarrow\quad
  \big(p\cdot q\big)(\bar x)=0$
\item[b)] $p(\bar x)\neq0 \quad\Leftrightarrow\quad \exists y:y\cdot p(\bar x)-1=0$
\end{enumerate}
\end{fact}
Note that, other than a \BSS machine,
even a nondeterministic Turing machine cannot in general
guess assignments over $\Ring$ in case the ring is infinite. 
Hence it is far from clear that $\FEAS_{\Ring,\Ring}$ even be decidable.
Not just in view of the Church-Turing Hypothesis
the question naturally arises of how the \BSS model
compares to the Turing model. Of course the latter is not
fitted to process inputs containing arbitrary, say, real numbers:
Instances of $\FEAS_{\Ring,\Ring}$ consist of both discrete (e.g. 
a list of monomials with their multi-degrees) and algebraic
(e.g. coefficients from $\Ring$) information, technically being 
words from a formal language over alphabet $\{0,1\}\cup\Ring$
while Turing machines work over $\{0,1\}$. 

\begin{myexample} 
The decision problem
$\big\{ \langle z,\bin(n)\rangle \mid z\in\IC, n\in\IN, z^n=1 \big\}$
for complex roots of unity belongs to $\calP_{\IC}$ but 
its instances cannot naturally be presented to a Turing machine.
\end{myexample}
This has suggested
restricting $\FEAS_{\Ring,\Ring}$, and more generally $\calNP_\Ring$,
to binary instances:

\begin{mydefinition} \lab{d:BP}
For a ring $\Ring\supseteq\IZ$, call
$\calBP(\calNP_\Ring):=\big\{\calL\cap\{0,1\}^*\mid
\calL\in\calNP_\Ring\big\}=\big\{\calL'\mid\calL'\subseteq\{0,1\}^*,
\calL'\in\calNP_\Ring\big\}$
the \textsf{Boolean part} of $\calNP_\Ring$;
cmp. \mycite{Definition~3.2}{Meer}.
\end{mydefinition}
Note that this 
complexity class is indeed closed under
Turing machine reduction.

\begin{myexample} \lab{x:BP}
For $\IZ\subseteq\Ring$, the following problem,
$\FEAS_{\IZ,\Ring}$, belongs to $\calBP(\calNP_\Ring)$:
\begin{quote}
\it Given (the coefficients, encoded in binary, of)
finitely many polynomials $p_1,\ldots,p_k\in\IZ[X_1,\ldots,X_n]$, 
do they admit a common root in $\Ring$ ?
\end{quote}
\end{myexample}
Although the input polynomials have integer coefficients,
their roots may in general `live' in $\Ring$:
consider for example $x^2-2$ or $x^2+1$.

\begin{fact} \lab{f:BSSBP} 
\begin{enumerate}
\item[a)]
$\FEAS_{\IZ,\Ring}$ 
is complete for $\calBP(\calNP_\Ring)$.
\item[b)]
Let $\calNP_\Ring^0$ denote the class of languages decidable
by a nondeterministic polynomial-time \BSS machine with\emph{out}
constants. It holds
$\calBP(\calNP_\IC^0)=\calBP(\calNP_\IC)$ 
\mycite{Proposition~7.9}{BCSS}
and $\calBP(\calNP_\IR^0)=\calBP(\calNP_\IR)$
\mycite{Theorem~4.1}{Buergisser}.
\item[c)]
It holds $\calNP\subseteq\calBP(\calNP_\IR)\subseteq\PSPACE$; 
cf. e.g. {\rm\cite{Grigoriev,Heintz,Canny,Matera}}.
\item[d)]
Subject to the \textsf{Generalized Riemann Hypothesis},
it holds $\calBP(\calNP_\IC)\subseteq\calcoRP^\calNP$;
cf. {\rm\cite{Koiran}}. 
\item[e)]
$\calBP(\calNP_\IZ)$ coincides with the class 
of all binary languages recursively enumerable
by a Turing machine;
cf. {\rm\cite{Matiyasevich}}.
\item[f)]
The decidability of $\FEAS_\IQ$ is
a long-standing open question of
extending Hilbert's tenth Problem 
from integers to rationals;
cf. e.g. {\rm\cite{Poonen}}.
\item[g)]
$\calBP(\calP_\IR^0)$ 
belongs to the Turing counting hierarchy {\rm\cite{Allender}}.
\end{enumerate}
\end{fact}
It remains an open challenge to tighten the relations
in Items~c)+d)+g).

In particular the 
class $\calBP(\calNP_\IR)$ has turned out
to be of interest of its own with several further
complete problems \cite{Shor,Rossello2,Zhang,Koiran99,Richter,Schaefer2}.
Higher \BSS complexity classes arise naturally for
problems in algebraic geometry \cite{BuergisserCuckerII}. 

We also mention \emph{analytic machines} as variants
of \BSS machines which may approximate their output
with/without error bounds \cite{Chadzelek}.

\subsection[\protect\pagebreak
Overview: Related Work, Objectives, Results, and Methods]{%
Overview: Related Work, Objectives, Results, and Methods}

Recently, there has been  new interest in
the classical topic of quantum logic,
the system of closed  subspaces  of a Hilbert space
and  axiomatic investigations  of  it;
see e.g. \cite{Mayet07,Megill}. 
Finite-dimensional  Hilbert spaces   are 
ubiquitous  in quantum computation,
which motivated \person{Hagge} et al. \cite{Hagge1,Hagge2} 
to  take up the  programme
of \person{Birkhoff} and \person{von Neumann} \cite{Birkhoff}
to study the associated  subspace ortholattices 
as  particularly tractable objects
of quantum logic. 
They pointed out  that  \person{Tarski}'s famous result applies
to  obtain decidability of the first 
order theory of  any single such ortholattice.
This has been used  in \cite{hn} 
to prove decidability of the equational
 theory of the class comprising all these ortholattices.
This theory  was  also shown to coincide
with  the equational theory of all  projection ortholattices
of finite von Neumann algebras. 

Our subject, here, is the counterpart of validity,
namely satisfiability of equations in a fixed $\Gr(\calH)$.
As in the Boolean case, satisfiability of a system can be 
compiled into the
(strong)  satisfiability of a single equation $t(\bar x)=\One$.
But, this is no longer the complement of validity 
of an identity $\forall \bar x\;t(\bar x)=\Zero$;  this complement 
$\exists\bar x\;t(\bar x)\neq\Zero$ 
will be called weak satisfiability.

Although computational complexity has become a standard
topic of investigation in logic since \cite{Cook}---cmp.
e.g. \cite{Graedel,Modal}---it seems to have passed on quantum logic.
We have taken upon this direction of research in \cite{LiCS}
and arrived at well-known complexity
classes $\calNP$, $\calNP_\IR$, and $\calBP(\calNP_\IR)$:
For these, the satisfiability problems for real \emph{and} complex
quantum propositional terms of appropriate dimensions
turn out to be complete. The present work recalls and 
presents the full proofs in a self-contained way.
One-dimensional quantum logic coinciding with
the classical Boolean one,
Section~\ref{s:2D} considers the satisfiability problem in 2D,
Section~\ref{s:Arith2QL} the case of dimensions
three and higher but fixed, \ref{s:Indef} the question
of satisfiability in \emph{some} finite dimension,
and \ref{s:Neumann} closes with some aspects of
the infinite-dimensional cases.

\begin{results}[] \lab{r:MainResults}
\begin{enumerate}
\item[a)] In fixed dimension,
weak and strong satisfiability (in general differ but) 
are polynomial-time equivalent (Theorem~\ref{t:WeakStrong})
\item[b)] Satisfiability in 2D quantum logic is
as hard as its classical Boolean (i.e. 1D) variant:
$\calNP$-complete, independent of the underlying field $\IF$
(Theorem~\ref{t:2DNPc})
\item[c)] Whereas starting with dimension three, (strong)
satisfiability over both real \emph{and complex} quantum logic
is complete for \emph{real} nondeterministic polynomial-time \BSS machines
(Theorem~\ref{t:IntStarField})
\item[d)] and remains so even when restricting to terms of the form 
  $\bigwedge\bigvee\bigwedge$ but becomes polynomial-time 
  decidable for $\bigwedge\bigvee$-terms (Theorem~\ref{t:Syntax}a+e).
\item[e)] Another syntactic variant complete for \emph{complex} nondeterministic
 polynomial-time \BSS machines is presented in Theorem~\ref{t:Syntax}c).
\item[f)] Satisfiability over rational 3D quantum logic is
 equivalent to Hilbert's tenth problem over $\IQ$ 
 (Corollary~\ref{c:realBSSNPc})
\item[g)]
 and validity over 3D rational quantum logic of a $\Sigma^0_3$-formula 
 is undecidable (Corollary~\ref{c:Poonen}). More generally, 
 quantified quantum logics correspond to the Boolean 
 and the \BSS polynomial hierarchy
 (Theorem~\ref{t:PolyHierarchy}).
\item[h)] Weak satisfiability over \emph{some} finite real or complex dimension
(i.e. asking for both a $d$ and a $d$-dimensional assignment)
is decidable by real nondeterministic polynomial-time \BSS machines
but not known hard (Theorem~\ref{t:WeakIndef}).
\item[j)] Strong satisfiability over some finite dimension is 
hard for polynomial-time \BSS nondeterminism
but not known decidable 
yet polynomial-time equivalent to the feasibility of
\emph{non}commutative polynomial equations
(Theorem~\ref{t:StrongIndef} and Proposition~\ref{p:StrongIndef}).
\end{enumerate}
\end{results}
We find that satisfiability in quantum logic may be
even more natural a generalization of the classical 
Boolean satisfiability problem than the feasibility
of a system of ring equations from Fact~\ref{f:BSSBP}a).
Moreover these results together provide a unified view 
on both Turing and real \BSS complexity theory.
They exhibit some resemblance to those concerning 
realizability questions for chirotopes \mycite{\S8}{Chirotopes} 
and to \textsf{descriptive complexity theory} 
where complexity classes are captured in
appropriate logics. Machine-independent 
characterizations of some \BSS complexity classes have
been obtained in \cite{GraedelMeer,Bournez}. 

We also emphasize the broad range of aspects of
logic in computer science
naturally joining in this contribution:
quantum logic, 
lattice theory, 
projective geometry, 
universal algebra, and computational complexity,
particularly in the Blum-Shub-Smale Model.

Our methods heavily draw from \person{John von Neumann}'s legacy.
We review, tailor, combine, and 
add new quantitative asymptotical and computational perspectives
to,  a variety of techniques known in quantum logic in order to yield
the polynomial-time reductions underlying the aforementioned results.
For reasons of self-containedness, proofs have been included also of 
the more technical tools and generally boil down their claims 
to facts in, e.g., basic linear algebra.

\section{Modular Quantum Logic}


Quantum Logic studies a wide spectrum of structures,
from the concrete $\Gr(\IC^d)$ to  abstract orthomodular lattices
and more general orthomodular structures (cmp. \cite{Handbook})
--- in  great variety of methods: from purely axiomatic 
to analysis of particular structures, with
representation results in between. 
The intermediate concept of \emph{Hilbert lattice}
comprises all  ortholattices of closed subspaces    
of generalized Hilbert spaces ---  the classical case
being characterized by 
\person{Sol\'{e}r} \cite{Prestel}.
In finite dimension $d$, Hilbert lattices are 
just  modular ortholattices as
introduced by
 \person{Birkhoff} and \person{von Neumann}   who showed 
 that  the irreducible ones  
correspond to   $d$-dimensional 
inner product spaces $\calH$ over  division $\ast$-rings $\IF$
(requiring $d \geq 4$ or the Arguesian Law). 
We deal with  the special  case where $\IF$ is a
$\ast$-subfield of $\IC$   and refer to Hilbert
lattices $\Gr(\calH)$  only in this context.
Some basic properties are quite well 
captured by the abstract concept of a finite
dimensional modular ortholattice.
We first recall some basic facts about modular lattices
and dimension --- and  give the easy proofs 
which  also shed some
new light on proofs in ordinary linear algebra.

\begin{fact} \lab{f:mod}  Let $L$ be a modular lattice with $\Zero$
  and $\One$.
In Items c)-e) suppose that $L$ 
contains a finite  maximal \textsf{chain}, that is a linearly
ordered subset.
\begin{enumerate} 
\item[a)]
For $u,v\in L$, the \textsf{interval} $[u,v]:=\{x:u\leq x\leq v\}$
is a sublattice and, in particular, modular.
\item[b)]
Between the intervals $[a \wedge b,\, a]$ and 
$[b,\,a \vee b]$  one has  mutually inverse  lattice
isomorphisms
 $x \mapsto b \vee x$ and $y \mapsto a \wedge y$
(this generalizes the canonical isomorphism
$A/(A\cap B) \cong (A+B)/B$  for linear subspaces $A,B$ of $\calH$). 
\item[c)]
All maximal chains $C$ in $L$ have the same
cardinality; $\dim(L):=|C|-1$ is  
 the \textsf{dimension} or \textsf{height} of $L$.
In particular, $\dim(u) :=\dim [\Zero ,u]$ is well defined for any 
$u \in L$. We call $u$ an \textsf{atom} if $\dim(u)=1$.
\item[d)]  
For $U \in \Gr(\calH)$,
$\dim(U)$  coincides with the vector space dimension.
\item[e)] One has the \textsf{dimension formulas}   
\[ \dim(v)= \dim(u)+\dim[u,v] \mbox{ for } u \leq v \;\mbox{ and } 
       \dim(a)+dim(b) =\dim(a\wedge b) +\dim(a\vee b) .\]  
\end{enumerate} 
\end{fact}

\begin{proof} 
a)   is obvious.  
b)  That the  two maps compose to identity, both ways,
is a simple application of modularity in each case.
Both maps being order preserving, it follows
 that they are  lattice isomorphisms.  

c) We use order induction,  considering all 
modular lattices $L$ admitting an $(n+1)$-element maximal chain $C$.
If $n=1$ then $L=\{\Zero ,\One \}$. Assume $n>1$. Consider any
$a \in C\setminus \{\Zero ,\One \}$, any maximal chain $C'$, and
 $b \in C'\setminus\{\Zero ,\One \}$. Then the inductive hypothesis
applies to any of the intervals $[\Zero ,\,a \wedge b]$,
$[a \wedge b,\,a]$,  $[a,\,a\vee b]$, and 
$[a\vee b,\,\One ]$. Due to the  isomorphisms, we have also
$\dim [a \wedge b,\,b] =\dim [a,\,a \vee b]$ and
 $\dim [b,\,a\vee b] =\dim [a\wedge b,\,a]$. 
This implies $|C'|=|C|$. 

d)  Given a basis $\vec v_1, \ldots ,\vec v_m$ of $U$ 
one obtains the $(m+1)$-element  maximal chain 
$U_0=\{\vec 0\}$,
$U_{i} =U_{i-1} +\IF\cdot\vec v_{i}$ \;$(i=1, \ldots ,m)$.  
   
 e) The first dimension formula
follows  since the union of two maximal
chains, one in  $[\Zero ,u]$, the other in $[u,v]$
is a maximal chain in $[\Zero ,v]$.
This, and the earlier   argument,  also   prove the second dimension
formula if one chooses for given $a,b$ maximal chains $C,C'$ in $L$
such that $a \wedge b,\,a,\,a\vee b \in C$ and
 $a \wedge b,\,b,\,a\vee b \in C'$. 
\end{proof}
Let us call the one-element lattice \emph{trivial}.
A homomorphism $\varphi:L\to L'$ between ortholattices is a
map having $\varphi(x\vee y)=\varphi(x)\vee\varphi(y)$
and $\varphi(x\wedge y)=\varphi(x)\wedge\varphi(y)$
and $\varphi(\neg x)=\neg\varphi(x)$.

\begin{fact} \lab{f:Foulis}
Abbreviate $C(x,y):=
(x\wedge y)\vee(x\wedge\neg y)\vee(\neg x\wedge y)\vee(\neg x\wedge\neg y)$,
the so-called \textsf{commutator}.
For $(L,\wedge,\vee,\neg,\Zero ,\One )$ an   ortholattice,
$a,b\in L$ are said to \textsf{commute} iff
$C(a,b)=\One $ holds. Since $C(a,b)=C(b,a)$, 
this means that $b,a$ commute, too.
Now, let $L$ be modular. 
\begin{enumerate} 
\item[a)] \textsf{Principle of Duality}:
If a statement is true for all modular
ortholattices, then so is its \textsf{dual} 
which arises by interchanging $\leq$ with $\geq$,
 $\wedge$ with $\vee$, and $\Zero $ with $\One $.
\item[b)]
Let $L$ be an MOL and $u \in L$.
Then $L_u=[\Zero ,u] \cup [\neg u,1]$ 
is a sub-ortholattice and
$([\Zero ,u],\wedge,\neg_u,\vee,\Zero ,u)$ is an MOL under the 
\emph{relative}
orthocomplement $\neg_u x = u \wedge \neg x$
and a homomorphic image of $L_u$  under the map 
$\phi(x)= u \wedge x$. 
For $U \in \Gr(\calH)$,
$[\Zero ,U] \subseteq \Gr(\calH)$ is (isomorphic to) the ortholattice $\Gr(U)$ 
given by  the scalar product induced on $U$.
\item[c)] For all $a,c \in L$: 
$C(a,b)=\One  \;\Leftrightarrow\;  a=(a\wedge b)\vee(a\wedge\neg b)
\;\Leftrightarrow\;   a =(a \vee b) \wedge (a \vee \neg b)
\;\Leftrightarrow\; 
a\wedge (\neg a \vee b) =a \wedge b \;\Leftrightarrow\; \neg a \wedge( a \vee b)= \neg a \wedge b$. 
In particular, comparable $a,b$ commute.
\item[d)] If $c$ commutes with all $x \in L$ 
(i.e. if $c$ is \textsf{central}) then  one has the isomorphism
\[ x \mapsto (x\wedge c,\;x \wedge \neg c),\;\;L \rightarrow
[\Zero ,c] \times [\Zero , \neg c]\]
with inverse $(y,z) \mapsto y \vee z$. 
In  particular, $L_c \cong [\Zero ,c] \times [\Zero , \neg c]$  for any $c\in L$. 
\item[e)]
If at least one of $a,b,c\in L$ commutes with the other two,
then the distributive laws 
$a\wedge(b\vee c) \;=\; (a\wedge b) \vee (a\wedge c)$ and
$a\vee(b\wedge c) \;=\; (a\vee b) \wedge (a\vee c)$
do apply. 
\item[f)]
If $a_1,\ldots,a_n\in L$ pairwise commute, 
they generate a Boolean algebra 
isomorphic to $\{\Zero ,\One \}^k$ for some $k\leq2^n$.
\item[g)]
The sublattice generated by $a,b\in L$ is 
isomorphic to $\{\Zero ,\One \}^k$ or to
$\{\Zero ,\One \}^k\times\calMO_2$ for some $k\leq 4$
where $\calMO_{n}$ is depicted in Figure~\ref{f:MO}b).
\end{enumerate}
\end{fact}
\begin{myremark} \lab{Beran} 
The first equivalence in c)  and e), f), g) 
also hold for orthomodular lattices --- see
 {\rm\mycite{\textsection3}{Beran}},
\mycite{Theorem~II.3.10}{Beran},
{\rm\cite[top of p.14 and \textsc{Theorem~II.4.5}]{Beran}}, and
\mycite{Theorem \textsection 3.9}{Kalmbach}.
\end{myremark}
Although the above  claims are folklore in the appropriate community
(and thus do not qualify for a lemma), 
we chose to include some proofs in the simpler case of 
modular $L$. These are to serve as illustration of 
reasoning in this non-Boolean logic, namely employing the
modular in place of the distributive law.  

\begin{proof} 
a) Given an ortholattice,  the mentioned exchange   
 yields an ortholattice due to the self-dual 
character of the axioms.  The same applies to modularity.

b) $x \leq y$ implies $\neg y \leq \neg x$  whence
$\neg_u y \leq \neg_ux$.  Also, $\neg_u(\neg_u x) 
= u \wedge  \neg(u \wedge \neg x) =
u \wedge (\neg u \vee x) = (u \wedge \neg u) \vee x =x$
by modularity. Thus, $x \mapsto \neg_u x$ is an involution on
$[\Zero ,u]$. Finally, $x \vee \neg_u x= 
u \wedge (x \vee \neg x)= u$ by modularity, i.e.
$[\Zero ,u]$ is an ortholattice.  
 
c)  We prove the second equivalence, first.
Indeed, in one direction, substituting $(a \wedge b)\vee (a \wedge \neg
b)$ for $a$ in the third equation one has
$a \leq (a \vee b)\wedge (a \vee \neg b)
=   \bigl((a \wedge \neg b) \vee b\bigr)\wedge 
 \bigl((a \wedge  b) \vee \neg b\bigr) 
= (a \wedge  b) \vee  \bigl( \bigl((a \wedge \neg b) \vee b\bigr)\wedge
\neg b \bigr)
= (a \wedge  b) \vee  (a \wedge \neg b) \vee (b\wedge
\neg b) \leq a $ by modularity. The converse direction follows by
duality.
Now, to prove the first equivalence  assume $C(a,b)=\One $.
Then, by modularity, $a = a \wedge C(a,b) 
= (a \wedge b) \vee (a \wedge \neg b) 
\vee \bigl(a \wedge\bigl( (\neg a \wedge b) \vee (\neg a \wedge  \neg b)
\bigr) \bigr)  = (a \wedge b) \vee (a \wedge \neg b)$.
Conversely, if $a= (a \wedge b) \vee (a \wedge \neg b)$
then by de Morgan
$\neg a = ( \neg a \vee \neg b) \wedge (\neg a  \vee b)$
whence, by the second equivalence,  $\neg a = (\neg a  \wedge b) \vee
(\neg a \wedge \neg b)$  and
$C(a,b)= a \vee \neg a =\One $. 
If $b,a$ commute then  $b=(b \wedge a)\ vee (b \wedge \neg a)$ 
whence  $a\wedge( b \vee \neg a) 
= a \wedge \bigl( (b \wedge a) \vee \neg a)\bigr)
= a \wedge b$ by   and modularity.
Conversely, if $a \wedge (\neg a \vee b)  =a \wedge b$ 
then by modularity  $b \leq (a \vee b) \wedge (\neg a \vee b) 
 = b \vee \bigl( a \wedge ( \neg a \vee b) \bigr)
= b \vee (a \wedge b)  =b$.   The last equivalence 
is due to the fact that $a,b$ commute if and 
only if $\neg a,b$ commute. Finally if, say, $a \leq b$ then
$\neg a \geq \neg b$ and 
$C(a,b) \geq \neg a \vee b \geq \One $.

d) Next, we show that if $c$  commutes with   $x,y$ then also
with  $x\vee y$, $x \wedge y$  and $\neg x$. Indeed,
$x \vee y \geq \bigl((x \vee y)\wedge c \bigr) 
\vee \bigl((x \vee y) \wedge \neg c \bigr) 
\geq  (x \wedge c) \vee (y \wedge c) \vee (x \wedge \neg c) \vee
(y \wedge \neg c)  \geq x \vee y$  and
$(\neg x \wedge c) \vee (\neg x \wedge \neg c) 
= \neg(x \vee \neg c) \vee \neg(x \vee  c) 
=  \neg \bigl((x \vee \neg c) \wedge (x \vee  c) \bigr) = \neg x$.
The claim about $x \wedge y$ follows by duality. 
 These  calculations   show that  the map
$x \mapsto (x \wedge c,\,x \wedge \neg c)$ is an homomorphism
 --- with inverse $(u,v) \mapsto u\vee v$.
Indeed, $\bigl((u \vee v) \wedge c\bigr)  \vee   
\bigl((u \vee v) \wedge \neg c\bigr)
= \bigl((u \vee (v \wedge c)\bigr)  \vee   
\bigl(( v \vee (u \wedge \neg c)\bigr)
 = u \vee v$ by modularity.

e)  Assume that $c$ commutes with $a$ and with $b$.
Then $c \wedge (a \vee b) = (c \wedge a)\vee (c \wedge b)$
due to that fact that $x \mapsto x \wedge c$ is a lattice
homomorphism --- the  isomorphism in g) followed by projection onto the 
first direct factor. 
It follows by c)  and modularity
$(a \wedge b)\vee (a \wedge c) \leq 
a \wedge (b \vee c) =
(a \vee c) \wedge (a \vee \neg c) \wedge \bigl( (b \wedge \neg c) \vee
c\bigr) = \bigl( b \wedge \neg c \wedge(a \vee c) \bigr)
\vee \bigl( c \wedge (a \vee \neg c)\bigr) 
= (b \wedge \neg c \wedge a) \vee (c \wedge a) 
\leq (a \wedge b) \vee (a \wedge c)$.
The remaining distributive relations follow by 
symmetry and duality. 
 
f)  The calculations in d)  also show that if $c$  commutes with 
all $x \in X$ then it is central in the sub-ortholattice
generated by $X \cup \{c\}$ (since $c$  commutes with  $c$). 
Thus, under the hypothesis of d),
any  two elements of the sub-ortholattice  $B$
generated by the $a_i$ commute with each other whence
$B$ is a Boolean algebra  by e).  

g) We show that $a$ commutes with $C(a,b)$
Indeed, $ a \wedge (\neg  a \vee  (C(a,b))
= a \wedge
 \bigl( \neg a  \vee (a \wedge b) \vee (a \wedge \neg b) \bigr)
=  (a \wedge \neg a) \vee (a \wedge b) \vee (a \wedge \neg b)
= (a \wedge b) \vee (a \wedge \neg b)
=a \wedge C(a,b)$  by modularity.      
By symmetry, $b$ commutes with $C(a,b)$, too.
Due to d) we have $L$ isomorphic to the
direct product $L_1 \times L_2$  where $L_i$ has generators
$a_i =\pi_i(a),\,b_i=\pi_i(b)$ 
 under the homomorphism
$\pi_1(x)= x \wedge c$  resp. $\pi_2(x)= x \wedge \neg c$ 
where  $c:=C(a,b)$ and $\pi_i$ denotes the projection 
onto the $i$-th component.
It follows  $C(a_1,b_1)= 
C(\pi_1(a),\pi_1(b))= \pi_1(c)=
\One _{L_1}$ and $C(a_2,b_2)= 
C(\pi_2(a),\pi_2(b)) =\pi_2(c)= \Zero _{L_2}$.
By f), $L_1$ is a Boolean algebra  with $2$ generators.
In $L_2$ we have $a_2 \wedge b_2 =a_2 \wedge \neg b_2
=\neg a_2 \wedge b_2 =\neg a_2 \wedge \neg b_2=\Zero $,
and, by de Morgan,  all joins $=\One $ which implies  that $L_2$ 
is either trivial or   a copy of 
$\calMO_2$.
\end{proof} 
In carrying out calculations in
modular ortholattices we often  write 
$\neg x = x ^\bot$,  
$x \wedge y = x \cap y$, 
and $x \vee y =x+y$ (and save brackets according to priority in
this order)  to  improve readability and
to appeal at the geometric intuition.
\cx{
\begin{digression} \lab{mol} 
 Modular ortholattices  of finite height are
isomorphic to direct products of irreducibles
cf.  \mycite{Ch. IV}{Birkhoff2}. 
The irreducibles  in turn can be understood 
in terms of irreducible projective spaces
with anisotropic orthogonality. Irreducibles of height
$d\geq 3$ are necessarily infinite \cite{Ball}, for $d\geq 4$ 
they correspond to inner product spaces over  division
$\ast$-rings \cite{Birkhoff}. On the other hand, there are examples where
the division ring has unsolvable word problem
which leads to a finitely presented MOL
having unsolvable word problem
\cite{Roddy}   and 
the undecidability of the equational theory  of MOLs of height $\leq d$ 
(where $d \geq 14$ is fixed) \cite{Micol}. 
\end{digression} }
%

\subsection{Truth, Equivalence, and Satisfiability} 
The classical Boolean satisfiability problem 
extends straightforwardly to quantum propositional 
terms---although \textsf{truth} (``$=\One$'') now
has to be distinguished from \textsf{non-falsity} (``$\neq\Zero$'');
and we may or may not permit constants (from a fixed MOL $L$).

\begin{mydefinition} \lab{d:QL}
\begin{enumerate}
\item[a)]
A (\textsf{ortholattice} or \textsf{quantum logic}) \textsf{term} 
is a syntactically correct expression 
over certain variables $x_1,\ldots,x_n$
with operations
 $\wedge$, $\neg$, and $\vee$ and  constants $\Zero $ and $\One $
(the latter  rather for formal reasons).
We may write $t(\bar x)$ to emphasize the
dependence on $(x_1,\ldots,x_n)=:\bar x$.
\\ The syntactic \emph{length} of $t$ is denoted by $|t|$,
defined recursively as $|x|=1$, $|\neg t|=|t|+1$, and
$|s\vee t|=|s|+|t|+1=|s\wedge t|$.
\item[b)]
For an ortholattice $(L,\wedge,\neg,\vee,\Zero ,\One )$
and $\bar a=(a_1,\ldots,a_n)\in L^n$, 
let $t_L(a_1,\ldots,a_n)
=t_L(\bar a)\in L$ denote the \textsf{value}
of $t$ in $L$ when substituting $a_i$ for $x_i$. 
\item[c)]
A term \textsf{with constants} from $L$
is one where some variables already
have been fixed to (i.e. substituted for)
elements of $L$.
\item[d)]
An $n$-variate term $t$ is \textsf{strongly satisfiable in $L$}
if there is $\bar a \in L^n$ such that
$t_L(\bar a)=\One $. It is \textsf{weakly satisfiable in $L$}
if there is $\bar a \in L^n$ such that
$t_L(\bar a)\neq\Zero $.
\item[e)]
Two $n$-variate terms $s$ and $t$ are 
\textsf{equivalent over $L$} if $s_L(\bar a)=t_L(\bar a)$
for every $\bar a\in L^n$.
\item[f)]
\textsf{Strong} and
\textsf{weak satisfiability over $L$} are the respective decision problems
\begin{eqnarray*}
 \SAT_L &:=& \big\{\langle t(x_1,\ldots,x_n)\rangle\mid
  n\in\IN, t\text{ term}, \exists \bar a\in L^n:t_L(\bar a)=\One\big\} 
\;\subseteq\;\{0 ,1\}^*
\qquad\text{and} \\ \sat_L &:=& \big\{\langle t(x_1,\ldots,x_n)\rangle\;\big|\;
  n\in\IN, t\text{ term}, \exists \bar a\in L^n:t_L(\bar a)\neq\Zero\big\} 
\;\subseteq\;\{0 ,1 \}^*
\end{eqnarray*}
\item[g)]
More generally, for a class $\calC$ of MOLs, consider the 
question of whether a given term $t$ is strongly/weakly satisfiable 
over \emph{some} $L\in\calC$:
$\SAT_\calC:=\bigcup_{L\in\calC}\SAT_L$,
$\sat_\calC:=\bigcup_{L\in\calC}\sat_L$.
\item[h)]
Returning to  single MOLs $L$,
strong satisfiability \emph{with constants from $C\subseteq L$}  is
\begin{multline*}
\SAT_{C,L}\;:=\;
\Big\{\big\langle t(x_1,\ldots,x_n,y_1,\ldots,y_m),c_1,\ldots,c_m\big\rangle\;\Big|
\\[-0.9ex]
  n,m\in\IN, t\text{ term}, c_1,\ldots,c_m\in C,
\exists \bar a\in L^n:t_L(\bar a,\bar c)=\One\Big\} 
\end{multline*}
and similarly for weak satisfiability with constants.
(The encoding $\langle\:\cdot\:\rangle$
of terms involving $c_i\in C$ will be specified later\ldots)
\end{enumerate}
\end{mydefinition}
Note that weak satisfiability of $t$
means \emph{in}validity of the identify ``$t=\Zero $''
in the model-theoretic sense. Moreover,
$\SAT_{\{\Zero ,\One \}}=\sat_{\{\Zero ,\One \}}$ 
coincides with the classical Boolean satisfiability problem.
We also record

\begin{observation} \lab{o:Product}
\begin{enumerate}
\item[a)]
For $U,V\in\Gr(\calH)$ with $U\bot V$ 
it holds that $\omega:\Gr(U)\times\Gr(V)\ni (A,B)\mapsto A\vee
 B\in\Gr(U\vee V)$ 
 is an embedding of the 
 product $\Gr(U)\times\Gr(V)$ of two ortholattices  into
 $\Gr(U\vee V)$.
\item[b)]
For the product $L\times L'$ of ortholattices $L$ and $L'$,
it holds $\SAT_{L\times L'}=\SAT_L\cap\SAT_{L'}$ and
$\sat_{L\times L'}=\sat_L\cup\sat_{L'}$.
\end{enumerate}
\end{observation}
For a) we may assume that $U \vee V =\calH$. Then
Fact~\ref{f:Foulis}d) applies.
Concerning b) note that, by the very definition of the product,
$\One_{L\times L'}=(\One_L,\One_{L'})$
and $\Zero_{L\times L'}=(\Zero_L,\Zero_{L'})$ 
and $t_{L\times L'}\big((x_1,y_1),\cdots,(x_n,y_n)\big)=\big(t_L(\bar x),t_{L'}(\bar y)\big)$.

\begin{myremark} \lab{r:Boolean}
The connective ``$\vee$''
satisfies the \emph{disjunction property} for weak truth:
$x\vee y\neq\Zero$ holds ~iff~ $x\neq\Zero$ or $y\neq\Zero$ holds.
In dimensions $>1$, however, strong truth generally fails this property:
$x\vee y=\One$ may well hold with neither $x=\One$ nor $y=\One$.
Similarly for the dual connective ``$\wedge$''.
Furthermore, Boolean negation has to be distinguished from complement:
$x\neq\Zero \LNRarrows \neg x=\Zero$.
\end{myremark}
Whenever this distinction needs emphasis,
we denote the traditional propositional connectives 
as in the \texttt{C} programming language:
$\Bor$ for ``\emph{or}'', $\Band$ for ``\emph{and}'',
$\Bneg$ for ``\emph{not}''. 
Lemma~\ref{l:Boolean2} below explains  how to
express them within \emph{quantified} quantum logic.

\begin{myexample} \lab{x:QL}
Fix some MOL $L$. 
\begin{enumerate}
\item[a)] $\neg x\vee\neg y$ and $\neg(x\wedge y)$
are terms over variables $x,y$.
According to the de Morgan laws,
they are equivalent.
\item[b)] Suppose $a\leq b$ and $a\vee \neg b=\One$
holds for some $a,b\in L$. Then $a=b$.
\item[c)] 
Recall the commutator $C(x,y)$ and, for $a,b\in L$, consider the 
sub-ortholattices $\{\Zero ,\One ,a,\neg a\}=:L(a)$ and $L(b)$ they span.
Whenever $a\in L(b)$ or $b\in L(a)$, 
it follows $C_L(a,b)=\One $.
In particular, $C(x,y)$ is equivalent to $\One$ on $\Gr(\IF^1)$.
Over $\Gr(\IF^2$), however,
$C(\IF\binom{1}{r},\IF\binom{1}{s})$ evaluates to $\Zero$
whenever $r\not\in\{s,-1/\bar s\}$. 
\item[d)] Let $t(x,y):=C(x,y)\vee x\vee y$ and
$s(x,y,z):=t(x,y)\wedge t(x,z)\wedge t(y,z)$.
Then $s$ is equivalent to 1 on $\Gr(\IF^2)$.
Over $\IF^3$, however, $s(
\IF\trinom{1}{0}{0},
\IF\trinom{1}{1}{0},
\IF\trinom{1}{1}{1})=\Zero$.  
\item[e)] 
More generally,
$\neg C(x,y)=(x\vee y)\wedge(\neg x\vee y)\wedge(x\vee\neg y)\wedge(\neg x\vee\neg y)$
is strongly satisfiable over each $\Gr(\IF^{2d})$, $d\in\IN$;
but not over $\Gr(\IF^{2d-1})$.
\item[f)] Over any MOL resp. complemented modular lattice it holds:
\[ \begin{array}{cccclcl} 
y\vee \neg(x\vee  y)=\One
&\Leftrightarrow&  x\leq y 
&\Leftrightarrow&
\exists z: \; &y \vee z=\One \;&\Band\; (x\vee y)\wedge z=\Zero \enspace ;\\
(x\wedge y)\vee\neg(x\vee y)=\One
&\Leftrightarrow&   x=y 
&\Leftrightarrow&
\exists z: \; &(x\wedge y)\vee z=\One \;&\Band\; (x\vee y)\wedge
z=\Zero \enspace . \end{array} \]
\end{enumerate}
\end{myexample}
\begin{proof}
Claim~b) is immediate by the orthomodular law.
Claim~c) is straightforward to verify. 
Claim~d) has been generalized in \cite{Hagge2} to (the complement of) terms
equivalent to $\One $ over $\IF^{d-1}$, but not over $\IF^{k}$ for every $k\geq d$;
cmp. Corollary~\ref{c:Hagge} below.
For Claim~e), consider $a:=\IF^d\times\{\Zero \}^d$ and
$b:=\{(\vec x,\vec x):\vec x\in\IF^d\}$ in the even-dimensional case.
Conversely, $A\vee\neg B=\One $ over $\IF^D$ implies 
$D=\dim(A\vee\neg B)\leq\dim(A)+\dim(\neg B)=\dim(A)+D-\dim(B)$;
hence $\dim(B)\leq\dim(A)$. And, by symmetry of $C(x,y)$, 
also $\dim(B)\leq\dim(\neg A)$, $\dim(A)\leq\dim(B)$,
and $\dim(A)\leq\dim(\neg B)$: Hence $D=\dim(A)+\dim(\neg A)$
must be even.
In   Claim~f), observe that  $y \leq x \vee y$ and that
$x \leq y$ if and only if 
$y =x \vee y$. Thus, the first equivalence follows from Claim~b)
(actually, it is valid in any orthomodular lattice
 \mycite{\S15 Theorem 3}{Kalmbach})
and the second by modularity.
The same arguing applies to the second line 
observing that $x \wedge y \leq  x \vee y$ 
and that equality holds if and only if $x=y$.
\qed\end{proof}
Example~\ref{x:QL}f) allows to reduce non-equivalence 
to weak satisfiability. More generally, we record

\begin{fact} 
\lab{f:clear} 
Let 
$s_1(\bar x),t_1(\bar x),\ldots,s_m(\bar x),t_m(\bar x)$ denote
$n$-variate terms and $L$ a MOL.
\begin{enumerate}
\item[a)]
An assignment $\bar a$ in $L$ simultaneously satisfies
all equations $s_i(\bar x)=t_i(\bar x)$, $1\leq i\leq m$,
~iff~ it strongly satisfies the single term
\[ \bigwedge\nolimits_{i=1}^m 
\big([s_i(\bar x)\wedge t_i(\bar x)]\vee\neg[s_i(\bar x)\vee t_i(\bar x)]\big) \enspace . \]
\item[b)]
The system $\{s_i(\bar x)=t_i(\bar x) \mid i=1, \ldots m\}$  
of equations is simultaneously satisfiable in $L$ ~iff~
the following two equations in $n+m$ variables $(\bar x,\bar y)$ are:
\[ \One=\bigwedge\nolimits_{i=1}^m \big([s_i(\bar x)\wedge t_i(\bar x)]\vee y_i\big)
\quad\Band\quad \Zero=\bigvee\nolimits_{i=1}^m \big([s_i(\bar x)\vee t_i(\bar x)]\wedge y_i\big)
\enspace . \] 
\end{enumerate}
\end{fact}

\subsection{$\calNP$--Hardness of Strong Satisfiability}
\begin{myproposition} \lab{p:NPhard} 
$\SAT_L$ is $\calNP$--hard for any nontrivial modular
ortholattice $L$, uniformly in $L$.
\end{myproposition}
In particular, for any non-empty class $\calC$ of nontrivial
MOLs, $\SAT_\calC$ is $\calNP$--hard. Theorem~\ref{t:2DNPc}
below shall extend this to weak satisfiability.

\begin{proof}
Convert a given term $t(x_1,\ldots,x_n)$
to the term  $s(\bar x):=t(\bar x)\wedge\bigwedge_{1\leq i<j\leq n} C(x_i,x_j)$
with the commutator from Fact~\ref{f:Foulis}:
this is clearly computable in polynomial time.
Moreover a satisfying Boolean assignment $\bar b \in \{\Zero,\One\}^n$ of $t$
is also one of $s$  in any non-trivial $L$
since $C(\Zero,\Zero)=C(\Zero,\One)=C(\One,\Zero)=C(\One,\One)=\One$.
Conversely a satisfying assignment of $g$ in $L$ 
consists of pairwise commuting elements $b_1,\ldots,b_n\in L$;
hence `lives' in a Boolean algebra (isomorphic to) 
$\{\Zero ,\One \}^k$ according to Fact~\ref{f:Foulis}f):
A satisfying Boolean assignment of $t$ is thus obtained by
projecting the $b_i$ onto their respective first components.
\end{proof} 

\begin{myexample}  \lab{x:2D}
Recall the modular ortholattices $\calMO_m$ from Figure~\ref{f:MO}b).
\begin{enumerate}
\item[a)]
Strengthening Example~\ref{x:QL}f),
$\bigwedge_{1\leq i<j\leq 2n } (x_i \vee x_j) \wedge (\neg x_i \vee
\neg x_j)$ is weakly/strongly
satisfiable in  $\calMO_m$ if and only if $m \geq n$.
\item[b)]
For $a\in\calMO_m$ with $m\geq2$, it holds 
$a=\One\Leftrightarrow\exists y,z:\neg C(y,z)\wedge C(a,y)\wedge C(a,z)\wedge a\neq\Zero$.
\end{enumerate}
\end{myexample} 
By a), both concepts of satisfiability thus depend on $L$
at least as far as finite two-dimensional $L$  are concerned.
The complexity of $\SAT_L$, however, will turn
out to not depend on $L$ as long as it has dimension two.

\subsection{Two-Dimensional Satisfiability is in $\calNP$}
\lab{s:2D}

\begin{mylemma} \lab{l:2D}
Let $L,L'$ be modular ortholattices 
of dimension $2$. 
\begin{enumerate}
\item[a)]
If $A \subseteq L$ then $A \cup \{\neg a\mid a \in A\} 
\cup\{\Zero,\One\}$  is a subortholattice of $L$.
In particular, $|L| \leq 2n+2$ if $L$ has
an $n$-element generating set.
Moreover,  $L$ is isomorphic to $L'$ if
and only if $|L|=|L'|$.
And
$L$ maps isomorphically onto a sub-ortholattice 
of $L'$ if and only if $|L| \leq |L'|$.
\item[b)]
A term $t(x_1,\ldots,x_n)$ is weakly/strongly satisfiable 
over $L$ ~iff~ it is so over $\calMO_m$ for $m:=\min\{n,\Card(L)/2-1\}$
with the convention that $\infty/2-1=\infty$.
In particular for infinite  $L$ we have
$\SAT_L=\SAT_{\calMO_\omega}$ and $\sat_L=\sat_{\calMO_\omega}$.
\end{enumerate}
\end{mylemma}
\begin{proof}
Indeed, since all maximal chains in
a modular ortholattice of dimension $2$   are 3-element
 (Fact~\ref{f:mod}), they are of the form $\Zero<a<\One$; and
$a\vee\neg a=\One$ and $a\wedge\neg a=\Zero$
require $\neg a\neq a$:
hence  the set $A$ of atoms of $L$ 
is a disjoint union $A_0 \cup A_1$  where $a \in A_0$ if and 
only if $\neg a \in A_1$. This establishes a).

Concerning b), observe that any homomorphism $\varphi:L\to L'$
commutes with the evaluation of a term: 
$t_{L'}\big(\varphi(a_1),\ldots,\varphi(a_n)\big)=\varphi\big(t_L(a_1,\ldots,a_n)\big)$
for every $a_1,\ldots,a_n\in L$.
Now an assignment $a_1,\ldots,a_n\in L$ 
`lives' in a sublattice $L'$ of $L$
isomorphic to $\calMO_m$ according to a).
And, conversely, $\calMO_m$ maps isomorphically 
onto a sublattice $L'$ of $L$ in case $m\leq\Card(L)/2-1$.
\end{proof}
Every two-dimensional MOL $L$ is thus isomorphic to $\calMO_m$ 
for some cardinal $m$; and
$\Gr(\IF^2)$ is isomorphic to $\calMO_{\Card(\IF)}$:
every one-dimensional subspace of $\IF^2$ corresponds to an atom.
Hence, checking weak/strong satisfiability according to 
Definition~\ref{d:QL}e) naively involves
an infinite choice of possible arguments (=subspaces of $\calH$).
However from Lemma~\ref{l:2D} we conclude
the following extension of Example~\ref{x:EvalSat}:

\begin{myproposition} \lab{p:2DinNP}
\begin{enumerate}
\item[a)]
The following evaluation problem over $\calM_\omega$
is decidable in polynomial time:
\begin{quote}
Given a term $t(X_1,\ldots,X_n)$ 
as well as an assignment $x_1,\ldots,x_n\in\calMO_\omega$ 
and value $y\in\calMO_\omega$
(with $\Zero,\One$ encoded as integers $0,1$ and
atoms $a_m,\neg a_m$ as $2m,2m+1$, say),
is $t_{\calMO_\omega}(x_1,\ldots,x_n)=y$?
\end{quote}
\item[b)]
For any nonempty class $\calC$ of 2-dimensional MOLs,
both $\SAT_\calC$ and $\sat_\calC$ are in $\calNP$.
\end{enumerate}
\end{myproposition}
\begin{proof}
\begin{enumerate}
\item[a)] Disregarding parsing details,
$t$ can be evaluated by recursion on its subterms.
Concerning the recursion bottom
observe that $b\vee c=\One$    
holds in $\calMO_\omega$  for $\Zero\neq b\neq c\neq\Zero$;
and the other cases are trivial anyway.
\item[b)]
Consider a nondeterministic
Turing machine which, on input of an $n$-variate term $t(\bar x)$,
calculates from fixed $\max\{\Card(L):L\in\calC\}$
the $m$ according to Lemma~\ref{l:2D}b), and
then guesses and verifies
an assignment $\bar b$ in $(\calMO_m)^n$ as in a).
\qed\end{enumerate}
\end{proof}
2-dimensional quantum satisfiability is thus 
computationally as hard as 1-dimensional (i.e. Boolean) 
satisfiability:

\begin{theorem} \lab{t:2DNPc}
For any 2-dimensional modular ortholattice $L$,
both $\SAT_L$ and $\sat_L$ are $\calNP$--complete.
\end{theorem}
\begin{proof}
In view of Propositions~\ref{p:NPhard} and \ref{p:2DinNP}
it suffices to show $\calNP$--hardness of $\sat_L$.
For $L=\calMO_m$ for $m\geq2$, this follows from 
Example~\ref{x:2D}b), observing that strong satisfiability $t(\bar x)=\One $ 
here reduces to the weak satisfiability of
$\neg C(y,z)\wedge C\big(t(\bar x),y\big)\wedge C\big(t(\bar
x),z\big)\wedge t(\bar x)$: 
clearly in polynomial time.
For the remaining two-dimensional modular ortholattice
$\calMO_1\cong\{\Zero ,\One \}^2$ on the other hand,
the claim follows from $\sat_{L\times L'}=\sat_L\cup\sat_{L'}$.
\end{proof}

\subsection{Fixed-dimensional Satisfiability over Hilbert Lattices is in \BSS--$\calNP$}
\lab{ss:UpperComplexity}  
It has been observed in \mycite{Corollary~7}{Hagge1}

\begin{fact} \lab{f:Tarski}
Based on Tarski's real (!) quantifier elimination,
$\sat_{\Gr(\IC^d)}$ is decidable, uniformly in $d$.
\end{fact}
Intuitively, a nondeterministic \BSS machine can decide
satisfiability of a term $t$ by guessing an
assignment $\bar X$ in $\Gr(\IF^d)$ and evaluating $t_{\Gr(\IF^d)}(\bar X)$.
However for $\IF=\IC$,
the latter requires separate access to real and imaginary
parts---which a $\IC$-machine does not have (Definition~\ref{d:BSS}).
\cxx{
\begin{digression} \lab{d:ImagField}
Recall that we required all $\IF$ 
to be $\ast$-subfields of $\IC$, i.e.
subfields closed under conjugation
(since conjugation enters into the concept
of scalar product).
\begin{longenum}
\item[a)]   
A mere subfield of $\IC$, as $\IQ(\sqrt[3]{2}e^{2\pi/3})$,
may fail to be $\ast$-subfield.
\item[b)]
Observe that $\IF \not\subseteq \IR$ is a $\ast$-subfield
of $\IC$ if and only if 
 $\IF= \IF' +i \IF'$
with subfield $\IF'=\IF\cap \IR$ of $\IR$; 
equivalently: with $\IF'=\Re(\IF):=\{x\mid x,y\in\IR,x+iy\in\IF\}$. 
\item[c)]
Strictly speaking, the concept of an $\ast$-field 
defines the imaginary unit $i$ and/or the imaginary part
only up to the isomorphism $i\mapsto -i$.
\end{longenum}
\end{digression} }
Furthermore note that every $U \in \Gr(\IF^d)$ 
is of the form $U=\range A$ for (many) 
$A \in \IF^{d \times d}$
where $\range A$ 
denotes the linear subspace of $\IF^d$ spanned by
the columns of $A$.

\begin{myproposition} \lab{p:UpperComplexity}
Fix $\IF$.
\begin{enumerate}
\item[a)]  
Given $d\in\IN$ and matrices $A,B\in\IF^{d\times d}$
\begin{enumerate}
\item[i)] a matrix $C\in\IF^{d\times d}$ with $\range(C)=\range(A)\vee\range(B)$ and
\item[ii)] a matrix $C'\in\IF^{d\times d}$ with $\range(C')=\range(A)\wedge\range(B)$
\end{enumerate}
can be calculated by a constant-free \BSS-machine 
over $\IF$ in time $\calO(d^3)$. Similarly,
\begin{enumerate}
\item[iii)] a matrix $C''\in\IF^{d\times d}$ with $\range(C'')=\neg\range(A)$
\end{enumerate}
can be calculated by a constant-free \BSS-machine 
$\calM$ over $\IF$ in time $\calO(d^3)$:
in case $\IF\subseteq\IR$; otherwise $\calM$ 
in general needs both real part  $\Re(A)$ and imaginary part $\Im(A)$.
\item[b)]
First suppose $\IF\subseteq\IR$.
Then, given an $n$-variate term $t$ 
and matrices $A_1,\ldots,A_n\in\IF^{d\times d}$
representing $U_j=\range(A_j)\in\Gr(\IF^d)$,
a constant-free \BSS-machine $\calM$ over $\IF$ can calculate
$C\in\IF^{d\times d}$ with
$\range(C)=t_{\Gr(\IF^d)}(U_1,\ldots,U_n)$
in time polynomial in both $d$ and $|t|$.
\\ 
In the general case $\IF\subseteq\IC$,
given matrices $\Re(A_1),\Im(A_1),\ldots,\Re(A_n),\Im(A_n)\in\Re(\IF)^{d\times d}$,
representing $U_j=\range(A_j)\in\Gr(\IF^d)$,
a similar machine over $\Re(\IF)$ can calculate
$\Re(C),\Im(C)\in\Re(\IF)^{d\times d}$ with
$\range(C)=t_{\Gr(\IF^d)}(U_1,\ldots,U_n)$.  
\item[c)]
Both weak and strong satisfiability over $\Gr(\IF^d)$
of a given term $t$ can be decided by a 
nondeterministic constant-free \BSS-machine 
over $\Re\IF$ in time polynomial in $|t|$ and in $d$.
In particular, it holds
\[ \sat_{\Gr(\IR^d)},\,\SAT_{\Gr(\IR^d)},\:
\sat_{\Gr(\IA^d)},\,\SAT_{\Gr(\IA^d)},\:
\sat_{\Gr(\IC^d)},\,\SAT_{\Gr(\IC^d)}\:
\;\in\:\calBP(\calNP_\IR^0)\:\subseteq\:\PSPACE \enspace . \]
\item[d)]
Concerning satisfiability \emph{with constants}, 
we make Definition~\ref{d:QL}h) more precise:
\begin{multline*}
\SAT_{\Gr(\IF^d),\Gr(\IF^d)}\;:=\;
\Big\{\big\langle t(x_1,\ldots,x_n,y_1,\ldots,y_m),C_1,\ldots,C_m\big\rangle\;\Big|\;
C_1,\ldots,C_m\in\IF^{d\times d},
\\[-0.9ex]
\exists A_1,\ldots,A_n\in\IF^{d\times d}:
t_{\Gr(\IF^d)}(\range A_1,\ldots,\range A_n,
\range C_1,\ldots,\range C_m)=\One\Big\}
\end{multline*}
where matrices $C_j$ are encoded as $d^2$-element 
sequences over $\IF$. Then it holds
\[ \sat_{\Gr(\IF^d),\Gr(\IF^d)},\;
\SAT_{\Gr(\IF^d),\Gr(\IF^d)},\;
\sat_{\Gr(\Re(\IF)^d),\Gr(\Re(\IF)^d)},
\;\SAT_{\Gr(\Re(\IF)^d),\Gr(\Re(\IF)^d)}
\;\in\calNP^0_{\Re(\IF)}. \]
\end{enumerate}
\end{myproposition}
Claims~b) and c) can be regarded as natural a generalization of
Proposition~\ref{p:2DinNP}a) and b) to fixed, higher dimensions.

\begin{proof}
\begin{longenum}
\item[a\,i)] 
Gaussian Elimination (Example~\ref{x:BSS}) 
works for calculation with columns, as well, and allows 
to turn a matrix into column echelon form.
Recall, that Gaussian Elimination 
leaves the linear span of the columns invariant.
Given $A$ and $B$ apply this 
to the  $(d\times 2d)$--matrix $(A|B)$ 
and put the first $d$ columns
(the remaining ones are zero columns)  of the
echelon form into $C$ to obtain
 $\range(A)\vee\range(B)=\range(C)$.
\COMMENTED{Gaussian Elimination (Example~\ref{x:BSS}) turns
the $(d\times2d)$--matrix $(A|B)$ 
with $\range(A)\vee\range(B)=\range(A|B)$
into reduced 
row echelon form;
so that 
 $\dim\range(A|B)\leq d$ spanning columns can simply be read off
$(A|B)$ 
and put into $C$. }  
\item[a\,ii)] To compute the meet,
we use an algorithm due to \person{Zassenhaus} \cite{Zassenhaus}. 
Form the block $2d \times 2d$--matrix 
 $\binom{A \;\; B}{A \;\; O}$ of rank 
$r=\rank (A) +\rank(B)$ (which can be seen using row transformations).   
Use Gaussian Elimination
on columns (according to the upper half of the matrix)
 to transform  this matrix  into $\binom{A' \;\; O}{A'' \; C'}$
with the same block structure
and rank $r=\rank(A')+ \rank(C')$ since $\rank(A'')=\rank(A) \leq \rank(A')$.
Then $\range(C')= \range(A) \wedge \range(B)$.
Indeed, as observed in i), the columns of $A'$ span
$\range(A) \vee \range(B)$. It follows by the
dimension formula  that
$\rank(C') = \dim(\range(A) \wedge \range(B))$
and it remains to show that the columns
$\vec c$  of $C'$ are in $\range(A) \wedge \range(B)$.
But, by construction, the
column $\binom{\vec 0}{\vec c}$ 
is  of the form  $\binom{\vec a}{\vec a} +\binom{\vec b}{\vec 0}$
where the first summand is a linear  combination
of columns of $\binom{A}{A}$, the second of $\binom{B}{O}$.
From $\vec a +\vec b=\vec 0 $ it follows
$\vec c= \vec a \in   \range(A) \wedge \range(B)$.
\item[a\,iii)]
It is well known that 
$\ker(A^\ast)=\neg(\range A)$ where $A^\ast$ is
the \textsf{adjoint} of $A$, i,.e. 
transposed and complex conjugate; indeed, $ A \perp \vec v$
(i.e. $\vec a \perp \vec v$ for all columns $\vec a$ of $A$) if
and only if
 $A^\ast\vec v =0$ i.e. $\vec v \in  \ker(A^\ast)$.
Hence let machine $\calM$ calculate $A^\ast$
 (using
separate access to both $\Re(A)$ and $\Im(A)$);
then apply Gaussian Elimination to obtain 
a basis of its kernel. \\
We now show that calculating to
given $A\in\IC^{2\times 2}$ some 
$C''\in\IC^{2\times 2}$ with $\range(C'')=\neg\range(A)$
entails complex conjugation.
To this end consider $A:=\binom{0\;\;z}{0\;\;1}$.
Then $\neg\range(A)=\IC\binom{-1}{\bar z}$ and therefore
necessarily $C''=\binom{-a\;\;-b}{a\bar z\;\;b\bar z}$
for some $a,b\in\IC$ not both 0. A \BSS machine over $\IC$ 
can identify a non-zero column of $C''$; 
and divide its lower by the upper component to obtain $-\bar z$,
from which both $\Re(z)=(z+\bar z)/2$ and $\Im(z)=(z-\bar z)/2$ are readily obtained.
\item[b)] 
Neglecting parsing details as in Proposition~\ref{p:2DinNP},
recursive application of a) to each subterm of $t$ yields the claim.
\item[c)]
In the real case,
let the nondeterministic \BSS-machine `guess'
$A_1,\ldots,A_n\in\IF^{d\times d}$;
then apply b) to evaluate $t$ thereon and accept iff
the result is the identity ($\SAT$) or is not the zero matrix ($\sat$), respectively.
In the complex case, similarly `guess'
$\Re(A_1),\Im(A_1),\ldots,\Re(A_n),\Im(A_n)\in\Re(\IF)^{d\times d}$.
The final claim follows from Fact~\ref{f:BSSBP}b)+c). 
\item[d)]
Given $(n+m)$-variate $t$ as well as
$C_1,\ldots,C_m\in\IF^{d\times d}$,
guess $A_1,\ldots,A_n\in\IF^{d\times d}$
and then proceed as before to evaluate
$t_{\Gr(\IF^d)}(\range A_1,\ldots,\range A_n,\range C_1,\ldots,\range C_m)$.
\qed\end{longenum}\end{proof}
Observe that the calculations in Proposition~\ref{p:UpperComplexity}a) 
can be performed by a BSS-machine over the ring $\IF$, i.e. not using inversion.
An alternative proof of Proposition~\ref{p:UpperComplexity}c)
replaces Gaussian Elimination in a)  
by a non-deterministic approach based on
Fact~\ref{f:lat} below.

\subsection{Abstract Unitary Spaces and Bases} 
The results of the preceding subsection  have been most conveniently formulated
for the $\IF$-unitary spaces  $\IF^d$  and their ortholattices.
Transition to abstract spaces $\calH$  is possible
via  bases. A  basis
$\vec v_1, \ldots , \vec v_d$ 
is \textsf{orthogonal} if 
$\langle \vec v_i \mid \vec v_j \rangle =0$ 
for all $i \neq j$.
It is \textsf{orthonormal} if, in addition,
$\langle \vec v_i \mid \vec v_i \rangle =1$
for all $i$. 
An  \textsf{isometry} $\omega:\calH\to\calH'$
between $\IF$-unitary spaces $\calH$ and $\calH'$
is an $\IF$-linear isomorphism
preserving scalar products:
$\langle\omega(\vec x)|\omega(\vec y)\rangle=\langle\vec x|\vec y\rangle$
for all $\vec x,\vec y\in\calH$.

\begin{convention} 
Unless stated otherwise, $\calH$ is assumed to be a finite
dimensional  $\IF$-unitary space.
$\IF$ is \textsf{Pythagorean}  if
 in $\IF\cap \IR$ 
any   sum of squares is a square.
\end{convention}

\begin{fact} \lab{f:GramSchmidt}
\begin{enumerate}
\item[a)]
Any  $\IF$-unitary space
admits  an orthogonal basis.
\item[b)]
The basis can be chosen orthonormal if
$\IF$ is Pythagorean. 
\item[c)]
An   $\IF$-unitary space
admitting  an orthonormal
basis of $d$ vectors
is isometric to   $\IF^d$ (with the basis mapped onto
the canonical one).
\item[d)]
If $\calH$ and $\calH'$ are isometric,
then $\Gr(\calH)$ and $\Gr(\calH')$ are isomorphic.
\end{enumerate}
\end{fact} 
\begin{proof} 
\begin{longenum}
\item[a)] 
follows from the   \emph{Gram-Schmidt} Process, 
cmp. e.g. \mycite{\S2-3}{Gelfand},
  \mycite{\S73-74}{Halmos}, or \mycite{6.19}{Axler}.  
\item[b)] 
Having all $(a+bi)(a-bi)=a^2+b^2$  squares in $\IF \cap \IR$,
any $\langle \vec v \mid \vec v \rangle$ is a square
in $\IF \cap \IR$, too,  so that one can normalize.
\item[c)]
Matching two orthonormal bases extends to an isometry.
\item[d)]
A linear homomorphism $\omega:\calH\to\calH'$ satisfies
$\omega[A]+\omega[B]=\omega[A+B]$ for all $A,B\in\Gr(\calH)$; 
an injective one also $\omega[A]\cap\omega[B]=\omega[A\cap B]$;
and an isometry even $\neg\omega[A]=\omega[\neg A]$.
In the latter case,
$\Gr(\calH)\ni A\mapsto\omega[A]\in\Gr(\calH')$ thus
constitutes an isomorphism.
\qed\end{longenum}\end{proof}
Fact~\ref{f:GramSchmidt}d) and its converse can be traced back to
\cite{Birkhoff}, its relevance for Quantum Logic   is also pointed out
in \cite[line after \textsc{Definition2}]{Hagge2}.
We weaken Item~c) in the following

\begin{mydefinition} \lab{d:Equinormal}
A family
$\vec v_1, \ldots ,\vec v_d$ of vectors is 
\textsf{equinormal} if
$\|\vec v_i\|^2=\|\vec v_j\|^2$ for all $i,j$. 
\end{mydefinition} 
For instance,
the one-dimensional $\IQ$-unitary vector space 
$\{(x,x):x\in\IQ\}$ 
has the single vector $(1,1)$ as an (trivially equinormal and orthogonal) 
basis but contains no unit vector.

\cx{
\begin{digression}
The two-dimensional $\IQ$-unitary vector space 
$\{(x,x,y):x,y\in\IQ\}$ has vectors 
$(1,1,1)$ and $(1,1,-1)$ constituting
an equinormal basis but does not admit an
orthogonal equinormal one: 
Equations
$2x^2+y^2=2x'^2+y'^2$ and $2xx'+yy'=0$ have
(other than the trivial $x=y=x'=y'=0$)
no simultaneous solution over $\IZ$ nor, dividing
by a common denominator, over $\IQ$.
\end{digression}}
\begin{mylemma} \lab{l:Equinormal}
Let $\calH$ and $\calH'$ denote finite-dimensional
$\IF$-unitary spaces.
\begin{enumerate}
\item[a)]
If both $\calH$ and $\calH'$ admit equinormal orthogonal bases
and have coinciding dimensions,
then $\Gr(\calH)$ and $\Gr(\calH')$ are isomorphic.
\item[b)]
If both $\calH$ and $\calH'$ admit equinormal orthogonal bases
and $\dim(\calH)$ divides $\dim(\calH')$,
then $\Gr(\calH)$ embeds into $\Gr(\calH')$.
\end{enumerate}
\end{mylemma}
Partial converses can be found in Corollary~\ref{c:Equinormal} below.
\begin{proof}
\begin{longenum}
\item[a)]
Let $\vec u_1,\ldots,\vec u_d$ be the
equinormal orthogonal base for $\calH$
and $\vec v_1,\ldots,\vec v_d$ that for $\calH'$
with isomorphism $\omega:\calH\to\calH'$,
$\vec u_j\mapsto\vec v_j$. In view of the proof of
Fact~\ref{f:GramSchmidt}d) it remains to show that 
$\vec x\perp\vec y\Leftrightarrow\omega(\vec x)\perp\omega(\vec y)$ holds.
Indeed, exploiting orthogonality and equinormality,
\[
\Big\langle \sum\nolimits_k \alpha_k\vec u_k,\sum\nolimits_\ell \beta_\ell \vec u_\ell\Big\rangle
\;=\;
\sum\nolimits_k \langle \alpha_k\vec u_k,\beta_k \vec u_k\rangle
\;=\; \|\vec u\|^2 \cdot\sum\nolimits_k \alpha_k\bar\beta_k \]
vanishes iff 
$\|\vec v\|^2 \cdot\sum\nolimits_k \alpha_k\bar\beta_k
=\big\langle \omega\big(\sum\nolimits_k \alpha_k\vec u_k\big),
\omega\big(\sum\nolimits_\ell \beta_\ell \vec v_\ell\big)\big\rangle$
does.
\item[b)]
Let $d=\dim(\calH)$ and $\dim(\calH')=kd$, $k\in\IN$.
Let $(\vec v_{j\ell})_{1\leq j\leq d,1\leq\ell\leq k}$ denote
an orthogonal equinormal basis of $\calH'$.
Then $\calH'=V_1\obot\ldots\obot V_k$ where
$V_j:=\IF\vec v_{j1}+\cdots+\IF\vec v_{jd}$.
Now by a), there is an isomorphism 
 $\omega_\ell:\Gr(\calH) \rightarrow\Gr(V_\ell)$
for each  $1\leq\ell\leq k$. On the other hand, 
from  Observation~\ref{o:Product}a)
it follows by induction that  $\Gr(V_1)\times\Gr(V_2)\times\cdots\times\Gr(V_k)$
embeds into $\Gr(\calH')$.
Thus, $\omega(X)= \big(\omega_1(X), \ldots, \omega_k(X)\big)$ 
 defines an embedding of $\Gr(\calH)$ into  $\Gr(\calH')$.
\qed\end{longenum}\end{proof}
Since the standard basis for $\IF^d$ is clearly equinormal orthogonal,
we conclude from Proposition~\ref{p:UpperComplexity}c) and Lemma~\ref{l:Equinormal}a)

\begin{corollary}
Suppose $\calH$ is an $\IF$-unitary  vector space
admitting an equinormal orthogonal basis.
Then $\sat_{\Gr(\calH)},\SAT_{\Gr(\calH)}\in\,\calBP(\calNP_{\Re(\IF)}^0)$.
\end{corollary} 

\section{Strong Satisfiability is Complete for \BSS-$\calNP$
in Dimensions $\geq3$}\lab{s:Arith2QL} 

The main result of this section
shows that, for fixed $\calH$ of dimension $d \geq 3$,
 strong satisfiability in
$\Gr(\calH)$   is hard for the complexity class
$\calBP(\calNP_{\Re(\IF)}^0)$.  This is shown by
recapturing $\IF$ within $\Gr(\calH)$ 
using coordinatization methods for modular
lattices due to \person{von Neumann} \cite{Neumann2}.
Subsections~\ref{ss:Affine} and \ref{ss:Projective}
convey the elementary (affine and projective, respectively) 
geometric intuition behind the quantum logic terms 
then presented in Subsection~\ref{ss:vonNeumann}.
 
\subsection{Background from Elementary Affine Geometry} \lab{ss:Affine}
Let us start with Elementary Geometry  in the
(affine)  plane. 
 Following \person{Hilbert} \cite{Hilbert}, only the
structure matters. The basic objects are  points and lines, the basic
relation is incidence: a given  point  is or is not on a given 
line. A line may be 
considered as the set of all points incident with it.
Following Descartes, to introduce coordinates,
we need two lines $\ell_2,\ell_3$  (the choice for these
indices will become clear, later) 
intersecting
in a point $A_1$. Then for any point $P$ and $i\neq j$ in
$\{2,3\}$ 
  we have
the intersection $P_i$ of $\ell_i$ with
the parallel  to  $\ell_j$  through $P$. Thus, 
one obtains a bijective correspondence between points $P$ 
and their coordinate pairs $(P_2,P_3)$. 

In the sequel, let us consider points
$P$ on the coordinate line $\ell_2$ 
as \textsf{scalars}  and do calculations with them. 

\begin{observation} \lab{o:sub} 
Subtraction and multiplication of scalars
can be defined geometrically:
\end{observation} 
Choosing $A_{12}\neq A_1$ on $\ell_2$ 
we fix an orientation on $\ell_2$ 
which allows to understand $P$ 
as the length of the segment $A_1P$ 
with  sign given by   orientation.
Then the difference $P\ominus Q$
can be obtained by the geometric construction
of Figure~\ref{f:mult}a): choose $A_{13}\neq A_1$ on
$\ell_3$ to determine the line $\ell$ through
$A_{12}$ and $A_{13}$,
 let $Q_{13}$ the intersection with $\ell_3$
of the parallel to $\ell$  through $Q$,      
let $S$ the intersection of $\ell$ with the
parallel to $\ell_2$ through $Q_{13}$ and
$P \ominus Q$ the intersection of $\ell_2$
and  the parallel to $\ell_3$ through $S$.

\cx{
\begin{figure}[htb]
\includegraphics[width=0.45\textwidth]{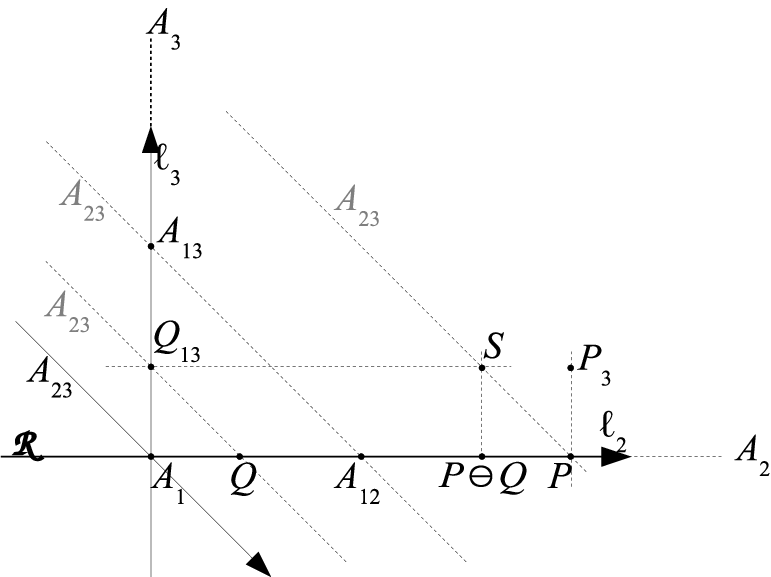}
\hfill\includegraphics[width=0.52\textwidth]{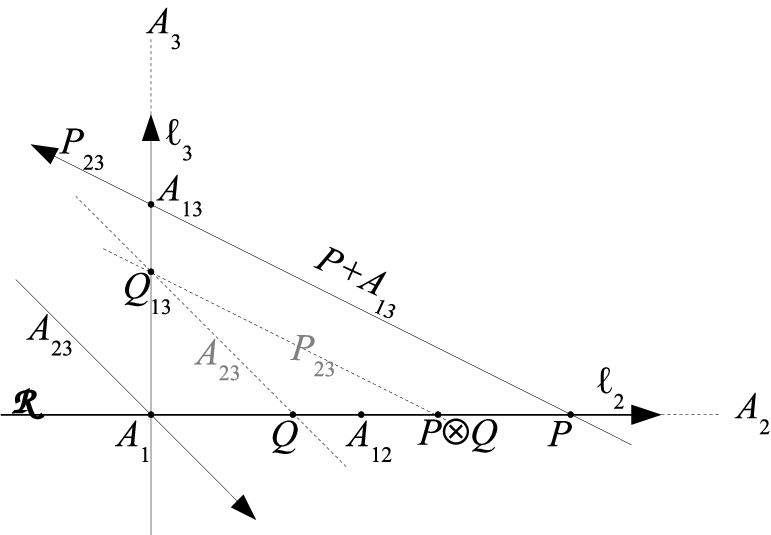}

\caption{\label{f:mult}Geometry of ~ a) subtraction and ~ b) multiplication}
\end{figure} }

Multiplication is based on the \emph{Intercept Theorem}.
Having  determined $Q_{13}$  as before, let $k$ the line through
$P$ and $A_{13}$ and $P\otimes Q$ the intersection 
of $\ell_2$ and the parallel to $k$ through $Q_{13}$.
Observe, that $A_{12}\otimes P=P$.

 Notice, that 
it is not required that $\ell_3$ be perpendicular
to $\ell_2$ or that  the distances $A_1A_{12}$ and $A_1A_{13}$ are
 the same.  Operations  with scalars on $\ell_3$ 
can be defined, symmetrically, and one can show 
that 
the correspondence $\omega$  between scalars on  $\ell_2$ and $\ell_3$,
respectively,   given by $\ell$ is indeed an isomorphism.
Using  Desargues' Theorem resp. axiom 
one can further show that  the scalars form a division ring
- which is rendered as $\IR$ if we   have all 
axioms for incidence and orientation in  Elementary Geometry available.

\begin{observation} \lab{o:coorsy} 
A \textsf{coordinate system} of the
plane is given by the \textsf{origin} $A_1$ and
 the \textsf{coordinate lines} $\ell_2,\ell_3$
(resp. their directions $A_2,A_3$)
and \textsf{unit points} 
$A_{1i}\neq A_i$ on $\ell_i$, $i=2,3$.   
\end{observation} 
Suppose that $\ell_2$ and $\ell_3$ are
perpendicular.
Consider the  lines $g$ and $h$ through
$A_1$ and the points with coordinates
$(A_{12},A_{13})$ and
 $(\ominus A_{12}, A_{13})$, respectively.
Thus,  $h$ is parallel to $\ell$.

\begin{observation} \lab{o:orthogframe}     
$\ell$ is 
 perpendicular to $g$   if and only if 
$\omega$ is an isometry,
i.e. if the 
``unit lengths'' $A_1A_{12}$ and $A_1A_{13}$  agree,
see Figure~\ref{f:orthogframe}. 
\end{observation} 
\cx{
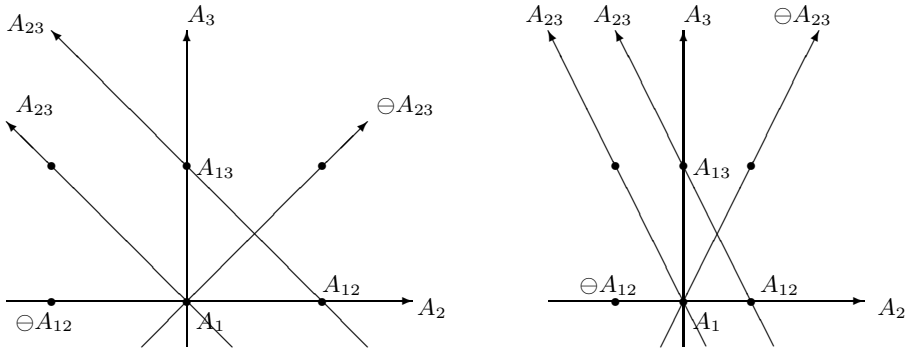
\begin{figure}[htb]
\hspace*{5ex}
\setlength{\unitlength}{6mm} 
\begin{picture}(13,10)  
\put(0,2){\vector(1,0){9}}  
\put(4,1){\vector(0,1){7}}  
\put(4,2){\circle*{0.2}} 
\put(1,2){\circle*{0.2}} 
\put(1,5){\circle*{0.2}} 
\put(4,5){\circle*{0.2}} 
\put(9.1,1.7){$A_2$} 
\put(7,2){\circle*{0.2}} 
\put(7,5){\circle*{0.2}} 
\put(7,2.2){$A_{12}$} 
\put(4.2,4.8){$A_{13}$}  
\put(8,1){\vector(-1,1){7}} 
\put(4.2,1.4){$A_1$} 
\put(4,8.2){$A_3$}
\put(0,8){$A_{23}$}
\put(0.2,6.2){$A_{23}$}
\put(0.2,1.4){$\ominus A_{12}$}
\put(3,1){\vector(1,1){5}} 
\put(5,1){\vector(-1,1){5}} 
\put(8.2,6.2){$\ominus A_{23}$}

\put(12,2){\vector(1,0){7}}  
\put(15,1){\vector(0,1){7}}  
\put(13.5,2){\circle*{0.2}}
\put(13.5,5){\circle*{0.2}}
\put(15,2){\circle*{0.2}} 
\put(19.3,1.7){$A_2$} 
\put(16.5,2){\circle*{0.2}} 
\put(16.5,5){\circle*{0.2}} 
\put(14.5,1){\vector(1,2){3.5}}
\put(15,5){\circle*{0.2}} 
\put(16.7,2.2){$A_{12}$} 
\put(12.7,2.2){$\ominus A_{12}$} 
\put(15.2,4.8){$A_{13}$}  
\put(17,1){\vector(-1,2){3.5}} 
\put(15.5,1){\vector(-1,2){3.5}} 
\put(15.2,1.4){$A_1$} 
\put(15,8.2){$A_3$}
\put(13,8.2){$A_{23}$}
\put(11.5,8.2){$A_{23}$}
\put(17,8.2){$\ominus A_{23}$} 

\end{picture}
\caption{\label{f:orthogframe} Orthogonal frames: Orthonormal versus
  non-orthonormal} 
\end{figure}}

\subsection{Background from Projective Geometry} \lab{ss:Projective}

 Projective Geometry    originated from an
analysis of perspective drawing.  The basic idea
(in the plane) 
is to add  
 a \emph{point at infinity} for any parallel  pencil
of lines  (and make it  incident with  each line of the pencil)
and to put all these points on a new \emph{line}  $\ell_\infty$ 
\emph{at infinity} --- just imagine several straight railroads
in the  plains. Thus, the following applies:

\begin{mydefinition} \lab{d:pplane} 
A \textsf{projective plane} $\calP$
is an incidence structure  such that   any two distinct lines
meet in exactly one point and any  two
distinct points $P,Q$ determine a unique line $P\vee Q$
incident with both $P,Q$.
\end{mydefinition} 
The original affine  plane
$\calA$ is  recovered  within  this projective plane
$\calP$ as $\calP \setminus \ell_\infty$ --- 
its lines are given by those lines of $\calP$
which are incident with at least two points from $\cal A$.  
We still may consider lines just as sets of points.  
But now there are
only two basic operations in $\calP$
 and these may also be used
 to recapture $\calA$, completely:
 drawing the parallel $g$ 
to $h$ through $P$ is  replaced by 
 determining the intersection $S$ of $h$ with $\ell_\infty$ 
and to join $S$ with $P$ to obtain $g$.
Notice, that choosing any line $\ell$ 
 of $\calP$ we obtain an
affine plane $\calP \setminus \ell$ for which $\ell$ 
is the line at infinity.

In particular, if a line $\ell_\infty$  at infinity is designated, 
\emph{directions} are just points on $\ell_\infty$.
Thus, the concept of a coordinate system 
in Observation~\ref{o:coorsy}  now translates into the following
(cmp. Figure~\ref{f:projective}b).

  \cx{
\begin{figure}[htb]
\centerline{\includegraphics[width=0.6\textwidth]{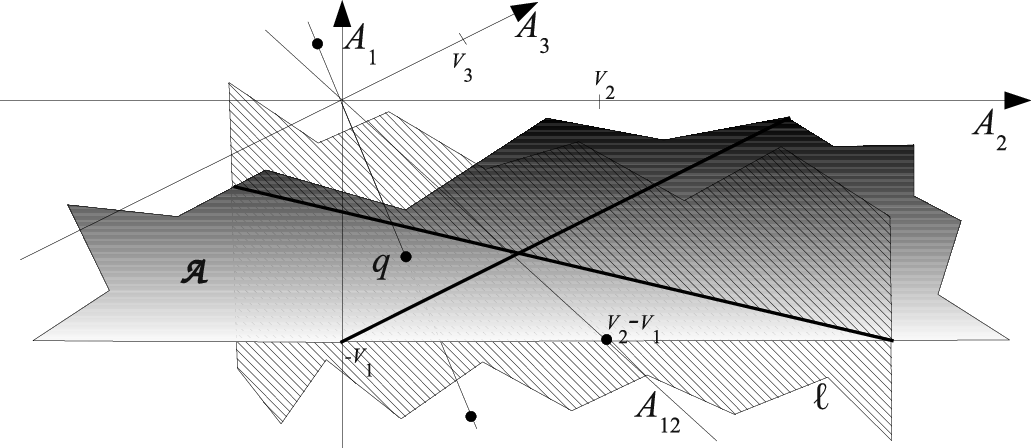}\quad
\includegraphics[width=0.28\textwidth]{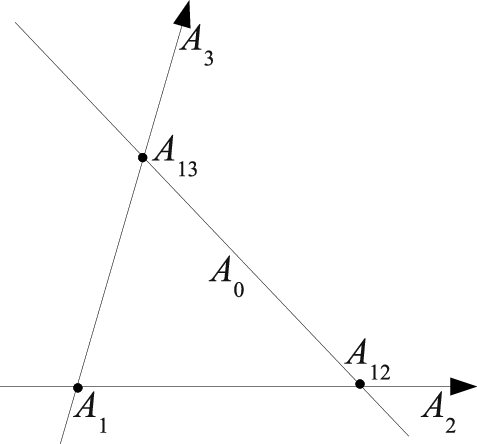}}
\caption{\label{f:projective}a) Projective plane and affine space:
orthogonal vectors $\vec v_i$ induce projective points $A_i$;
the 2-dimensional affine space $\calA$ with an affine line
$\ell\cap\calA$, $\ell$ a line of projective space, i.e. a
homogeneous plane.~ b) Non-orthogonal coordinate system.}
\end{figure}}

\begin{observation} \lab{o:coopro} 
\begin{longenum}
\item[a)]
A coordinate system is given by three non-collinear
\textsf{basis points}
$A_1,A_2,A_3$ and  \textsf{unit points} $A_{1i}$ on $A_1 \vee A_i$ for
$i=2,3$ such that $A_{1i} \neq A_1,A_i$.  Alternatively,  we may choose
 units points $A_{ij}=A_{ji}$ on $A_i \vee A_j$ for
$i \leq j$ such that $A_{12},A_{13},A_{23}$  are  on the same line $\ell$ 
 and $A_{ij} \neq A_i,A_j$. 
The latter choice can be also implemented by choosing $\ell$.     
\item[b)]
How does this relate to the lattices
$\Gr(\calH)$?  We first observe, that due
to the Dimension Formula
  $\Gr_1(\calH)$  and $\Gr_2(\calH)$  
as the sets of  `points' and `lines' 
with incidence relation $\subseteq$ 
form a projective plane
if $\dim(\calH)=3$.
More generally,   $\Gr_1(\calH)$  and $\Gr_2(\calH)$
form a `projective space' $\calP$ of dimension $\dim(\calH)-1$,
where the condition that any two lines meet
 in a some point is  required only for lines
in a `subspace'  which is a plane.
\item[c)]
Of course, with $\calH$ we may also
associate its \emph{affine space}
 with point set $\calH$, lines and planes of
 the form $\vec v+U$ where $U \in \Gr(\calH)$ and
 $\dim(U)=1,2$, respectively.
Choosing \emph{origin} $\vec 0$,  
the  $\vec 0+U$  are  (homogeneous) line and
plane through this origin in the affine space 
$\calH$,  points   and lines
of the projective space $\calP$, respectively. 

With respect to a basis $\vec v_1, \ldots ,\vec v_d$,
the points of $\calP$ are  
 given by their homogeneous coordinates
$\bar x=(x_1, \ldots ,x_d)$,  where $\bar x$ and 
$\lambda \bar x$ ($\lambda \neq 0$ in $\IF$) 
describe the same point. Then the points with
coordinate $x_1=-1$  (the reason why we do not choose
$x_1=1$ will become transparent, later) 
 form a $(d-1)$-dimensional
affine subspace $\calA$  with  lines $U \cap \calA$,
where $U$ is a 2-dimensional linear subspace of $\calH$ 
not parallel to $\calA$.
(Thus, in case $\IF^3$ we may imagine  $\calA$ as  the plane parallel
to the $x_2$--$x_3$--plane and meeting the $x_1$--axis at $-1$;
cf. Figure~\ref{f:projective}a).
This space is canonically
  coordinatized by $\IF^{d-1}$. 
The points with coordinate $x_1=0$ 
form then  a  $(d-2)$--dimensional projective subspace, 
of $\calP$, the
 `hyperplane at infinity''.
\end{longenum}
\end{observation} 

\subsection{Von Neumann Frames for Coordinatizing Lattices} 
\lab{ss:vonNeumann}

Frames were introduced by \person{von Neumann} to recapture
linear algebraic structure within certain modular lattices,
including all Hilbert lattices $\Gr(\calH)$. The intuition came
from Projective Geometry
but the concept can be explained within Linear Algebra, too.
A basis  
 $\vec v_1, \ldots ,\vec v_d$ 
of $\calH$ is turned into certain elements of $\Gr(\calH)$
satisfying relations which, conversely,
guarantee that they have been induced by a basis---sufficiently 
unique to  construct an isomorphism of $\calH$  onto $\IF^d$ 
via  an isomorphism from $\Gr(\calH)$ onto  $\Gr(\IF^d)$.

From the 
$A_i=\IF \vec v_i$  we may recover the $\vec v_i$ up to
scalar multiple. Moreover, relations 
\[ A_1 \vee \ldots \vee A_d= \One
\quad\text{ and }\quad
A_i \wedge \bigvee\nolimits_{j\neq i}  A_j =\Zero \
\quad (1\leq i\leq d)\]
make sure that $\calH$ is a direct sum 
$\calH= A_1 \oplus \ldots \oplus A_d$.

In particular
 $\dim(A_i)=1$ since $\dim(\calH)=d$; and one recovers a basis 
of $\calH$ by choosing non-zero $\vec v_i \in A_i$ ($1\leq i\leq d$).
What has to be added is a simultaneous scaling  of the $\vec v_i$ 
in terms of the lattice.  
The basis vectors $\vec v_i$ are naturally related via the
isomorphisms
$\eta_{ij}: A_i \rightarrow A_j$ such that
$\eta_{ij}(\vec v_i) =\vec v_j$. Clearly,
$\eta_{jk} \circ \eta_{ij} =\eta_{ik}$.

To capture these isomorphisms within the lattice,
one might associate with  each $\eta_{ij}$ 
its graph $\IF(\vec v_i +\vec v_j)=\{\vec x +\eta_{ij}(\vec x)\mid
\vec x \in A_i\}$.
Though, as e.g. the graphs of $\eta_{12}, \eta_{13}$ and
$\eta_{23}$ give a direct decomposition of $\IR^3$, 
there is no meaningful
lattice relation on these graphs.    
 \person{Von Neumann}  realized, that there are meaningful
relations if one uses  \emph{negative graphs} 
$A_{ij}= A_{ji}= \IF(\vec v_i -\vec v_j)$.
 The relevant relations are 
\[ A_{ij} \wedge A_i =\Zero,\quad A_{ij} \vee A_i =\One,\quad
A_{ik}=(A_i\vee  A_k)\wedge ( A_{ij}\vee A_{jk})
\qquad (i\neq j\neq k\neq i) \]     
(cmp. Observation~\ref{o:coopro}a).  
The first two express the fact that
$A_{ij}$ defines an isomorphism from $A_i$ onto $A_j$,
the third guarantees the above  composition laws. 
Observe that the linear automorphisms  leaving the
$A_i$ and $A_{ij}$ invariant   are exactly the
multiplications with a scalar.
 Following  \cite{Neumann2}
this leads to the 

\begin{mydefinition} \lab{d:frame} 
A system
$\bar a= (a_{ij}\mid 1 \leq  i,j \leq d)$
of elements of a modular lattice,  $L$,
is  a  $d$-\textsf{frame}
if it satisfies  (abbreviating $a_i=:a_{ii}$)
\[\begin{array}{lcll}  
\One&=& a_1\vee \ldots \vee a_d,\\
\Zero&=& a_i \wedge \bigvee_{j \neq i}  a_j&\mbox{ for all } i,j =1, \ldots ,d,\,i\neq j,   \\
\Zero&=&a_i \wedge a_{ij} &\mbox{ for all } i,j =1, \ldots ,d,\,i\neq j,\\
a_i \vee a_j &=& a_i \vee a_{ij} &\mbox{ for all } i,j =1, \ldots ,d,\,i\neq j,\\
a_{ij}&=&a_{ji}  &\mbox{ for all } i,j =1, \ldots ,d,\,i\neq j,\\
a_{ik}&=& (a_i\vee a_k) \wedge (a_{ij}\vee a_{jk}) 
& \mbox{ for all } i,j,k =1, \ldots ,d,\,i\neq j\neq k \neq i.   
 \end{array}\]
\end{mydefinition} 
The frame relations and Fact~\ref{f:mod} imply 
the following.

\begin{fact} \lab{f:frame}   Let $\bar a$ be a $d$-frame
in the modular lattice $L$.
\begin{enumerate}
\item[a)] 
$\dim a_i =\dim a_j =\dim a_{ij}$
for all  $i,j>0$. 
\item[b)] $\dim \bigvee_{i \in I}a_i =|I|\cdot \dim a_1$  
for any $I \subseteq \{1, \ldots ,d\}$.
\item[c)] $\dim(L)=d\cdot\dim(a_1)$.
\item[d)] 
 Let $\dim(\calH)=d$. Then
$\bar A=(A_{ij}\mid 1\leq i,j \leq d)$ is a $d$-frame
of $\Gr(\calH)$ if and only if there exists a basis
$\vec v_1, \ldots ,\vec v_d$ of $\calH$ such that 
$A_i=A_{ii}= \IF\vec v_i$ and $A_{ij}=\IF(\vec v_i -\vec v_j)$
for $i \neq j$.
\item[e)] In particular, every $d$-dimensional $\Gr(\calH)$ admits a $d$-frame. 
\end{enumerate}
\end{fact}
\begin{proof} a) to c) are immediate by Fact~\ref{f:mod}.
Concerning d), assume that $\bar A =(A_{ij}\mid 1\leq i,j \leq
 d)$ is a  $d$-frame
in  $L=\Gr(\calH)$.
Since  $\calH=A_1 \oplus \ldots \oplus A_d$ 
with $\dim(A_i)=1$, we may choose a basis 
$\vec v_1, \ldots, \vec v_d$ of $\calH$ 
such that $A_i=\IF \vec v_i$. Since 
$\dim A_{1i}=1$  for $i>1$,  we may choose
$\vec w_{1i}$ such that $A_{1i}= \IF \vec w_{1i}$.
From $A_{1i} \subseteq A_{1} + A_i$ 
we get $\vec w_{1i}=r_i\vec v_1 +s_i \vec v_i$
and $r_i \neq 0$ since $A_{1i}+A_1 =A_1+A_i \neq A_i$.
Thus, by scaling, i.e.  replacing $\vec v_i$ 
by $r_i^{-1} s_i \vec v_i$, we may assume that
$\vec w_{1i}  = \vec  v_1 - \vec v_i$ for $i>1$.    
Then, it follows     
$A_i=\IF \vec v_i$, $A_{ij} =\IF(\vec v_i -\vec v_j)$
for all $i\neq j$.
Indeed, by the sixth line of equations in Definition~\ref{d:frame}, 
for $1<i,j$ and $i \neq j$,
$ \vec 0\neq \vec v_i -\vec v_j
=\vec v_1-\vec v_j-(\vec v_1-\vec v_i)  \in (A_i +A_j) \cap
(A_{1i}+A_{1j})=A_{ij}$.
\end{proof}

\subsection{Recovering the Ground Field in Hilbert Lattices of Dimension $\geq3$} \lab{ss:IntField} 
The next step is to recover the field
$\IF$ within $\Gr(\calH)$ where
$d=\dim(\calH)\geq 3$. Let a $d$-frame
$\bar A$  be given as above. 
 Observe that with each of
the isomorphisms $\eta_{ij}$ from above and
each scalar $r \in \IF$ we also have the
linear map $r \eta_{ij}: A_i \rightarrow A_j$ 
(which is an  isomorphism if $r \neq 0$).   For simplicity of notation
we just consider the $r \eta_{12}$; but the discussion applies
 to any other choice as well.  Since we have seen above
that the $\eta_{ij}$ should be captured by their negative
graphs, associate with each scalar $r \in \IF$ 
a 1D subspace via
\begin{equation} \label{e:Ring}
\Interpret_{\bar A}:\IF\ni r\mapsto \IF(\vec v_1 - r\vec v_2) 
\;\in\; \IntRing_{\bar A} \;:=\;
\{ X \in \Gr(\calH)\mid X \cap A_2 =\Zero,\;X+A_2=A_1+A_2\} \enspace . 
\end{equation}
Conversely, if $X \in \IntRing_{\bar A}$ then $\dim(X)=1$
(since $X\not\in\{ 
\Zero , A_1+A_2\}$ and $\dim(A_1+A_2)=2$)
whence $X=\IF(x_1\vec v_1+x_2 \vec v_2)$
for some $x_1,x_2 \in \IF$. Assuming $x_1=0$ it would follow
$X \subseteq A_2$ and $A_1 \subseteq A_2$.
Thus $x_1 \neq 0$ and $X=\Interpret(r)$ with
$r= -x_1^{-1}x_2$. This demonstrates that
$\Interpret_{\bar A}:\IF \rightarrow \IntRing_{\bar A}  \subseteq  \Gr_1(\calH) $
is a bijection.  Observe that $\Interpret_{\bar A}(0)=A_1$ and 
$\Interpret_{\bar A}(1)=A_{12}$.

\begin{comment} 
Imagining $\calA$ and its projective
extension $\calP$ represented within $\calH$ as in
Observation~\ref{o:coopro}c), we  may consider
$\Interpret$ a bijection 
from $\IF$ onto    the  coordinate line $\ell_2$
of $\calA$: mapping $r$ onto the point with
coordinate $r$ on this line (recall that the
point with coordinate $1$, the unit point on this line,
is $A_{12}=\Interpret(1)$).  
Indeed, the (homogeneous)  line $\Interpret(r)$ of $\calH$ 
 meets the affine subspace  $\calA$ in this point.
For the following definition compare Observation~\ref{o:sub}. 
\end{comment} 
Now define  binary  operations $\ominus$ and $\otimes$ on $\Ring$
in terms of $\bar A$ and lattice operations  
such that $\Interpret$ becomes an isomorphism
of $(\IF,0,1,-,\cdot)$ onto $(\Ring,A_1,A_{12},\ominus,\otimes)$:
\begin{eqnarray*}
P \ominus Q &=& \bigl( [(Q_{13} +A_2) \cap (P+A_{23})] + A_3\bigr) \cap
(A_1+A_2),\quad\;\;
Q_{13}:= (Q+A_{23})\cap (A_1+A_3) \\     
P \otimes Q&=& (Q_{13} +P_{32})\cap (A_1+A_2)  
,\qquad\qquad\qquad\qquad\qquad\;\;\:
P_{32}:=(P+A_{13})\cap (A_2+A_3)
\end{eqnarray*}
Indeed, referring to the above basis one has
\[ A_i+A_j=\{x_i \vec v_i +x_j \vec v_j\mid x_i,x_j \in \IF\},\quad
 P=\IF(\vec v_1 -p \vec v_2),\quad   Q=\IF(\vec v_1 -q \vec v_2)
\] 
for unique $p,q \in \IF$ and calculates
\[\begin{array}{lclcl} 
Q_{13}&=&\{x\vec v_1 -xq \vec v_2 +y\vec v_2 -y \vec v_3\mid
 x,y \in \IF\}\cap (A_1+A_3)&=& \IF(\vec v_1-q\vec v_3)\\
P_{32}&=&\{x\vec v_1 -xp \vec v_2 +y\vec v_1 -y \vec v_3\mid
 x,y \in \IF\}\cap (A_2+A_3)&=& \IF(-p\vec v_2+\vec v_3)
\end{array} \]
  since 
$y-xq=0$ respectively $x+y=0$. It follows
for 
$S:=(Q_{13} +A_2) \cap (P+A_{23})$ and $P\ominus Q =
(S+A_3)\cap(A_1+A_2)$ that
\[\begin{array}{lcl} 
S&=&\{x\vec v_1+y \vec v_2- x\vec qx_3 \mid x,y \in \IF\}
\cap \{r\vec v_1-(rp-s)\vec v_2 -s\vec v_3\mid r,s \in \IF\}\\
&=&\{r\vec v_1-r(p-q)\vec v_2 -rq\vec v_3\mid r \in \IF\}\\
P \ominus Q&= &\{t\vec v_1-t(p-q)\vec v_2\mid t \in \IF\} = \Interpret(p-q)
\end{array}     \]
since $r=x$ and $s=xq=rq$. Finally,
observing   $y=xq$, one gets
\[ P\otimes Q=\{x\vec v_1 -yp\vec v_2- (xq-y)\vec v_3)\}
\cap (A_1+A_2)= \{x\vec v_1- xqp\vec v_2\mid x\in \IF\}=
\Interpret_{\bar A}(p\cdot q).    
\]
\begin{comment} 
The  term for multiplication has been defined
in order to generalize in Subsection~\ref{ss:MatinGrH}
to matrices resp. endomorphisms
when writing applications of maps on the left.
Doing so, one should  write scalars to the
right of vectors,
 i.e consider right vector spaces,  if 
commutativity of $\IF$ is not to be presumed.
\end{comment} 
We define $P \oplus Q=P \ominus (A_1 \ominus Q)$ 
and summarize the above. 

\begin{fact} \lab{f:IntField} 
For $\calH$ of $\dim(\calH)=d \geq 3$ with
orthogonal basis $\vec v_1,\ldots,\vec v_d$ and 
$\bar A$ a $d$-frame of $\Gr(\calH)$,
 $\Interpret_{\bar A}(r)=\IF(\vec v_1-r\vec v_2)$
constitutes an isomorphism  of the field $\IF$ with
operations $+,-,0,\cdot,1$ 
onto $\IntRing_{\bar A}$ with operations $\oplus_{\bar A} ,\ominus_{\bar A}, A_1, \otimes_{\bar A} , A_{12}$
defined as above  in terms of $\bar A$. 
\end{fact}
There is of course nothing particular about
the indices $123$. 
So Fact~\ref{f:IntField} holds also
for any triple $ijk$ of pairwise distinct indices
between $1$ and $d$.
\begin{comment} We will refer, informally,  to results
as in Fact~\ref{f:IntField} as \emph{interpretations} 
- established by  maps $\Interpret_{\bar c}$ (with or without parameters
  $\bar c$)  
and formulas relating the two structures
(cmp. Subsection~\ref{ss:Interpret}).
\end{comment}

\subsection{Strong Satisfiability is Hard for \BSS--$\calNP$ in Dimension $\geq3$} \lab{ss:realBSSNPc}
Having Fact~\ref{f:IntField},
a well-known method applies  to translate  information about  
the ring   $\IF$ into such about
$\Gr(\calH)$, given a  $d$-frame $\bar A$
($d =\dim \calH \geq 3$); and we observe that both 
the translation and the conditions for a $d$-frame
can be encoded in polynomial time:

\begin{myproposition} \lab{p:IntField}
Let $\calH$ denote a $d$-dimensional $\IF$-unitary vector space,
$d\geq3$.
\begin{enumerate}
\item[a)] $\SAT_{\Gr(\calH)}$ is $\calBP(\calNP_{\IF}^0)$--hard
\item[b)] $\SAT_{\Gr(\IF^d),\Gr(\IF^d)}$ is $\calNP_{\IF}$--hard.
\item[c)]  
Claim~a) holds uniformly in the following sense: 
Given $3\leq d\in\IN$
and a finite system $\bar p=0$ of 
of polynomial equations $p_j(X_1,\ldots,X_n)=0$
with $p_j\in\IZ[X_1,\ldots,X_n]$, 
a Turing machine can produce within 
time polynomial in $d$ and in the length of the binary encoding
of the $p_j$'s coefficients and monomials,
an orthologic term $t_{\bar p,d}$ such that,
for every field $\IF\subseteq\IC$ and every $d$-dimensional
$\IF$-unitary vector space $\calH$,
the $p_j$ admit a common root in $\IF$
~iff~ $t_{\bar p,d}$ admits a
strongly satisfying assignment in $\Gr(\calH)$.
\end{enumerate}
\end{myproposition}
\begin{proof}
Claim~c) asserts $\calBP(\calNP_{\IF}^0)$--complete
$\FEAS_{\IZ,\IF}$ to be
polynomial-time reducible to $\SAT_{\Gr(\calH)}$,
hence a) follows from c).
Concerning c) itself, consider the following algorithm:
\begin{enumerate}
\item[i)] Replace each integer coefficient $c$ 
  by an expression over $(1,+,-,\cdot)$ 
  according to Observation~\ref{o:MultChain} below,
  thus obtaining $n$-variate terms $p_j'$ in the
  language of rings equivalent to $p_j$ over $\IF$.
\item[ii)] In these terms, replace each constant/operation
  $+,-,0,\cdot,1$ by the corresponding term
  $\oplus_{\bar A} ,\ominus_{\bar A}, A_1, \otimes_{\bar A} , A_{12}$
  according to Fact~\ref{f:IntField},
  thus obtaining terms $t_{j}(X_1,\ldots,X_n;\bar A)$ 
  in the language of lattices.
\item[iii)] Now take the lattice equations
  $t_{j}(X_1,\ldots,X_n;\bar A)=A_1$ and
\item[iv)] for each variable $X_i$ add conditions 
  $X_i\wedge A_2=\Zero$ and $X_i\vee A_2=A_1\vee A_2$ 
  from Equation~(\ref{e:Ring})
\item[v)]
  as well as those from Definition~\ref{d:frame}
  for $\bar A$ to constitute a $d$-frame
\item[vi)] and combine the resulting system of 
  lattice equations into a single ortholattice
equation  ``$t_{\bar p}=\One$''
  according to Fact~\ref{f:clear}a).
\end{enumerate}
First, a Turing machine can indeed perform each of the above
steps in time polynomial in $d$ and in
the joint binary length of the $p_j$.

Secondly, a $d$-frame $\bar A$ exists according to Fact~\ref{f:frame}e)
and hence complies with (v).
Structural induction on (ii) yields for any $b_1, \ldots ,b_n$
\begin{equation}  \label{e:Interpret}
 t_j\big(\Interpret_{\bar A}(b_1),\ldots \Interpret_{\bar A}(b_n)\big)
 \;=\; \Interpret_{\bar A}\big(p_j'(b_1,\ldots b_n)\big)
  \overset{\text{(i)}}{=} \Interpret_{\bar A}\big(p_j(b_1,\ldots b_n)\big)
  \enspace. 
\end{equation}
Thus, a common root $b_1,\ldots,b_n\in\IF$ of all $p_j$
gives rise via $B_i:=\Interpret_{\bar A}(b_i)$ to an assignment 
$B_1,\ldots,B_n$ in $\Gr(\calH)$ satisfying the requirements of 
(iv) as well as (iii)  and
\begin{equation}  \label{e:Interpretb}
 t_j\big(\Interpret_{\bar A}(b_1),\ldots \Interpret_{\bar A}(b_n)\big)
  \;=\; \Interpret_{\bar A}(0)\;=\; A_1 \enspace.  
\end{equation}
 leading to a strongly satisfying assignment
$\bar B,\bar A$ over $\Gr(\calH)$ of $t_{\bar p}$.

Conversely, a strongly satisfying assignment over $\Gr(\calH)$ 
of $t_{\bar p}$ consists by (v) of a $d$-frame $\bar A$ and,
due to (iv), of $B_1,\ldots,B_n\in\IntRing_{\bar A}$, that is
$B_i=\Interpret_{\bar A}(b_i)$ for certain $b_1,\ldots,b_n\in\IF$.
Moreover, $t_j(\bar B)=A_1$ due to (iii) implies  
$p_j(\bar b)=0$ by reading Equations~(\ref{e:Interpret}) and (\ref{e:Interpretb})  backwards:
$b_1,\ldots,b_n$ a common root of all $p_j$.

Concerning the case (b) with constants,
let $p_j$ denote polynomials over $\IF$.
Now, instead of (v), fix $\bar A$ to denote the $d$-frame
corresponding to the standard basis 
$\vec v_i=(0,\ldots,0,1,0,\ldots,0)$ of $\IF^d$ according to
Fact~\ref{f:frame}d); and, instead of (i), replace the
polynomials' coefficients $c\in\IF$ by constants 
$C:=\Interpret_{\bar A}(c)=\IF(0,\ldots,c,\ldots,-c,\ldots,0)\in\Gr_1(\IF^d)$.
\qed\end{proof}

\begin{observation} \lab{o:MultChain}
To each $c\in\IN$ there exists a term $t_c$ over $(1,+,\cdot)$
of length $|t_c|\leq\calO(\log c)$ evaluating to $c$
over each ring containing $\IN$.
Moreover, such $t_c$ can be computed from $c$ in time
polynomial in the binary length of $c$.
\end{observation}
Note that an \emph{addition chain} 
`expresses' $c$ as unique solution to a system of equations
of logarithmic size; 
however we are interested in a single term
and thus employ also multiplication.

\begin{proof}
By induction, claiming $|t_c|\leq 2+7\log_2(c)$:
Indeed $2c=(1+1)\cdot t_c=:t_{2c}$ 
and $2c+1=(1+1)\cdot t_c+1=:t_{2c+1}$ 
both have length at most $7+|t_c|\leq 
7+2+7\log_2(c)=2+7\log_2(2c)$ by induction hypothesis.
\end{proof}
\begin{corollary} \lab{c:realBSSNPc}
In case $\IF\subseteq\IR$ and for every $d\geq3$,
$\SAT_{\Gr(\IF^d)}$ is $\calBP(\calNP^0_{\IF})$--complete;
and $\SAT_{\Gr(\IF^d),\Gr(\IF^d)}$ is $\calNP_{\IF}$--complete.
\end{corollary}
In particular, the decidability of $\SAT_{\Gr(\IQ^3)}$ is as
open as that of $\FEAS_{\IQ}$, recall Fact~\ref{f:BSSBP}f).
And $\SAT_{\Gr(\IR^d)}$ is complete for
$\calBP(\calNP_\IR^0)$ and perhaps closest in spirit 
to the Boolean satisfiability problem complete for $\calNP$.

Note however the gap between 
Proposition~\ref{p:UpperComplexity}c) and Proposition~\ref{p:IntField}
in case $\IF\not\subseteq\IR$: to be closed in the sequel.

\subsection{Orthonormal Frames} \lab{ss:OrthonFrame}
In order to recapture $\IF \not\subseteq \IR$ as a $\ast$-ring
within $\Gr(\IF^d)$, i.e.
taking into account the involution of $\IF$,
  we have 
to include orthogonality in our concept
 of frames. The  first step is again due
 to \person{von Neumann}.

\begin{mydefinition} \lab{d:OrthogFrame}
For $L$ a modular ortholattice, 
call a $d$-frame $\bar a$ of $L$ \textsf{orthogonal}
if $a_i \leq \neg a_j$ holds for all $i \neq j$ 
(which makes the second line of conditions in
Definition~\ref{d:frame} redundant).
\end{mydefinition}
In the situation of Fact~\ref{f:frame}d),
orthogonality of  $\bar A$ 
corresponds to orthogonality of the basis $\vec v_i$.
What we want, to make things a simple as possible,
is a concept of frames which corresponds to
equinormal orthogonal  bases. So, assume that 
an orthogonal 
basis $\vec v_1, \ldots ,\vec v_d$ is given
with associated orthogonal $d$-frame  $\bar A$.  We
refer  to Fact~\ref{f:IntField} and claim that
$\|\vec v_1\|^2=\|\vec v_2\|^2   \quad\Leftrightarrow\quad
\Interpret_{\bar A} (1) = A_{12}  \;\perp \;
\Interpret_{\bar A}(-1)=\ominus_{\bar A} A_{12}$:
Recall Observation~\ref{o:orthogframe} for the geometric intuition
or calculate, with $\langle\vec v_1 \mid\vec v_2 \rangle =0$ in mind,
 \[0 =\langle \vec v_1 -\vec v_2\mid \vec v_1+\vec v_2\rangle
=   \langle\vec v_1\mid\vec v_1 \rangle  +\langle\vec   v_1\mid\vec v_2 \rangle
-  \langle\vec v_2\mid\vec v_1   \rangle- \langle\vec v_2\mid\vec v_2
  \rangle
= \langle\vec  v_1\mid\vec v_1 
\rangle -
\langle\vec  v_2 \mid\vec v_2 \rangle.  \] 
By Fact~\ref{f:IntField}, it holds
$\ominus_{\bar A }A_{12} =A_1 \ominus  A_{12}
= \bigl([(A_{13}+A_2)\cap (A_1+A_{23})] +A_3\bigr)\cap (A_1+A_2)  $.
In view of the last paragraph of Subsection~\ref{ss:IntField},
 similarly to $\ominus_{\bar A 12} \overset{\text{def}}{=}
 \ominus_{\bar A}$ we may  consider $\ominus_{\bar A 1k} A_{1k}$  for all $2<k\leq d$.
This motivates  defining the following
terms  $\ominus_{\bar y 1k}y_{1k}$ 
in variables  $\bar y=(y_{ij}\mid 1\leq i,j \leq d)$
where index $1<\ell\leq d$ is distinct from $1,k$ but otherwise arbitrary:
\[ \ominus_{\bar y 1k}y_{1k} :=
 \bigl([(y_{1\ell}\vee y_k)\wedge (y_1\vee y_{k\ell})] \vee y_\ell\bigr)\wedge
 (y_1\vee y_k).\]

\begin{mydefinition} \lab{d:OrthonFrame}
An orthogonal $d$-frame
$\bar a$  of a modular ortholattice  is \textsf{orthonormal} 
if $\ominus_{\bar a 1k}a_{1k} \leq \neg a_{1k}$ for all $k\geq 2$. 
\end{mydefinition} 
Fact~\ref{f:frame}d) and the above calculation   prove

\begin{mylemma} \lab{l:OrthonFrame}
For $d\in\IN$, a  system 
$\bar A =(A_{ij}\mid 1\leq i,j \leq d)$  
is an orthonormal $d$-frame of $\Gr(\calH)$   
if and only if there is an equinormal orthogonal 
basis $\vec v_1, \ldots ,\vec v_d$ of
$\calH$ such that 
$A_{ii}=\IF\vec v_i$ and
$A_{ij}=\IF(\vec v_i-\vec v_j)$   for $i\neq j$, $i,j=1,\ldots, d$.  
\end{mylemma} 
For $\calH=\IF^d$ with the canonical basis this frame is
the \textsf{canonical} orthonormal $d$-frame of $\Gr(\IF^d)$.

Now,  the image of an orthonormal $d$-frame under an isomorphism of
 $\Gr(\IF^d)$ onto $\Gr(\calH)$ is again an orthonormal 
$d$-frame; 
which yields an equinormal orthogonal basis for $\calH$
according to Lemma~\ref{l:OrthonFrame}.
We have thus concluded the promised partial converse
to Lemma~\ref{l:Equinormal}a):

\begin{corollary} \lab{c:Equinormal}
\begin{enumerate}
\item[a)]
Let $\calH$ be a $d$-dimensional unitary
$\IF$-space and $\Gr(\calH)$ isomorphic to $\Gr(\IF^d)$.
Then $\calH$ admits an equinormal orthogonal base
(i.e. is a `scaled' isometric copy of $\IF^d$).
\item[b)]
Suppose $\calH$ and $\calH'$ are finite-dimensional
unitary $\IF$-spaces such that $\Gr(\calH)$ embeds
into $\Gr(\calH')$. Then $\dim(\calH)$ divides
$\dim(\calH')$.
\end{enumerate}
\end{corollary}
In a) we even have $\calH$ and $\IF$ `isometric up to scaling'.
Concerning Claim~b), take a $d$-frame $\bar a$ of $\Gr(\calH)$,
$d:=\dim(\calH)$. Its image $\bar b$ in $\Gr(\calH')$ under
the embedding thus also constitutes a $d$-frame;
hence $\dim(\calH')=d\cdot\dim(b_1)$ according to
Fact~\ref{f:frame}c). \qed

\subsection{Recovering the $\ast$-Field in Dimension $\geq 3$}\lab{ss:IntStarField}
Based on the above preparation, the interpretation 
of $\IF$ in $\Gr(\calH)$  from Subsection~\ref{ss:IntField}
can be extended to respect also involution.
In the more general and abstract setting
of projective spaces with an orthocomplementation
on the subspace lattice,
this is due to Birkhoff and von Neumann
\cite{Birkhoff}.  Yet, in our setting there
is no need for abstract calculations.
Recall Fact~\ref{f:IntField}. 

\begin{mylemma} \lab{l:IntStarField} 
Let $d \geq 3$.
For any   $\IF$-unitary space $\calH$
of dimension $d$ admitting an equinormal orthogonal basis,
 there is some
orthonormal $d$-frame $\bar A$ of $\Gr(\calH)$.
Given such,  there is 
a bijection $\Interpret_{\bar A}: \IF \rightarrow
\IntRing_{\bar A}$ 
 such that the following hold
\begin{enumerate}  
\item[ a)] $\Interpret_{\bar A}$ is an isomorphism
of the  ring  $(\IF,+,-,0,\cdot,1)$ onto 
$(\IntRing_{\bar A},\oplus_{\bar A},\ominus_{\bar A},A_1,\otimes_{\bar A}, A_{12})$.
\item[ b)] For  $P=\Interpret_{\bar A}(p)$ and
$Q=\Interpret_{\bar A}(q)$  with $p,q \neq 0$ in $\IF$ one has
$q =\bar p$, the conjugate of $p$, if and only if
there is $R\in \IntRing_{\bar A}$ such that
$Q\otimes_{\bar A} R=\ominus_{\bar A} A_{12}$ and $P \perp R$ or,
equivalently, 
 $Q\otimes_{\bar A} ((A_1 \vee A_2) \wedge \neg P)=\ominus_{\bar A} A_{12}$.
\COMMENTED{
if and only there is $R\in \IntRing_{\bar A}$ such that
$Q\otimes_{\bar A} R=\ominus_{\bar A} A_{12}$  and $P= \neg R \wedge
(A_1\vee A_2)$.}  
\end{enumerate} 
\end{mylemma}
\begin{proof} 
The existence of $\bar A$ 
follows from Lemma~\ref{l:OrthonFrame};
and conversely it is no loss of generality to
assume that $\bar A$ is induced by an equinormal orthogonal basis.
Now a) is just a special case of Fact~\ref{f:IntField}. 
In b) we first show that $q=\bar p$ 
if and only there is $R\in \IntRing_{\bar A}$ such that
$Q\otimes_{\bar A} R=\ominus_{\bar A} A_{12}$  and $P\perp R$.
Namely, 
 for $R=\Interpret_{\bar A}(r)$ we have 
$Q\otimes R=\ominus A_{12}$ if and only if $qr=-1$, 
i.e. $R=\IF(q \vec v_1+\vec v_2)$. 
Moreover,   $P \perp R$ amounts to
\[0= \langle \vec v_1 -p \vec v_2 \mid q\vec v_1 +\vec v_2
 \rangle 
= \langle \vec v_1\mid q\vec  v_1 \rangle 
+\langle \vec v_1\mid \vec v_2 \rangle 
-\langle p\vec v_2 \mid q \vec v_1  \rangle
- \langle p \vec  v_2  \mid \vec v_2 \rangle 
 = q \|v_1\|^2 - \bar p \|v_2\|^2 \]
i.e. to $q=\bar p$. 
Conversely,  if $q=\bar p$ choose
$R=\Interpret_{\bar A}(-q^{-1})$.  
Now, observe that here $R$ is uniquely determined 
as the inverse of $\ominus_{\bar A}Q$ and that
 $P \perp R$ means $R=(A_1\vee A_2) \wedge \neg P$.
\end{proof} 
We can now extend Proposition~\ref{p:IntField}
to match with Proposition~\ref{p:UpperComplexity}c)
also in the complex case. 
Here, a  \textsf{$\ast$-polynomial} 
is a term in the language of (commutative) $\ast$-rings, 
i.e. may involve involution as in
$p(X,Y)=(X+Y^\dagg)^2/4+(X-X^\dagg)^2/2+(Y+Y^\dagg)^2/4+(Y-Y^\dagg)^2/2$.
Considering substitutions in the  $\ast$-field $\IF$,
it is convenient to define
such polynomials as members of
the polynomial ring $\IF[X_1,X_1^\dagg,\ldots,X_n,X_n^\dagg]$
where $X^\dagg$ is a variable to be interpreted in $\IF$ 
 as the conjugate  $a^\ast$  if $X$ is interpreted as $a$.

\begin{theorem} \lab{t:IntStarField}
\begin{enumerate}
\item[a)]
Given $d\geq3$ and finitely many $\ast$-polynomials 
$p_j\in\IZ[X_1,X_1^\dagg,\ldots,X_n,X_n^\dagg]$,
a Turing machine can within 
time polynomial in $d$ and in the length of the binary encoding
of the $p_j$'s coefficients and monomials
produce 
an orthologic term $t_{\bar p,d}$ such that,
for every $\ast$-field $\IF\subseteq\IC$ and every 
$d$-dimensional $\IF$-unitary vector space $\calH$
admitting an equinormal orthogonal basis,
the $p_j$ have a common root in $\IF$
~iff~ $t_{\bar p,d}$ admits a
strongly satisfying assignment in $\Gr(\calH)$.
\item[b)] For every $\IF$-unitary vector space $\calH$ of
dimension $d\geq3$ admitting an equinormal orthogonal basis,
$\SAT_{\Gr(\calH)}$ is $\calBP(\calNP_{\Re\IF}^0)$--complete.
\item[c)] For $d\geq3$, satisfiability with constants
$\SAT_{\Gr(\IF^d),\Gr(\IF^d)}$ is $\calNP_{\Re\IF}$--complete.
\end{enumerate}
\end{theorem}
\begin{proof}
\begin{longenum}
\item[a)]
In view of Lemma~\ref{l:IntStarField} extending Fact~\ref{f:IntField},
modify the algorithm in the proof of Proposition~\ref{p:IntField} as follows:
\begin{enumerate}
\item[v')]
Add equations from Definitions~\ref{d:frame},
\ref{d:OrthogFrame}, and \ref{d:OrthonFrame}
for $\bar A$ to constitute an orthonormal $d$-frame.
\item[iv')] 
For each variable $X_i$ add conditions $X_i\wedge A_2=\Zero$ and $X_i\vee A_2=A_1\vee A_2$. \\
Do similarly for each of the
newly introduced variables $Y_1,\ldots,Y_n$ and $Z_1,\ldots,Z_n$ 
and add equations
$Y_k\otimes_{\bar A} Z_k=\ominus_{\bar A} A_{12}$ and
$X_k=\neg Z_k\wedge(A_1\vee A_2)$. Finally replace 
any starred occurrence $X_i^\dagg$ of $X_i$ with $Y_i$.
\end{enumerate}
\item[b)]
by reduction from $\FEAS_{\IZ,\Re\IF}$:
To the given polynomials $p_j(X_1,\ldots,X_n)$ over $\IZ$
add $\ast$-polynomial equations $X_i=X_i^\dagg$ ($1\leq i\leq n$),
thus restricting solutions from $\IF$ to $\Re\IF$.
Now apply a) to conclude $\calBP(\calNP^0_{\Re\IF})$--hardness
and combine with Proposition~\ref{p:UpperComplexity}c).
\item[c)]
Combine b) with the (proof of) Proposition~\ref{p:IntField}b),
i.e. use the constants permitted in the term to encode the
fixed orthonormal standard frame $\bar A$ as well as the 
polynomials' coefficients' $c\in\IF$ by constants 
$C:=\Interpret_{\bar A}(c)$:
this proves $\calNP_{\Re\IF}$--hardness,
for the upper bound recall
Proposition~\ref{p:UpperComplexity}d).
\qed\end{longenum}\end{proof}

\begin{corollary} \lab{c:complexBSSNPc}
$\SAT_{\Gr(\IC^d)}$ is $\calBP(\calNP^0_{\IR})$--complete;
and $\SAT_{\Gr(\IC^d),\Gr(\IC^d)}$ is $\calNP_{\IR}^0$--complete.
For any fixed $\IF$ and $k,d\geq3$, 
$\SAT_{\Gr(\IF^d)}$ is polynomial-time equivalent to $\SAT_{\Gr(\IF^k)}$.
\end{corollary}
%


\cx{
\begin{digression} \lab{d:ax} 
\begin{longenum}
\item[a)]
From 
Theorem~\ref{t:IntStarField} one can derive, for any $d \geq 4$,  an
 axiomatization of the first order theory of 
$\Gr(\IR^d)$  --- namely  as that of  modular
 ortholattice of height $d$ admitting an orthonormal
$d$-frame such that the coordinate $\ast$-ring $\IntRing$
is a real closed  field with identity
involution.  This follows  as a special case of the
coordinatization 
of projective spaces  and  the result of \cite{Birkhoff} 
cmp. \mycite{Thm.13.2}{coord}.
For $\Gr(\IC^d)$ the coordinate $\ast$-ring $\IntRing$ 
is required to be a field  having real closed subfield
$\IntRing_0=\{x \in \IntRing\mid x=x^*\}$ and an element $i$ 
such that  $i^2=-1$, $\IntRing=\IntRing_0+i\IntRing_0$, 
and $(a+ib)^*=a=-ib$.  Of course, there is no finite axiomatization.
\item[b)] Concerning axiomatization, 
for $d=3$ the modular law has to be replaced
by the stronger Arguesian law. The easiest
cases are $d=1$ (all height $1$ ortholattices) and
$d=2$ (all infinite height $2$ ortholattices). 
\item[c)] The question of axiomatizating the universal or 
equational theory in fixed dimension $d\geq 3$ leads to open problems cmp.
Digression~\ref{d:uclass}.  In contrast,
the equational theory of $\{\Gr(\IF^d)\mid d<\infty\}$  
 can be  recursively axiomatized (whence  decided) 
 for $\IF \in \{\IR,\IC\}$ cmp. Fact~\ref{f:Tarski},  
Section~\ref{s:Neumann}, and \cite{Marina}.  
On the other hand, the universal (Horn) theory
of this class is not recursively axiomatizable (cmp. \mycite{Corollary 7.2}{hn}).  
\end{longenum}
\end{digression} }

\section{More Background from Linear Algebra, Logic,
and Modular Orthologic}
\lab{s:Facts}
We recall some definitions and results from Linear Algebra
(Subsection~\ref{ss:endo})  and Logic (Subsection~\ref{ss:tool}) 
 to be
referred to in the sequel. Also, we  discuss 
 (Subsection~\ref{ss:Interpret}) 
a particularly
simple instance of the concept of interpretation, relational interpretation,
which captures the idea  of evaluating a term by
keeping track of the intermediate results 
(and hint upon the general concept and
 the fact that many of our results may be seen as
instances).   
The main contribution is adapting (in Subsection~\ref{ss:Huhn}) a variant of the concept of 
$d$-frames,
\person{Huhn}'s $d$-diamonds,  
to capture  (in Subsection~\ref{ss:Quantifiers})  the quantifier free  part  of the first order
theory  of $d$-dimensional MOLs  by equations, quantified
existentially  or, by choice, universally.

\subsection{Endomorphism Rings, Adjoints, and Matrix Units} \lab{ss:endo} 

Matrices form a non-commutative ring
naturally equipped with transposition
(and possibly complex conjugation).
For any  $\IF$-vector spaces $U$, $V$ 
let $\Hom(U,V)$ denote the set of all $\IF$-linear
maps $\phi:U \rightarrow V$.
In particular, $\End(U)=\Hom(U,U)$ 
is the endomorphim ring  of $U$.
Given a basis $\vec v_1, \ldots ,\vec v_d$ of
$\calH$, there is a canonical
isomorphism of $\End(\calH)$ onto  the matrix ring $\IF^{d \times d}$,
the matrix $A$ of $\phi$ has as its columns the
coordinate vectors of the $\phi(\vec v_j)$.

Recall that (cmp. \mycite{\S 2.13}{Rowen}, \mycite{\S 5.1}{Farenick})
a $\ast$-\textsf{ring} is a ring $R$ equipped
with an \textsf{involution} $r \mapsto  r^\ast$,
i.e. a map  such that  $(r^\ast)^\ast=r$,
$(r+s)^\ast=r^\ast+s^\ast$, and 
$(rs)^\ast =s^\ast r^\ast$.
Thus, the  ring $\IF^{d \times d}$ of $d$-by-$d$-matrices over 
the $\ast$-subfield $\IF$ of $\IC$  
becomes a $\ast$-ring with the 
involution $A \mapsto B= A^\ast$, the conjugate transpose
 of $A$, i.e. $b_{ij}= \bar a_{ji}$. Occasionally,
we write $A^\dagger =A^\ast$. 
In the setting of $\IF$-unitary spaces, 
the involution is given as adjunction.

\begin{mydefinition}  \lab{d:adjoints} 
(Cf. \mycite{\S11-13}{Gelfand}, \mycite{\S 73}{Halmos},
\mycite{Ch.7}{Axler},  
\mycite{1.20}{Farenick}).  
\begin{enumerate}
\item[a)]
Given subspaces $U,V$ of  $\calH$,
linear maps $\phi:U \rightarrow V$ and
$\psi:V \rightarrow U$ are \textsf{adjoint} to each other
if 
$\langle \vec u \mid \phi(\vec u) \rangle =
\langle \psi(\vec v)\mid \vec v\rangle$
for all $\vec u \in U,\,\vec v \in V$.
\item[b)]
$\phi:U \rightarrow V$ is an \textsf{isometry} if
$\langle \vec x \mid \vec y \rangle =
\langle \phi(\vec x)\mid \phi(\vec y)\rangle$
for all $\vec x,\;\vec y \in U$.
\item[c)]
An \textsf{idempotent} in a ring is an element $e$ such that
$e^2=e$. 
\item[d)]
Dealing with  $\ast$-rings,
 call $r$ \textsf{selfadjoint} if
$r=r^\ast$. 
\item[e)] An idempotent $e$    
 is a \textsf{projection} 
if it is also selfadjoint.
\end{enumerate}
\end{mydefinition}

\begin{fact} \lab{f:adjointa} 
Consider 
$U,V,W\in \Gr(\calH)$, $\calH$ a finite-dimensional $\IF$-unitary space.
\begin{enumerate} 
\item[ a)] If $\phi\in \Hom(U,V)$ has an adjoint $\psi$ 
then $\psi$ is unique and denoted by $\phi^*$.
\item[ b)] If $\phi\in \Hom(U,V)$ and $\psi\in \Hom(V,W)$ 
have adjoints, then so does $\psi \circ \phi \in Hom(U,W)$ and it 
holds $(\psi\circ \phi)^*= \phi^* \circ \psi^*$. 
\item[ c)] Any $\phi \in \End(\calH)$ has an adjoint $\phi^*$: 
If $\phi$ has matrix $A$ w.r.t.  the  orthogonal basis
$\vec v_1, \ldots, \vec v_d$ then 
$\phi^*$ has matrix  $D^{-1}\cdot A^\adjoint\cdot D$, 
$D$  diagonal with diagonal entries $\langle \vec  v_i \mid \vec v_i  \rangle$. 
\item[ d)] $\End(\calH)$ is a $\ast$-ring  with adjunction as involution.
\item[ e)]  Any equinormal orthogonal basis $\vec v_1,\ldots,\vec v_d$ of $\calH$  
establishes an isomorphism of $\ast$-rings
$\End(\calH)\ni\phi\mapsto A\in\IF^{d \times d}$
such that $A$ is the matrix of $\phi$ with respect to this basis. 
\item[ f)] $\psi \in \Hom(U,V)$ is an isometry if and only if
 it is a linear isomorphism and $\phi^{-1}=\phi^\ast$.
\item[g)] For any subspace $U$ of $\Gr(\calH)$ there is a
unique projection $\pi$ such that $U =\range(\pi)$;
moreover $\neg U =\range(\id -\pi)$ and,
 for any projection $\phi$,  $U \subseteq \range \phi$ 
if and only if  $\phi|U =\id_U$ if and only if 
 $\pi \circ \phi =\pi$ if and only if $\phi \circ \pi =\pi$.  
\end{enumerate} 
\end{fact} 

\begin{proof}  a) to g) are standard  Linear Algebra.
\end{proof} 
Matrix units provide  a well known method 
to establish  an isomorphism of a ring or $\ast$-ring
onto a suitable matrix ring. 
\begin{myexample} \lab{x:matunit}
For $1\leq i,j\leq d$ let $E_{ij}$ denote the $d\times d$-matrix
mapping the $i$-th canonical basis vector onto the $j$-th and
the others to $\vec 0$, i.e.
\[ E_{ij} = (\delta_{i\ell}\cdot\delta_{jk})_{_{k\ell}} \;\;\mbox{ where }   
\delta_{i\ell}=\Big\{ \begin{array}{ll} 1 &\mbox{ if\; } i=\ell\\
0 &\mbox{ if\; } i \neq\ell  \end{array} \]
\end{myexample}
Note the deliberate transposition implicit in this example
to make it comply with 

\begin{mydefinition} \lab{d:matunit}
In a fixed (not necessarily commutative) ring $\Ring$,
elements $(\vep_{ij}\mid 1\leq i,j \leq d)$ 
form a $d$--\textsf{system} $\bar \vep$  of \textsf{matrix units}  if
\begin{equation} \label{e:matunitr}
\sum\nolimits_{i=1}^d \vep_{ii}=1,\;\; 
\vep_{k \ell}\cdot\vep_{ij} = \delta_{jk}\vep_{i\ell}
\end{equation}
In particular, the $\vep_{ii}$ are idempotents.
For  $\Ring$ a   $\ast$-ring, $\bar\vep$ is a $d$-system 
of $\ast$-\textsf{matrix units} if, in
  addition, $\vep_{ij}^\ast= \vep_{ji}$ for all $i,j$; in particular, the
$\vep_{ii}$ are projections.  
\end{mydefinition} 
Following common (though not reasonable) use
to write maps $\phi:X \rightarrow Y$ as $\phi(x)=y$ 
and composition as $\psi\big(\phi(x)\big)=\big(\psi \circ \phi\big)(x)$,
on indexed domains $X_i$ this naturally leads to
$\phi_{ij}:X_i \rightarrow X_j$   and composition
$\phi_{ik}=\phi_{jk}\circ \phi_{ij}$. 

\begin{fact} \lab{f:matunit}
\begin{enumerate}
\item[ a)]
Given an idempotent (projection) $\vep \in \End(\calH)$ one has
a direct (orthogonal)  decomposition $\calH= \range(\vep) \oplus
\range(\id -\vep)$ and 
 obtains a ring ($\ast$-ring) 
$\vep\circ\End(\calH)\circ\vep:=\{\vep \circ \phi \circ \vep\mid \phi
 \in \End(\calH)\}$, with unit $\vep$ while
all other structure is  inherited from $\End(\calH)$, 
which is isomorphic to
$\End(\range \vep)$. 
\item[ b)] Given a  $d$-system $(\vep_{ij}\mid 1\leq i,j \leq d)$ 
of matrix ($\ast$-matrix)  units in $\End(\calH)$
one has 
$\calH$
a (orthogonal)  direct sum of the  $U_i=\range(\vep_{ii})$, and
the restrictions $\eta_{ij}:=\vep_{ij}|U_j$ 
isomorphisms (isometries) $\eta_{ij}:U_i \rightarrow U_j$
such that $\eta_{jk} \circ\eta_{ij} =\eta_{ik}$.
In particular $\eta_{ii}=\id_{U_i}$ and $\eta_{ji} 
=\eta_{ij}^{-1}$ ($=\eta_{ij}^*$) and $\dim(U_i)$ does not depend on $i$.
Conversely, given such $U_i$ and $\eta_{ij}$,
one obtains a $d$-system  of matrix ($\ast$-matrix) units via $\vep_{ij}:=\eta_{ij} \circ
\pi_{i}$ where $\pi_i(\vec v) =\vec u_i$ for $\vec v =\sum_j \vec
u_j$ with $\vec u_j in U_j$.   
 \item[ c)] Any (orthogonal equinormal)   basis
$\vec v_{ik}$ ($1\leq i\leq d;1\leq k\leq n$) of $\calH$ 
determines a $d$-system of matrix ($\ast$-matrix)  units where
$U_i=\sum_{k=1}^n \IF \vec v_{ik}$ and $\vep_{ij}(\vec v_{ik})
=\vec v_{jk}$. If $d=\dim \calH$,
every $d$-system of
matrix ($\ast$-matrix)  units arises in this way.
\item[d)] Given a $d$-system $\bar \vep$  of matrix
($\ast$-matrix) units in $\End(\calH)$ 
where $\dim(\calH)=m<\infty$, one has unique $\ell$ with
$m=\ell d$  and a ring ($\ast$-ring) isomorphism  of $\IF^{\ell \times \ell}$ onto
$\End(\range \vep_{11})$.  
\item[e)] If $\vep$ is a $d$-system of matrix
($\ast$-matrix)  units in 
$\End(\calH)$ and $\dim(\calH)=d$ then $a \mapsto a \id$   
is an ring ($\ast$-ring)  isomorphism of $\IF$
onto   $\{ \phi \in \End(\calH)\mid \phi \circ \vep_{ij}
= \vep_{ij} \circ \phi \mbox{ for all } i,j\}$.  
\end{enumerate} 
\end{fact} 
\begin{proof} a) Clearly,  $\id -\vep$ is an idempotent
(projection) , too, and 
yields the decomposition (which is orthogonal if $\vep$ is a 
projection). Moreover,
$\vep\circ\End(\calH)\circ\vep$ is 
a ($\ast$-)  subring (except unit) of $\End(\calH)$. 
 Mutually inverse isomorphisms between this   and $\End(U)$, 
$U:=\range(\vep)$,
  are given by
\[ \vep \circ \phi \circ \vep  \mapsto   (\vep \circ \phi \circ
\vep)|U,\quad \psi \mapsto \hat{\psi} \quad \mbox{ where } \hat{\psi}(\vec u
+\vec v):= \psi(\vec u) \mbox{ for } \vec u \in U,\,\vec v\in
\range(\id -\vep)\] 
b) Given a $d$-system of ($\ast$-)  matrix units, the first
condition yields $\vec x= \sum_{i=1}^d \vep_{ii}(\vec x)$ whence
$\calH =\sum_i U_i$ and this  sum is direct (orthogonal)  by the second
condition. Moreover, given $\ast$-matrix units, for 
 $i\neq j$ we have
$\big\langle \vep_{ii}(\vec x) \big|  \vep_{jj}(\vec y) \big\rangle
=\big\langle \big(\vep_{ji}\circ \vep_{jj}\circ \vep_{ij}\big)(\vec x) \big|
\vep_{jj}(\vec y) \big\rangle
= \big\langle \big(\vep_{jj}\circ    \vep_{ij}\big)(\vec x) \big|
\big(\vep_{ji}^* \circ \vep_{jj}\big)(\vec y) \big\rangle =0$
since $\vep_{ji}^*=\vep_{ij}$ and $\vep_{ij} \circ \vep_{jj}=0$:
$U_i \perp U_j$.  Also,  $\eta_{ij}^*=\eta_{ji}$
since $\vep_{ij}^*=\vep_{ji}$.
 The remaining claims follow from
similar calculations.  
 c) and d) are then   obvious.
 In e) choose a basis
according to c) thus turning $\vep_{ij}$ into $E_{ij}$. 
\end{proof}  

\begin{fact}  \lab{f:lat}
 Let $\IF$ be a $\ast$-subfield of $\IC$.
\begin{enumerate}
\item[a)] Any subspace $U$ of $\IF^d$ 
is the range of some $A \in \IF^{d \times d}$. 
\item[b)]
For any $A,B \in \IF^{d \times d}$ one has 
 $ \exists
X \in \IF^{d \times d},\;  A=BX \;\Rightarrow\;
\range(A) \subseteq \range(B).$
\item[c)]
For $A,B,C\in\IF^{d\times d}$ it holds:
$\range(C)=\range(A)\vee\range(B)
\quad\Leftrightarrow\quad$
\[ \exists W,X,Y,Z\in\IF^{d\times d}: \;
C=A\cdot X+B\cdot Y \;\Band\;
A=C\cdot W \;\Band\; B=C\cdot Z \enspace . \]
\item[d)]
For $A,C\in\IF^{d\times d}$ it holds
\[ \range(C)=\neg\range(A)
\quad\Leftrightarrow\quad
\exists X,Y\in\IF^{d\times d}: \;
A^\adjoint\cdot C=O \;\Band\; A\cdot X+C\cdot Y=I \enspace . \]
\end{enumerate}
Similarly,  for finite dimensional $\IF$-unitary spaces
$\calH$
and their endomorphism $\ast$-rings $\End(\calH)$.
\end{fact} 
\begin{proof}
In a) choose $\IF^d=U\oplus V$ and
let $A(\vec u +\vec v)=\vec u$ for $\vec u \in U$ and $\vec v \in V$.\\
b) This writes the generators of $\range(A)$ as linear combinations
of those of $\range(B)$.  
In c), the first condition in the formula  is equivalent to
the columns of $C$ being linear combinations of those
of $(A|B)$, i.e. to
$\range(C)\subseteq\range(A)\vee\range(B)$; 
 while
in view of b)  the converse
inclusion  is equivalent to the  conjunction of the last two conditions.
The condition in d)   tells $\range(C)$ orthogonal to
$\range(A)$ and joining with it to $\IF^d$. This characterizes the
orthocomplement of $\range(A)$. The  reasoning carries
over to abstract spaces and their endomorphism $\ast$-rings;
in a)--c) we may use any basis, then  in d) there
is no more  need to refer to a basis.      
\end{proof}

\subsection{Relational Interpretations} 
\lab{ss:Interpret}

A particularly useful kind of reduction   is  implicit  in the
proof of the Cook-Levin-Theorem and in its adaptation to the 
Blum-Shub-Smale Model \cite[p.403]{Cucker}: here the evaluation of a
term is decomposed into a polynomial-size, existentially quantified 
system of 
equations together with a simple  term providing the final value:

\begin{myexample} \lab{x:RingInterp}
Fix some commutative ring $\Ring$ with unity.
\\
The evaluation of a $k$-variate polynomial $p\in R[x_0,x_{-1},\ldots,x_{-k+1}]$
decomposes into a series of basic binary operations
$+$, $\times$, and constants (i.e. 0-ary) $c\in R$.
More precisely, a straight-line program $\Gamma$ of length $N$ over $R$
calculating $R^k\ni\bar r\to p(\bar r)\in R$ consists of a
sequence of assignments ``$x_n:=f_n(x_{n_1},\ldots,x_{n_{k_f}})$''
($n=1,\ldots,N$) each
applying a function $f_n$ of arity $k_{f_n}=:k_n$ from $R$'s signature 
to previous intermediate results $x_{n_i}$ ($-k<n_i<n$)
such that $x_N=p(\bar x)$ yields the final result.
That is, for each choice of $r_0,r_{-1},\ldots,r_{-k+1},y\in R$,
the following system of equations in variables
$x_1,\ldots,x_N$ is satisfiable over $R$ 
(and uniquely so) ~iff~ $p(\bar r)=y$ holds:
\begin{equation} \label{e:RelInterp}
y=x_N \quad\Band\quad x_n=f_n(x_{n_1},\ldots,x_{n_{k_n}}), 
\quad n=1,\ldots,N \enspace .
\end{equation}
Conversely if $\Ring$ is ordered, any such system can be
combined into a single equation:
$\Band_{n=1}^N p_n(\bar x)=q_n(\bar x)$ for polynomials $p_n,q_n$
is equivalent to $0=\sum_{n=1}^N \big(p_n(\bar x)-q_n(\bar x)\big)^2$.
In the case of a $\ast$-subfield $\IF\not\subseteq \IR$ of $\IC$ 
use  $0=\sum_{n=1}^N \big(p_n(\bar x)-q_n(\bar x)\big)^\ast\big(p_n(\bar x)-q_n(\bar x)\big)$.
\end{myexample}
This motivates the following.
\begin{mydefinition} \lab{d:RelInterp} (\mycite{\S 2.6.1}{Hodges}) The
 \textsf{relational interpretation} 
$\Psi_{\sigma}$ for a signature $\sigma$ associates with any term $t=t(\bar x)$
a conjunction $\Psi_\sigma(t)$ of  equations with \textsf{output variable} $X_t$  
according to the following recursive definition: 
$\Psi_\sigma(x)$ is $X_x=x$ for any variable $x$;  for any 
 $k$-ary operation symbol $f$   and 
$t=f(t_1,\ldots ,t_k)$ the formula $\Psi_\sigma(t)$ 
is the conjunction of $ X_t=f(z_1, \ldots ,z_k)$
and the $\Psi_\sigma(t_i)(z_i/X_{t_i})$   
with new auxiliary variables $z_i$   
substituted  for the $X_{t_i}$.
 The relational
interpretation $\Psi_\sigma(s=t)$  of an equation $s=t$ is then
the conjunction of $\Psi_\sigma(s)$, $\Psi_{\sigma}(t)$ and
$X_s=X_t$.
\end{mydefinition} 
For example, the relational interpretation of
$ t= \big(x_1 \wedge (x_2 \vee \neg x_1)\big) \vee  (x_2 \vee \neg x_1) $
is 
\[ X_t= z_1 \vee z_2\;\Band\; z_1= x_1 \wedge z_3\;\Band\; z_2= x_2 \vee z_4\;
\Band\;z_3= x_2 \vee z_5\;\Band\; z_4 =\neg x_1\;\Band\; z_5=\neg x_1.\]
\COMMENTED{The  universal variant  yields   here
\[(z_1= x_1 \wedge z_3\;\Band\; z_2= x_2 \vee z_4\;
\Band\;z_3= x_2 \vee z_5\;\Band\; z_4 =\neg x_1\;\Band\; z_5=\neg
x_1 ) \Rightarrow  X_t= z_1 \vee z_2. \]}

\begin{observation} \lab{o:RelInterp}
The relational interpretation $\Psi_\sigma$  
 captures the term $t=t(\bar x)$  within the class of all structures
of  signature $\sigma$: for all 
$a_0$ and $\bar a$ in $\calA$, $a_0=t_A(\bar a)$ 
if and only if $\calA \models \exists\; \Psi_\sigma(t)(a_0;\bar a)$
 where quantification is over all
auxiliary variables.
Given (a suitable encoding of) a system of equations 
$s_i(\bar x)=t_i(\bar x)$ ($1\leq i\leq I)$ 
with terms $s_i,t_i$ of some signature $\sigma$,
their image under the relational interpretation $\Psi_{\sigma}$
can be calculated by a polynomial-time Turing machine 
to produce a system $y_j=f_j(\bar y)$ ($1\leq j\leq J)$ 
of equations of signature $\sigma$ such that 
\begin{itemize}
\item[\textbullet] the original system is satisfiable in a structure $\calS$
  of signature $\sigma$ ~iff~ the new system is
\item[\textbullet] the new system is basic in the sense that each equation
  has a variable symbol on the left and one function symbol on the right.
\end{itemize} 
\end{observation}
In Example~\ref{x:RingInterp}, observing that
the $p_n$ in Equation~(\ref{e:RelInterp})
are just pairwise distinct variables
and the $q_n$ are either linear of quadratic polynomials
(w.l.o.g. with coefficients $0,\pm1$ according to Observation~\ref{o:MultChain})
it follows the well-known

\begin{fact} \lab{f:RingInterp}
For any ordered field $\IF$, $\FEAS_{\IZ,\IF}$ is polynomial-time equivalent
to the question of whether a list of quadratic integer polynomials
(w.l.o.g. with coefficients $0,\pm1$)
has a joint root; which in turn is polynomial-time equivalent 
to that of a single quartic integer polynomial $P(X_1,\ldots,X_N)$
(with coefficients $0,\pm1,\ldots,\pm N$) having a root.
\end{fact}
Indeed, the squared polynomials $(p_n-q_n)^2$
($1\leq n\leq N$)
arising at the end of Example~\ref{x:RingInterp} 
involve in expanded form only coefficients
$0,\pm1,\pm2$; which in $Q=\sum_{n\leq N} (p_n-q_n)^2$
cannot add up beyond $\pm N$.

Similar reformulations in the structure of ortholattices
will be employed extensively in Theorem~\ref{t:Syntax} below.

\begin{myremark} \lab{r:Interp}
Reductions respecting truth had been employed in logic at least
since \cite{Tarski} and have become ubiquitous in complexity theory 
with \cite{Cook,Karp}. For counting problems, they have been
generalized to parsimonious reductions \cite{Valiant}.
And for more general (i.e. not necessarily integer-valued)
function problems $f:S\to T$, e.g. in algebraic complexity theory
\mycite{Remark~4.4}{ACT}, morphisms give rise to a notion of 
reduction---provided that $S$ and $T$ share the same signature.
\end{myremark}
Subsection~\ref{ss:IntField} however `reduces'
$(\IR,+,-,0,\cdot,1)$ to $\big(\Gr(\IR^3),\vee,\wedge,\neg,\Zero,\One\big)$:
clearly structures of very different signatures. 
Moreover the mapping $\Interpret_{\Bar A}:\IR\to\Gr(\IR^3)$ 
depends on the frame $\bar A$.
This can be seen as an instance of
the quite general concept  \emph{interpretation}
regularly employed in Model Theory  \mycite{\S5.5}{Burris} 
and  \mycite{\S 5.3}{Hodges}. 
We have relied in Theorem~\ref{t:IntStarField} 
 on the  famous  \emph{coordinatization} due to \person{von Neumann},
one of many results in Geometry
and Algebra which can be understood as interpretations. 
Another such result  is the  isomorphism between 
$\Gr(\IF^d)$  and the ortholattice of
(principal) right ideals of the matrix $\ast$-ring $\IF^{d \times d}$
--- related to the  use of  Fact~\ref{f:lat}  to prove 
Proposition~\ref{p:UpperComplexity} by reduction to
$\FEAS_{\IZ,\Re\IF}$. 

\cxx{
\begin{digression} \lab{d:Interp}
Such interpretations also  provide translations from sentences
in one language (e.g. fields) into another (e.g. lattices) 
which often  are actually polynomial time reductions.
E.g. the case of fixed $d \geq 3$ 
in   Theorem~\ref{t:PolyHierarchy}  
may be  recovered  from  interpretations yielding 
polynomial time equivalence between any prenex  sentences
in the language of $\ast$-rings w.r.t validity in $\IF$ to such in 
the language of ortholattices  w.r.t validity in $\Gr(\IF^d)$
---  preserving alternation type of quantifications and uniform in
$\IF$.  Here, in view of  Lemma~\ref{l:Boolean2} 
we may assume that  ortholattice  sentences $\phi$ have
matrix of the form $t(\bar x)=\One$. Such
$\phi$ is translated to a prenex sentence in $\ast$-ring language --
refering to $\IF^{d \times d}$, first, and then  
proceeding to $\IF$. Depending on the innermost quantifier of $\phi$
one has to use the existentially respectively
the universally quantified variant of relational
interpretation (cmp. \mycite{\S2.5.1}{Hodges}),
the first combined with
the interpretation given by Fact~\ref{f:lat}, the second 
 with the interpretation of 
 $\Gr(\IF^d)$ as the ortholattice of projections in $\IF^{d \times d}$
 (cmp. Section~\ref{s:Neumann}). 
In the converse direction,
given a  sentence $\psi$  in $\ast$-ring 
language
 proceed as in the proof
 of Propositiom~\ref{p:IntField},    using in addition  
the term  $x^{\dagger_{\bar y}}$ 
of Lemma~\ref{l:MatinGrH} capturing conjugation,
to  translate  atomic subformulas whence  matrix 
of $\psi$  into the matrix of $\psi'$ in  ortholattice language.
 Then, as in the proof of
Theorem~\ref{t:PolyHierarchy}.
 introduce 
the orthonormal $d$-frame $\bar y$  by 
bounded existential  or universal  quantification
put in front of $\psi'$  
and
 convert quantifiers  of $\psi'$ into bounded 
ones   to
encode  the conditions $x \in \IntRing_{\bar y}$;
 pass from there  to prenex form by Lemma~\ref{l:bdquant}
or, more efficiently, avoid bounded quantifications
using  terms
based on \cite{Mayet}, evaluating in $\Gr(\IF^d)$ 
to an orthonormal $d$-frame or to $\Zero$, 
 and an ortholattice term $r(x,\bar y)$  such that
 $r(x,\bar a)\in\IntRing_{\bar a}$ and 
$x = r(x, \bar a) \Leftrightarrow x \in \IntRing_{\bar a}$
for all orthogonal $d$-frames.  
(cmp. \cite{Marina}).

Finally, we mention the particularly useful interpretations
of  the field $\Re \IF$ in the $\ast$-subfield $\IF$ of $\IC$ 
and vice versa. In order to comply with the BSS-machine concept,
we derive most results for real fields, first, and
then use this equivalence to read them as results about real parts.
But, working with interpretations more directly,
one could deal with $\ast$-fields $\IF$ 
(including the real case), primarily, and then derive
the results  involving  real parts $\Re \IF$.   
\end{digression} }

\subsection{Some Tools from Logic and  Orthologic}\lab{ss:tool} 

Recall the \textsf{bounded quantifiers} 
\[(\exists \bar x \alpha(\bar x)).\; \phi(\bar x,\bar y) 
\;:\Leftrightarrow \; \exists \bar x.\; \alpha(\bar x) \Band \phi(\bar
x,\bar y) \quad\text{and}\quad
(\forall \bar x \alpha(\bar x)).\; \phi(\bar x,\bar y) 
\;:\Leftrightarrow \; \forall  \bar x.\; \alpha(\bar x) \impl \phi(\bar
x,\bar y).\]
We use  the old fashioned
$ \phi \impl \psi$ as a shorthand notation for
 logical implication 
$\Bneg \phi \Bor \psi$.

\begin{mylemma} \lab{l:bdquant} 
Let $\phi(\bar x^1, \ldots ,\bar x^\ell, \bar y)$  be a first order
formula 
(all listed variables pairwise distinct) and
$\psi$ given as
\[ \big(\qf_1 \bar x^1 \alpha_1(\bar  x^1)\big) \;
\big(\qf_2 \bar x^2 \alpha_2(\bar x^1,\bar x^2)\big) \;
\big(\qf_3 \bar x^3 \alpha_3(\bar x^1,\bar x^2,\bar x^3)\big) \ldots
\big( \qf_\ell \bar x^\ell \alpha_l(\bar x^1, \ldots ,\bar x^\ell)\big) .\, 
\phi(\bar x^1, \ldots ,\bar x^\ell, \bar y)\]
where $\qf_i \in \{\forall,\exists\}$. Then there is a Boolean combination $\beta$
(to be determined in polynomial time) of the 
$\alpha_i$ and $\phi$ such that $\psi$ is logically equivalent to
\[  \qf_1 \bar x^1  \qf_2 \bar x^2 \ldots \qf_\ell x^\ell.\;
\beta( \bar x^1,  \ldots, \bar x^\ell, \bar y).\] 
\end{mylemma} 
\begin{proof}
  Recall the logical equivalences 
\[ \alpha \Band \qf  \bar x \beta \Leftrightarrow 
\qf \bar x\; (\alpha \Band \beta),\quad
\alpha \impl \qf  \bar x \beta \Leftrightarrow 
\qf  \bar  x \;(\alpha \impl\beta)\]
where $\qf \in \{\forall,\;\exists\}$ and 
no variable form  $\bar x$ does  occur freely
in $\alpha$ - for the easy proof one may assume that the
$\bar x$ are the only 
variables occuring freely.   Proceeding by induction on $\ell$ we have
\[ 
\big(\qf_2 \bar x^2 \alpha_2(\bar x^1,\bar x^2)\big) \;
\big(\qf_3 \bar x^3 \alpha_3(\bar x^1,\bar x^2,\bar x^3)\big) \ldots
\big( \qf_\ell \bar x^\ell \alpha_\ell(\bar x^1, \ldots ,\bar x^\ell)\big) .\, 
\phi(\bar x^1, \ldots ,\bar x^\ell, \bar y)\]
logically equivalent to 
\[ \qf_2 \bar x^2 \qf_3 \bar x^3 \ldots \qf_\ell \bar x^\ell.\;
\chi(\bar x^1, \ldots ,\bar x^\ell, \bar y)\]
where $\chi$ is a boolean combination of $\phi$ and the
$\alpha_2, \ldots, \alpha_\ell$.   Thus,
depending on its first quantifier, $\psi$ is logically
equivalent to  one of
\[ \begin{array}{lcl}
\exists  \bar x^1.\; \alpha_1(\bar x^1) &\;\Band\;& \qf_2 \bar x^2 \qf_3 \bar x^3 \ldots \qf_\ell \bar x^\ell.\;
\chi(\bar x^1, \ldots ,\bar x^\ell, \bar y)\\
\forall  \bar x^1.\; \alpha_1(\bar x^1) &\impl& \qf_2 \bar x^2 \qf_3 \bar x^3 \ldots \qf_\ell \bar x^\ell.\;
\chi(\bar x^1, \ldots ,\bar x^\ell, \bar y) \end{array}\]
and these in turn by the above equivalence rules to
\[ \begin{array}{lcl}
\exists  \bar x^1 \qf_2 \bar x^2 \qf_3 \bar x^3 \ldots \qf_\ell \bar x^\ell.\;
 \alpha_1(\bar x^1) &\;\Band\;&\chi(\bar x^1, \ldots ,\bar x^\ell, \bar y)\\
\forall  \bar x^1  \qf_2 \bar x^2 \qf_3 \bar x^3 \ldots \qf_\ell \bar x^\ell.\;
 \alpha_1(\bar x^1) &\impl&\chi(\bar x^1, \ldots ,\bar x^\ell, \bar
 y) \end{array}\]
having  a form as required.
\end{proof}

The following extends Lemma~\ref{l:2D}b)
and Observation~\ref{o:Product}b): 
\begin{mylemma} \lab{l:2Dx}
\begin{enumerate}
\item[a)]
For a first-order sentence $\phi$ 
in the language of ortholattices
with $n$ (bound) variables,
the following are equivalent: \\
i)~ $L\models \phi$ for some/any infinite MOL $L$ of dimension 2
\qquad ii)~ $\calMO_n\models \phi$ 
\item[b)]
For a term $t(x_1,\ldots,x_n)$ in the language of ortholattices
and MOLs $L,L'$, the sentences $\phi$ and $\chi$ given as
$\Quantifier_1 x_1\Quantifier_2 x_2\ldots\Quantifier_n x_n:t(\bar x)=\One$
and 
$\Quantifier_1 x_1\Quantifier_2 x_2\ldots\Quantifier_n x_n:t(\bar x)\neq\Zero$,
respectively, 
have 
\[ L\models \phi \;\mbox{ and } \; L'\models \phi\;\;\Leftrightarrow\;\; L\times L'\models \phi
\quad\text{ and }\quad
L\models \chi\;\mbox{ or }\; L'\models \chi\;\;\Leftrightarrow\;\; L\times L'\models \chi
\]
\end{enumerate}
\end{mylemma}

\begin{proof}
\begin{longenum}
\item[a)]
Fix an infinite MOL $L$ of dimension $2$ 
and add all its elements $a$ as constants $\underline{a}$  to the language.
For a sentence $\phi$ in this language 
let $n_\phi$ be total number of variables and new constants occurring in
$\phi$. Given $\phi$, consider  the system $\calS_\phi$ of all sub-ortholattices 
$S$ of $L$  of size $|S|=2n_\phi+2$ and
such that $S$ contains all  $a$ such that $\underline{a}$
occurs in $\phi$.  We consider $S \in \cal S_\phi$  in an extended
signature where 
 constants $\underline{a}$ occurring in $\phi$ 
are interpreted as $a \in S$. 
We show by induction on the length of the formula
$\phi$ (with logical connectives $\Bneg$, $\Band$, and  $\exists x$)
 that the following are equivalent 
\begin{quote} $L\models \phi$,\quad $S\models \phi$ for some $S \in \calS_\phi$,
\quad $S\models \phi$ for all $S \in \calS_\phi$.
\end{quote} 
The last equivalence is immediate by the fact that all $S \in \calS_\phi$ 
are isomorphic.  This follows from (the proof of)  Lemma~\ref{l:2D});
as does the remaining equivalence   
 if $\phi$  is an atomic sentence.  In the inductive step consider
$\phi$ given as  $\phi_1 \Band \phi_2$ and choose
$S \in \calS_{\phi_1} \cap \calS_{\phi_2} =\calS_\phi$; then the claim follows, 
readily. Even more so
  $\phi$ being $\Bneg \psi$ where $\calS_\phi=\calS_\psi$.
Finally, let $\phi$ be given as $\exists x.\, \psi(x)$.  If $L \models \phi$ 
then $L \models \psi(\underline{a})$ for some $a \in L$.
Choose $S \in \calS_{\psi(\underline{a})}$ and
apply the inductive hypothesis to conclude
$S \models  \psi(\underline{a})$ whence $S\models \phi$.
Conversely, assume that $S\models \phi$ for some
$S \in \calS_\phi$; then $S\models \psi(\underline{a})$ 
for some $a \in S$; now $S \in  \calS_{\psi(\underline{a})}$
 whence $L \models   \psi(\underline{a})$ 
by inductive hypothesis 
\item[b)]  Again, it is convenient to add constants to the language:
for each $(a,b) \in L\times L'$ let $\underline{(a,b)}$ a constant
interpreted as $(a,b)$ in $L \times L'$, as $a$ in $L$ and
 as $b$ in $L'$. Observe, that $L\times L'$ is the
direct product of $L$ and $L'$ also in the extended signature.
We prove the first equivalence for sentences  $\phi$ involving such constants
by induction on the number $n$ of quantifiers.
For $n=0$ the definition of direct product applies.
Let $n>0$ and $\phi=\qf_1 \psi(x_1)$. 
Then, using  induction,   $L\times L' \models \phi$ iff $L \times L' \models
\psi(\underline{(a,b)})$   (for some respectively all $(a,b)$) 
iff  $L \models   \psi(\underline{(a,b)})$ and
$L' \models   \psi(\underline{(a,b)})$ iff
 $L \models \phi$ and $L' \models \phi$. 
 The second equivalence follows by contraposition
and  converting $\Bneg \phi$ into prenex form,  
having replaced 
equation   $t(\bar x)=\One$ 
by $\neg t(\bar x) =\One$.      
\qed\end{longenum} 
\end{proof} 
%
\subsection{Huhn's Diamonds} \lab{ss:Huhn}

An equivalent concept of frame (more convenient for equational questions) 
 has been established by  \person{Huhn} \cite{Huhn}.
The idea is to  consider instead of the $A_{ij}$ 
just the single $A_0:= A_{12}+ \ldots + A_{1d}$.

From  the $A_i$ on can recover 
$A_{1i}$ as $A_0\wedge  (A_1 \vee A_i)$, $i>0$:
If $\vec x \in A_0\cap (A_1+A_i)$
then $\vec x = x_1 \vec v_1 +x_i \vec v_i
= \sum_{j=2}^d y_j(\vec v_1-\vec v_j)$  requires $y_j=0$ for
each $j\not\in\{1,i\}$.  And from $A_{1i}$ and $A_{1j}$ one
obtains $A_{ij}$ as $(A_i+A_j)\cap (A_{1i}+A_{1j})$. 
The relevant  relations between $A_0$ and the
$A_i$ ($i>0$) are
\begin{equation} \label{eq:frame2}
A_0 \vee A_i =\One,\;\; A_0 \wedge A_i =\Zero \;\;\mbox{ for }
i=1,\ldots ,d \end{equation}
Indeed observe that $\vec v_1 \in A_0+A_i$
whence $\vec v_j \in A_0+A_i$ for all $j$;
and that for  $\vec x = y_i(\vec v_1 - \vec v_i) \in A_k$, 
where $k\in\{1,i\}$, 
it follows $y_i= 0$.  
This motivates the following
\begin{mydefinition}
An  \textsf{orthogonal}  $d$-\textsf{diamond} in a MOL, $L$, is 
  a system $\bar a=(a_i\mid 0\leq i \leq d)$  of elements 
such that 
\[\begin{array}{lcll}  
\One&=& a_1\vee \ldots \vee a_d\\
a_i &\leq& \neg a_j&\mbox{ for all } i,j =1, \ldots ,d,\,i\neq j,   \\
\Zero&=&a_0 \wedge a_i &\mbox{ for all } i=1, \ldots ,d,\\
\One&=& a_0 \vee a_i &\mbox{ for all } i=1, \ldots ,d.
 \end{array}\]
\end{mydefinition} 
Actually, Huhn calls this a $(d-1)$--diamond and
requires  independence of the $a_1, \ldots, a_d$ in place of
orthogonality.

\newcommand{\ov}[1]{\overline{#1}}

\cx{
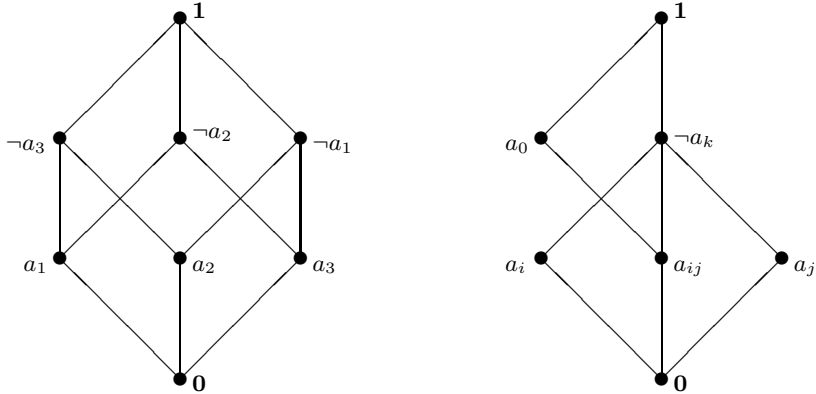
\begin{figure} 
\setlength{\unitlength}{8mm}
\begin{center}
   
\begin{picture}(10,6)  

\put(0,0){\circle*{0.2}} 
\put(0.2,-0,2){$\Zero$}

\put(-2,2){\circle*{0.2}} 
\put(-2.6,1.8){$a_1$} 

\put(0,2){\circle*{0.2}} 
\put(0.2,1.8){$a_2$}

\put(2,2){\circle*{0.2}}
\put(2.2,1.8){$a_3$}
 
\put(-2,4){\circle*{0.2}}
\put(-2.9,3.8){$\neg a_3$}

\put(2,4){\circle*{0.2}}
\put(2.2,3.8){$\neg a_1$}

\put(0,4){\circle*{0.2}}
\put(0.2,4){$\neg a_2$} 
\put(0,6){\circle*{0.2}}
\put(0.2,6){$\One$}

\put(0,0){\line(0,1){2}} 
\put(0,0){\line(1,1){2}} 
\put(0,0){\line(-1,1){2}} 

\put(-2,2){\line(0,1){2}} 
\put(-2,2){\line(1,1){2}}

\put(2,2){\line(0,1){2}} 
\put(2,2){\line(-1,1){2}}

\put(0,2){\line(-1,1){2}} 
\put(0,2){\line(1,1){2}}

\put(-2,4){\line(1,1){2}} 
\put(0,4){\line(0,1){2}} 
\put(2,4){\line(-1,1){2}}

\put(8,0){\circle*{0.2}}
\put(8.2,-0.2){$\Zero$}

\put(6,2){\circle*{0.2}}
\put(5.4,1.8){$a_i$}
\put(8,2){\circle*{0.2}}
\put(8.2,1.8){$a_{ij}$}
\put(10,2){\circle*{0.2}}
\put(10.2,1.8){$a_j$}
\put(6,4){\circle*{0.2}}
\put(5.4,3.8){$a_0$}
\put(8,4){\circle*{0.2}}
\put(8.2,3.9){$\neg a_k$}
\put(8,6){\circle*{0.2}}
\put(8.2,6){$\One$}

\put(8,0){\line(0,1){2}} 
\put(8,0){\line(-1,1){2}} 
\put(8,0){\line(1,1){2}} 

\put(6,2){\line(1,1){2}} 
\put(6,4){\line(1,1){2}} 
\put(8,2){\line(0,1){2}} 
\put(10,2){\line(-1,1){2}} 
\put(8,4){\line(0,1){2}} 
\put(8,2){\line(-1,1){2}} 

\end{picture}
\caption{ \label{f:diamond} Relations of a  $3$-diamond:  
Order diagrams of relevant sublattices} 
\end{center}
\end{figure}}

\begin{myexample} \lab{x:huhn} 
\begin{longenum}
\item[a)] For
a $d$-dimensional $\IF$-unitary space $\calH$, the 
 $(A_0,A_1,\ldots,A_d)$  as above   form  an orthogonal $d$-diamond 
if the basis $\vec v_1, \ldots ,\vec v_d$ is chosen orthogonal
(Fact~\ref{f:GramSchmidt}a).
\item[b)]
In every 2-dimensional MOL $L$, pairwise distinct atoms
$a_0,a_1,\neg a_1=:a_2$ give rise to an orthogonal 2-diamond.
\end{longenum}
\end{myexample} 

\begin{fact} \lab{f:Huhn} 
 $\dim[ a_{0},\One] =\dim a_i$  for all $i>0$.   
The  claims b) and c) of Fact~\ref{f:frame}     remain valid.
\end{fact} 
Consider the following ortholattice terms 
where $\bar y=(y_0,y_1, \ldots ,y_d)$.
\[ \begin{array}{rcl}
h_d(\bar y)&=&\big(\bigvee\nolimits_{i=1}^d (y_i \wedge \bigwedge\nolimits_{i \neq j >0}
\neg  y_j)\big) \wedge \bigwedge\nolimits_{i=1}^d
\big((y_0 \vee y_i)\wedge (\neg y_0 \vee \neg y_i)\big),\\
\tilde{h}_d(x,\ov{y})&=&
h_d(\bar y )\wedge  \hat{h}_d \;\mbox{ where } 
 \hat{h}_d=\bigvee\nolimits_{i,j >0} 
\Big(\big((x \vee \neg y_i) \wedge y_i) \vee y_0\big) \wedge y_j\Big), \\
 g_d(\bar y)&=&
\neg y_0 \wedge  \bigwedge\nolimits_{i=1}^d\big(y_0  \vee (y_i \wedge \bigwedge\nolimits_{i\neq
  j>0} \neg y_j)\big), \\
 \tilde{g}_d(x,\bar y)&=&  g_d(\bar y) \wedge  \bigwedge\nolimits_{i=1}^d
 \big( y_0\vee (x \wedge y_i)\big) 
\end{array} \]

\begin{mylemma} \lab{l:Huhn} 
 Let $L$ be any MOL, $1 \leq \dim(L) \leq d$; 
In a), b) and d) let $\dim(L)=d$.   
\begin{enumerate} 
\item[a)] $\bar a \in L^{d+1}$ is an
orthogonal $d$-diamond of $L$ if and only if $h_d(\bar a )=\One$.
\item[b)] Given an orthogonal $d$-diamond $\bar a$ and $b$ in $L$, one has 
$b>\Zero$ if and only if $\hat{h}_d(b, \bar a)=\One$.
In particular, if $L$ admits an orthogonal $d$-diamond,
then $b>\Zero$ if and only if $\tilde{h}_d(b,\ov{y})$
is strongly satisfiable in $L$.
\item[c)]  $\bar a \in L^{d+1}$ is
 an  orthogonal $d$-diamond of $L$  and $\dim(L) =d$ if and only if
$g_d(\bar a )>\Zero$. 
\item[d)] Given an orthogonal $d$-diamond $\bar a$ and $b$ in $L$, one has 
$b=\One$ if and only if $\tilde{g}_d(b, \bar a)>\Zero$.
In  particular, if   $L$ admits an orthogonal $d$-diamond,  
then $b=\One$ if and only if  $\tilde{g}_d(b,\ov{y})$
is weakly satisfiable in $L$. 
\end{enumerate}    
\end{mylemma} 
\mycite{Theorem~1}{Hagge1} and \cite{Hagge2} have presented
terms $t_d$ equivalent to $0$ over $\Gr(\IC^{d-1})$
but not over $\Gr(\IC^d)$. A careful analysis of their construction
reveals the length $|t_d|$ to be exponential in $d$. Based on
Huhn's diamonds, we improve that in the following

\begin{corollary} \lab{c:Hagge}
The term $g_d$ has length quadratic in $d$.
It is equivalent to $\Zero$ over $\Gr(\IF^{d-1})$ 
but not over $\Gr(\IF^d)$.  
\end{corollary}
Indeed, a lattice of height $<d$ does not admit an orthogonal
$d$-diamond whereas, by Example~\ref{x:huhn}a), $\Gr(\IF^d)$ does.

\begin{proof}[of Lemma~\ref{l:Huhn}] \begin{longenum}
\item[a)]
Observe that for any $i>0$ one has   $ a_i\leq\neg a_j$
for all $0<j\neq i$ if and only if
$a_i=a_i\wedge\bigwedge_{0< j\neq i}\neg a_j$.
\item[b)] 
Given an orthogonal
$d$-diamond $\bar a $ we have  $h_d(\bar a) =\One$. 
Consider $b>\Zero$. Then  $b \not\leq a_i^\bot $ for
some $i>0$ (since $\bigcap_{i=1}^n a_i^\bot =\Zero$)
and    $b+a^\bot _i=\One$ since $\dim[a_i^\bot ,\One]=1$. 
It follows $(b+a^\bot _i)\cap a_i+a_0 =a_i+a_0 =\One$ 
whence $((b+a_i^\bot )\cap a_i +a_0)\cap a_j=a_j$ for all $j>0$ and $
\tilde{h}_d(b,\bar a )=\One$. 
Conversely, assume $\tilde{h}_d(b,\bar a )=\One$.
Then $h_d(\bar a )=\One$ and   $\bar a $ is an orthogonal $d$-diamond.
Assuming  $b=\Zero$  one has
 $(b+a^\bot _i)\cap a_i =a^\bot _i\cap a_i=\Zero$ and
$((a+a_i^\bot )\cap a_i +a_0)\cap a_j=a_0\cap a_j=\Zero$ whence
$\tilde{h}_d(b,\bar a )=\Zero$, as contradiction. 
\item[c)]
Assume $g_d(\bar a )>\Zero$ and put 
$b_i:=a_i \cap  \bigcap_{i\neq j>0} a_j^\bot$ for $i>0$.
Assuming $b_i=\Zero$ we get
$g_d(\bar a ) \leq a_0^\bot \cap (a_0+b_i) =a_0\cap a_0^\bot =\Zero$, a contradiction.   
Thus $b_i >\Zero$.  Also, by definition, $b_i \leq a_i$ and 
$b_i \leq a_j^\bot \leq  b_j^\bot$ for $i\neq j>0$,
i.e. the $b_i$, are pairwise orthogonal. Thus
$d \leq \sum_{i=1}^d \dim b_i = \dim \sum_{i=1}^d b_i \leq d$
whence $\dim b_i =1$, $\dim(L)=d$,  and $\sum_{i=1}^d b_i=\One$.
From $\Zero<a_i \leq b_i$ it follows $b_i=a_i$ for $i>0$.
Assuming $a_i\leq a_0$ for some $i>0$ 
gives, $g_d(\bar a) \leq a_0^\bot \cap (a_0+a_i) =a_0 \cap a_0^\bot =\Zero$,
a contradiction. Thus, $a_0\cap a_i=\Zero$  and
  $\dim[a_0,\,a_0+a_i]=1$ for all $i>0$.
Assuming $a_0+a_i\neq a_0+a_j$ for some $i\neq j$,
$i,j>0$ gives $a_0=(a_0+a_i)\cap (a_0+a_j)$ whence
$g_d(\bar a ) \leq a_0^\bot\cap (a_0+a_i)\cap (a_0+a_j) =\Zero$,
a contradiction.  It follows $a_0+a_i =a_0+a_j$ 
for all  $i\neq j$,
$i,j>0$  whence $a_0+a_i \geq \sum_{j=1}^d a_j =\One$. 
\\
Conversely, given an orthogonal
$d$-diamond $\bar a$ one calculates from the  relations
that $g_d(\bar a) =a_0^\bot$. Assuming $a_0^\bot=\Zero$ 
would imply $a_0=\One$, $a_i=\One$ for all $i>0$, and $\One=\Zero$: contradiction.
\item[d)]
If $b=\One$  substitute into $\bar y$  an  orthogonal $d$-diamond
$\bar a $. Conversely, assume
 $g_d(\bar a ) \cap  \bigcap_{i=1}^d
 ( a_0+b\cap a_i) >\Zero$ for some $\bar a $.  
Then $g_d(\bar a )>\Zero$ and $\bar a $ is an orthogonal
$d$-diamond   by c).  
Assume $b<\One$. Then $b \not\geq a_i$ for some $i>0$,
whence $b \cap a_i=\Zero$ and it follows 
\[g_d(\bar a )\cap \bigcap\nolimits_{i=1}^d
 ( a_0+ b\cap a_i)  \leq a_0^\bot\cap  (a_0 +b \cap a_i) =a_0^\bot\cap
 a_0=\Zero \]
a contradiction. Thus, $b=\One$.
\qed\end{longenum}
\end{proof} 
%
\subsection{Expressing Boolean within First-Order Quantum Logic}
\lab{ss:Quantifiers}
The following  results are used in Subsection~\ref{ss:PolyHierarchy},
only, but are of some interest by themselves.
Recall that a quantifier free formula 
$\phi(\bar x)$  in the first order language of ortholattices
is a Boolean combination of equations $s_i(\bar x)=t_i(\bar x) $ 
of equations between ortholattice terms $s_i(\bar x) ,t_i\bar x)$. 
Recall, that to form Boolean combinations within Logic  we 
use  connectives $\Bneg, \Band, \Bor$  and
$\impl$   for `not, `and', `or'
and `implies'.  We have to carefully distinguish these
from ortholattice terms built from $\neg, \wedge$, and $\vee$.

The following Lemma
shows that within fixed dimension $d$  and in the 
presence of orthogonal $d$-frames any such
formula $\phi(\bar x)$ can be equivalently expressed
in the  form  $\qf \bar y. t(\bar x,\bar y)=1$  
with an  ortholattice term $t(\bar x,\bar y)$ 
and quantifier  $\qf$ which may be chosen both
as $\exists$ and as $\forall$.  In particular,
considering atomic formulae of the form
$x =\One$ and $x \neq \One$ (i.e.$\Bneg (x=\One)$)
this shows that the Boolean (propositional)  theory of
strong truth can be expressed within quantified modular
Quantum Logic  of fixed dimension.

  Following \mycite{\S15 Theorem 3}{Kalmbach} we define the 
ortholattice terms
\[ x \rightarrow y := y \vee \neg(x \vee y),\quad 
x \leftrightarrow y := (x \wedge y) \vee \neg(x\vee y)   \]
and recall (Example~\ref{x:QL}f) that in any MOL 
\[ a \rightarrow b =\One    \;\Leftrightarrow \; a \leq b,\quad 
a \leftrightarrow  b =\One    \;\Leftrightarrow \; a = b.\]
\begin{mylemma}  \lab{Mayet}
 For any $d\geq 2$  there is an
ortholattice term $f_d(\bar  y)$, $\bar y=(y_0, \ldots, y_d)$,
such that for any  MOL $L$ of $\dim(L)=d$ 
and any $\bar a$ in $L$
\begin{enumerate}
\item[a)]$f_d(\bar a)=\One$ if and only if  $\bar a$ is an
orthogonal $d$-diamond;  $f_d(\bar a)=\Zero$, otherwise. 
\item[b)]$f_d(\bar a) =\One \Rightarrow x=\One$ 
if and only if $f_d(\bar a) \rightarrow x =\One$. 
\end{enumerate} 
\end{mylemma} 
\begin{proof} 
Define
\[ f_d(\bar y) = \bigvee\nolimits_{i=1}^d \big((g_d(\bar y) \vee \neg y_i) \wedge y_i\big).\]   
If $\bar a$ is an orthogonal $d$-diamond, then 
$g_d(\bar a) =\neg a_0$ and $\neg a_0 \vee \neg a_i=\One$   
for $i>0$ since $a_0 \wedge a_i= \Zero$. It follows
$f_d(\bar a) =\bigvee\nolimits_{i=1}^d a_i=\One$.
Conversely, assume that $\bar a$ is not an orthogonal $d$-diamond.
Then $g_d(\bar a)=\Zero$ from Lemma~\ref{l:Huhn}c)
and its follows $f_d(\bar a)=  
\bigvee\nolimits_{i=1}^d \neg a_i \wedge a_i =\Zero$.
This proves  a) and b) follows immediately. 
\COMMENTED{
{\bf Sichere Variante}:
 In \cite{Mayet} (cmp. \cite{Marina})  there are provided
terms $t_i(\bar y)$ ($i=0, \ldots ,d$),
such that  for any $\bar a$ in an MOL $L$,
the $t_9(\bar a), \ldots ,t_d(\bar a)$ 
form an orthogonal $d$-frame  in some interval $[\Zero,u]$ of $L$.
By Fact~\ref{f:Huhn} we conclude that either
$u=\Zero$ or $u=\One$ if $\dim(L)=d$.   Now, define 
$f_d(\bar y)= \bigvee\nolimits_{i=1}^d t_i(\bar y)$. }
\end{proof}

\begin{mylemma} \lab{l:Boolean2}  
Fix $d \geq 2$.
For every quantifier free formula $\phi(\bar x)$
in the first order language of ortholattices  there 
are ortholattice terms $t'=t'_{d,\phi}(\bar x, \bar y)$ and  
$t''=t''_{d,\phi}(\bar x, \bar y)$ 
computable in polynomial time from $\phi$ such that
all modular ortholattices $L$, of $\dim(L)=d$ and
admitting an orthogonal $d$-diamond, and for all   
for all $\bar b$ in $L$
\[ 
L\models  \exists \bar y:\;t'_{d\phi}(\bar b, \bar y) =\One
\; \Leftrightarrow\;  L\models \phi(\bar b) \; \Leftrightarrow\;
L\models  \forall \bar y:\;t''_{d\phi}(\bar b, \bar y)=\One .\]
More precisely, we have $\bar y=(y_0, \ldots ,y_d)$ 
and ortholattice terms $t_{d \phi}(\bar x, \bar y)$ 
 such that 
\[ 
L\models  \exists \bar y:\;f_d(\bar y) \wedge t_{d\phi}(\bar b, \bar y)=\One 
\; \Leftrightarrow\; L\models \phi(\bar b)\; \Leftrightarrow\;
L\models  \forall \bar y:\;f_d(\bar y) \rightarrow   t_{d\phi}(\bar b, \bar y)=\One \]
for all $\bar b$ in $L$ --  where  $f_d$  is the term from
Lemma~\ref{Mayet}.   Moreover
\begin{enumerate} 
\item $t_{d\phi}(\bar x, \bar y) 
= s_1(\bar x) \leftrightarrow t_1(\bar x)$ for $\phi$ an equation
 $s_1(\bar x) =t_1(\bar x)$.
\item   $t_{d\phi}(\bar x, \bar y) =
t_{d\phi_1}(\bar x, \bar y) \wedge t_{d\phi_2} (\bar x, \bar y)$ for
$\phi$ given as $\phi_1 \Band \phi_2$. 
\item   $t_{d\phi}(\bar x, \bar y) =\hat{h}_d(\neg t_{d\psi}(\bar x,
  \bar y),\bar y)$
for $\phi= \Bneg \psi$ and term $\hat{h}_d(x,\bar y)$ from
Lemma~\ref{l:Huhn}.  
 \end{enumerate} 
\end{mylemma} 
\begin{proof} 
For atomic formulas $\phi$, i.e.
equations $s_1=t_1$ we refer to  Example~\ref{x:QL}f). 
Now, it suffices to 
proceed with structural induction  for formulas built
with    propositional connectives
$\Band$ and $\Bneg$ and to verify that (2) and (3) do the job for
these.

 Considering (2),
 assume that
$L\models \phi(\bar  b)$, i.e. that $L \models \phi_i(\bar  b)$ for $i=1,2$.
Due to inductive  hypothesis on the  second equivalence  we have
$L \models t_{d \phi_i}(\bar b,\bar a)=\One$ 
for any orthogonal $d$-diamond $\bar a$,
whence $L \models t_{d \phi}(\bar b,\bar a)=\One$ 
and $L \models \forall \bar y.\;f_d(\bar y) \rightarrow 
  t_{d \phi}(\bar b,\bar a)=\One$. 
Due to inductive hypothesis on the first equivalence,
there $\bar a_i$ such that
$L\models  f_d(\bar a_i) \wedge t_{d\phi_i}(\bar b, \bar a_i)=\One$;
in particular, both $\bar a_i$ are orthogonal $d$-diamonds.
As just observed, also 
$L\models  t_{d\phi_2}(\bar b, \bar a_1)=\One$ whence 
$L\models  \exists \bar y:\;f_d(\bar y) \wedge t_{d\phi}(\bar b, \bar
y)=\One$. 

Conversely, assume  $L\models \exists \bar y.\; f_d(\bar y) \wedge
t_{d\phi}(\bar b,\bar y)=\One$.  Choose a witnessing
(orthogonal $d$-diamond) $\bar a$, conclude 
$L\models t_{d\phi_i}(\bar b,\bar y)=\One$ and,
by induction,
$L \models \phi_i(\bar b)$. Now, assume 
 $L\models \forall \bar y.\; f_d(\bar y) \rightarrow
t_{d\phi}(\bar b,\bar y)=\One$. Choose any orthogonal
$d$-diamond $\bar a$  (we supposed that such exist) and
conclude $L\models t_{d \phi_i}(\bar b, \bar a)$ and,
by induction, $L\models \phi_i(\bar b)$ for $i=1,2$.

Coming to (3) let $\phi=\Bneg \psi$. Observe that, 
by Lemma~\ref{l:Huhn}, for any orthogonal $d$-diamond $\bar a$ 
one has 
\[t_{d\psi}(\bar b,\bar a) \neq \One \;\Leftrightarrow\; 
\neg t_{d\psi}(\bar b,\bar a) \neq \Zero\;\Leftrightarrow\; 
\hat{h}_d(\neg t_{d\phi}(\bar b,\bar a),\bar a)  =\One.\]
Now, assume that $L\models \phi(\bar b)$, i.e. that
$L \not\models \psi(\bar b)$. By inductive hypothesis,
for all orthogonal $d$-diamonds (and there are such) one has
$t_{d\psi}(\bar b,\bar a) \neq \One$. With the above
observation we derive the other two sentences in (3)
to be valid in $L$.  Conversely
if for  some, in particular if  for all, orthogonal
$d$-diamonds $\bar a$  one has
 $\hat{h}_d(\neg t_{d\phi}(\bar b,\bar a),\bar a)  =\One$ 
then
$t_{d\psi}(\bar b,\bar a) \neq \One$ and $L\not\models
\phi(\bar b)$ by induction. 
\end{proof}

\cx{
\begin{digression} \lab{d:Boolean} 
 Adding the constants of an orthogonal $d$-diamond  to a height $d$ MOL, $L$,
gives rise to a discriminator term on $L$, see \cite{Micol}.
This is the general reason behind Lemma~\ref{l:Boolean2}. 
For any $d\geq 3$,   Lemmas~\ref{l:IntStarField} 
and \ref{l:Boolean2} associate
with any 
quantifier free  formula $\psi(\bar x)$ in the language of
$\ast$-rings  an ortholattice term
$s_{d \psi}(\bar x,\, \bar y)$  such that
$\psi(\bar r)$ holds in $\IF$ if and only if
$s_{d \psi}(\Interpret_{\bar a}r _1, \ldots \Interpret_{\bar a}r_n,\,  \bar a) $ holds in
$\Gr(\IF^d)$ where  $\bar a$ is any orthonormal $d$-frame
 of $\Gr(\IF^d)$.
Cmp. \cite{Mnev2} for related results. 
\end{digression}}

\section{Variations of the Quantum Satisfiability Problem}
\lab{s:WeakStrong}
This section establishes (Theorem~\ref{t:WeakStrong}) the
polynomial-time equivalence of $\sat_L$ and $\SAT_L$
for any fixed modular ortholattice $L$ of finite height.
We investigate the complexity of satisfiability for
terms with a limited number of $\bigwedge\bigvee$-alternations
and for negation-free terms (Subsection~\ref{ss:Syntax}).
The effect of changing the ground field 
and the dimension is explored
in Subsections~\ref{ss:Fields} and \ref{ss:Dimensions}, respectively.
Generalizing satisfiability, Subsection~\ref{ss:PolyHierarchy}
determines the complexity of deciding the truth of quantified formulas.

\subsection{Weak versus Strong Satisfiability}\lab{ss:WeakStrong}  
Based on Lemma~\ref{l:Huhn} and Example~\ref{x:huhn} we can now conclude

\begin{theorem} \lab{t:WeakStrong}
Fix $d\in\IN$. Then for any modular ortholattice $L$ of height $d$
admitting an orthogonal $d$--diamond
and for any ortholattice term $t(\bar x)$, the following hold:
\begin{enumerate}
\item[a)] $t(\bar x)$ is weakly satisfiable in $L$ ~if and only if~
$\tilde{h}_d\big(t(\bar x),\bar y\big)$ is strongly satisfiable in $L$.
\item[b)] $t(\bar x)$ is strongly satisfiable in $L$ ~if and only if~
$\tilde{g}_d\big(t(\bar x),\bar y\big)$ is weakly satisfiable in $L$.
\end{enumerate} 
In particular, weak and strong satisfiability over Hilbert quantum
logics are mutually polynomial-time reducible: 
For every finite-dimensional unitary $\IF$-vector space $\calH$ it holds
$\sat_{\Gr(\calH)}\reduceq\SAT_{\Gr(\calH)}\reduceq\sat_{\Gr(\calH)}$.
\end{theorem}
Indeed, the translations indicated in a) and b) 
can easily be computed in time polynomial in 
the length of $t$ (and in $d$, otherwise independent of $L$).
Theorem~\ref{t:IntStarField}b) thus extends to

\begin{corollary} \lab{c:WeakStrong}
For every $\IF$-unitary vector space $\calH$ of dimension
$d\geq3$ admitting an equinormal orthogonal basis,
$\sat_{\Gr(\calH)}$ is $\calBP(\calNP_{\Re\IF}^0)$--complete.
\end{corollary}

\begin{comment}
\begin{longenum} \item[a)]
One can prove Theorem~\ref{t:WeakStrong}a) using the transitivity
of the orthogonal group: map $t$ to 
$\bigvee_{i=1}^d t_i$ where  $t_i$  arises from $t$,
replacing each variable $x$ in $t$ by  a variable
$x_i$  occurring only in $t_i$. Also,
in Theorem~\ref{t:WeakStrong}b)  one can use terms provided by Hagge \cite{Hagge2} 
and an iterative construction. 
\item[b)] 
According to Comment~\ref{mol}, any  MOL  of  height $d$
is  isomorphic to a direct product of  irreducibles
of heights $d_1, \ldots, d_k$ such that $d=d_1+\ldots +d_k$. 
And an MOL of height $d_i$  is irreducible if and only if it 
admits an orthogonal $d_i$-diamond. 
Since in Lemma~\ref{l:Huhn} the only requirement was the existence of such,
the equivalence stated in  Theorem~\ref{t:WeakStrong}  can be established for any
given MOL, $L$,  of finite height $d\geq 1$ 
via $\bigvee_i \phi_{d_i}$ and $\bigwedge_i \psi_{d_i}$.
 Though, there are $L$ where both problems are undecidable
and the same holds if one considers satisfiability
 within some (indefinite)  MOL of height $\leq d$ 
(resp. irreducible of height $d$) where $d \geq 14$ is fixed
(cmp. Digression~\ref{mol}).   
\end{longenum} 
\end{comment} 
%

\subsection{Quantified Quantum Propositions} \lab{ss:PolyHierarchy}
Stockmeyer's 
polynomial hierarchy for Turing machines is based on,
and extends, the complexity classes $\calP$ and $\calNP$ 
as well as the class $\co\calNP$ of 
complements of $\calNP$-problems. 
More specifically, starting with $\PolyS{0}=\calP=\PolyP{0}$
and $\PolyS{1}=\calNP$ and $\PolyP{1}=\co\calNP$,
higher classes $\PolyS{\ell}$ and $\PolyP{\ell}$ ($\ell\geq2$)
can be characterized both syntactically and semantically:
For the latter, $\PolyS{\ell}$ contains precisely those decision problems
accepted by nondeterministic polynomial-time Turing machines
with oracle access to some $V\in \PolyS{\ell-1}$ 
(equivalently: to some $V'\in \PolyP{\ell-1}$);
and $\PolyP{\ell}$ consists of the complements 
of members from $\PolyS{\ell}$ \cite{Papadimitriou}.
For the former characterization, $\PolyS{\ell}$ and $\PolyP{\ell}$ 
contain all decision problems of the form
\begin{align}
\big\{ \bar z\in\{\Zero,\One\}^n \;\big|\; n\in\IN, \;
\exists \bar y^{(1)}\in\{\Zero,\One\}^{p(n)} &\;
\forall \bar y^{(2)}\in\{\Zero,\One\}^{p(n)} \;
\exists \bar y^{(3)}\in\{\Zero,\One\}^{p(n)} \;
\forall \bar y^{(4)}\ldots   \label{e:PolyS}
\\[-0.5ex] &\ldots \nonumber
\Quantifier_\ell \bar y^{(\ell)}\in\{\Zero,\One\}^{p(n)} :\;
\langle\bar z,\bar y^{(1)},\ldots\bar y^{(\ell)}\rangle\in V\big\} \\
\big\{ \bar z\in\{\Zero,\One\}^n \;\big|\; n\in\IN, \;
\forall \bar y^{(1)}\in\{\Zero,\One\}^{p(n)} &\;
\exists \bar y^{(2)}\in\{\Zero,\One\}^{p(n)} \;
\forall \bar y^{(3)}\in\{\Zero,\One\}^{p(n)} \;
\exists \bar y^{(4)}\ldots \nonumber
\\[-0.5ex] &\ldots \nonumber
\Quantifier'_\ell \bar y^{(\ell)}\in\{\Zero,\One\}^{p(n)} :\;
\langle\bar z,\bar y^{(1)},\ldots\bar y^{(\ell)}\rangle\in V'\big\}
\end{align}
respectively, with $V,V'\in\calP$ and $p\in\IN[N]$ a polynomial.
Here $\Quantifier_\ell$ denotes the existential quantifier
when $\ell$ is odd and otherwise the universal one;
vice versa for $\Quantifier'$.

Generalizing the Cook-Levin Theorem,
a natural problem complete for $\PolyS{\ell}$
asks for the truth of a given Boolean formula with $\ell$
blocks of alternating quantifiers, starting with the existential one:
\begin{multline*}
\SAT^{\ell} \;=\;
\big\{ \langle t(\bar x^{(1)},\bar x^{(2)},\ldots
\bar x^{(\ell)}) \rangle \;\big|\;
n_1,\ldots n_\ell\in\IN, \;
\exists \bar a^{(1)}\in\{0,1\}^{n_1} \;
\forall \bar a^{(2)}\in\{0,1\}^{n_2} \;
\exists \bar a^{(3)}
\\ \ldots
\Quantifier_\ell \bar a^{(\ell)}\in\{0,1\}^{n_\ell}:
\; t(\bar a^{(1)},\bar a^{(2)},\ldots\bar a^{(\ell)})=1 \big\} \enspace .
\end{multline*}
So binarily encoded terms $\langle t(\bar x^{(1)},\bar x^{(2)},\ldots,
\bar x^{(\ell)}) \rangle$ replace $\bar z$ in Equation~(\ref{e:PolyS}).
$\UNSAT^{\ell}$ is defined similarly but starting 
with the universal quantifier---and complete for $\PolyP{\ell}$.
Moreover, the following problem $\QSAT$ is complete for $\PSPACE$:
\[
\big\{ \langle t(x_1,\ldots x_n)\rangle \;\big|\; n\in\IN\;
\exists a_1\in\{0,1\} \;\forall a_2\in\{0,1\}\;
\exists a_3 \ldots \Quantifier_n a_n\in\{0,1\}:
\; t(a_1,\ldots a_n)=1 \big\} \]
More generally, sequential polynomial space
corresponds to 
parallel polynomial time \cite{Borodin};
as well as to 
parallel alternating time \cite{Alternation}
hence $\PSPACE$ is sometimes also denoted as $\PAR=\PAT$. 

Both the polynomial hierarchy and its two characterizations translate
(although with notably different proofs) to the \BSS setting
\mycite{\S4}{Cucker}, \mycite{\S21.4}{BCSS}:
$\BssPolyS{\ell,\IF}$ contains precisely those decision problems
accepted by nondeterministic polynomial-time \BSS machines
over $\IF$ with oracle access to some $\IV\in\BssPolyS{\ell-1,\IF}$ 
(equivalently: to some $\IV'\in\BssPolyP{\ell-1,\IF}$);
and $\BssPolyP{\ell,\IF}$ consists of the complements 
of members from $\BssPolyS{\ell,\IF}$.
For the former characterization, $\BssPolyS{\ell,\IF}$ and 
$\BssPolyP{\ell,\IF}$ consist of all sets of the form
\begin{gather*}
\big\{ \bar z\in\IF^n \;\big|\; n\in\IN, \;
\exists \bar y^{(1)}\in\IF^{p(n)} \;
\forall \bar y^{(2)}\in\IF^{p(n)} \;
\exists \bar y^{(3)} \ldots
\Quantifier_\ell \bar y^{(\ell)}\in\IF^{p(n)} :
\langle\bar z,\bar y^{(1)},\ldots\bar y^{(\ell)}\rangle\in\IV\big\} \\
\big\{ \bar z\in\IF^n \;\big|\; n\in\IN, \;
\forall \bar y^{(1)}\in\IF^{p(n)} \;
\exists \bar y^{(2)}\in\IF^{p(n)} \;
\forall \bar y^{(3)}\ldots
\Quantifier'_\ell \bar y^{(\ell)}\in\IF^{p(n)} :
\langle\bar z,\bar y^{(1)},\ldots\bar y^{(\ell)}\rangle\in\IV'\big\}
\end{gather*}
respectively, with $\IV,\IV'\in\calP_\IF$ and $p\in\IN[N]$;
cmp. also \cite{Bournez}. And 
\begin{align*}
\FEAS^{2\ell-1}_{\IF,\IF} = & \big\{
 \langle p_1,\ldots p_k\rangle \;\big|\;
k,n_1,\ldots\in\IN, \; p_j\in\IF[\bar X^{(1)},\ldots \bar X^{(2\ell-1)}], \;
\exists\bar y^{(1)}\in\IF^{n_1}\;
\forall\bar y^{(2)}\in\IF^{n_2}\;
\\[-0.5ex] 
\exists\bar y^{(3)} \; \ldots &
\exists \bar y^{(2\ell-1)}\in\IF^{n_{2\ell-1}} :
p_1(\bar y^{(1)},\ldots \bar y^{(2\ell-1)})=0\;\Band\cdots\Band\;
p_k(\bar y^{(1)},\ldots \bar y^{(2\ell-1)})=0 \big\} \\
\FEAS^{2\ell}_{\IF,\IF} = & \big\{
\langle p_1,\ldots p_k\rangle \;\big|\;
k,n_1,\ldots\in\IN, \; p_j\in\IF[\bar X^{(1)},\ldots \bar X^{(2\ell)}], \;
\exists\bar y^{(1)}\in\IF^{n_1}\;
\forall\bar y^{(2)}\in\IF^{n_2}\;
\\[-0.5ex] 
\exists\bar y^{(3)} \; \ldots &
\forall\bar y^{(2\ell)}\in\IF^{n_{2\ell}} :
p_1(\bar y^{(1)},\ldots \bar y^{(2\ell)})\neq0\;\Bor\cdots\Bor\;
p_k(\bar y^{(1)},\ldots \bar y^{(2\ell)})\neq0 \big\} \\
\FEAS^{2\ell-1}_{\IZ,\IF} = & \big\{
\langle p_1,\ldots p_k\rangle \;\big|\;
k,n_1,\ldots\in\IN, \; p_j\in\IZ[\bar X^{(1)},\ldots \bar X^{(2\ell-1)}], \;
\exists\bar y^{(1)}\in\IF^{n_1}\;
\forall\bar y^{(2)}\in\IF^{n_2}\;
\\[-0.5ex] 
\exists\bar y^{(3)} \ldots &
\exists \bar y^{(2\ell-1)}\in\IF^{n_{2\ell-1}} :
p_1(\bar y^{(1)},\ldots \bar y^{(2\ell-1)})=0\;\Band\cdots\Band\;
p_k(\bar y^{(1)},\ldots \bar y^{(2\ell-1)})=0 \big\} \\
\FEAS^{2\ell}_{\IZ,\IF} = & \big\{
\langle p_1,\ldots p_k\rangle \;\big|\;
k,n_1,\ldots\in\IN, \; p_j\in\IZ[\bar X^{(1)},\ldots \bar X^{(2\ell)}], \;
\exists\bar y^{(1)}\in\IF^{n_1}\;
\forall\bar y^{(2)}\in\IF^{n_2}\;
\\[-0.5ex] 
\exists\bar y^{(3)} \ldots &
\forall \bar y^{(2\ell)}\in\IF^{n_{2\ell}} :
p_1(\bar y^{(1)},\ldots \bar y^{(2\ell)})\neq0\;\Bor\cdots\Bor\;
p_k(\bar y^{(1)},\ldots \bar y^{(2\ell)})\neq0 \big\} 
\end{align*}
are complete for 
$\BssPolyS{2\ell-1,\IF}$,
$\BssPolyS{2\ell,\IF}$,
$\calBP(\BssPolyS{2\ell-1,\IF})$, 
and $\calBP(\BssPolyS{2\ell,\IF})$, 
respectively.
Note that (Fact~\ref{f:Constructible}) the matrix form defining
$\FEAS_{\IF,\IF}^{\ell}$ and $\FEAS_{\IZ,\IF}^{\ell}$
can be relaxed to arbitrary finite Boolean combinations
of, and in the case admitting an order (bottom of Example~\ref{x:RingInterp})
also restricted to one single, polynomial equality or inequality:
depending on $\ell$'s parity!

Space complexity does not translate as
nicely to the \BSS setting \cite{Michaux,Briquel};
however $\PAT_{\IF}$ is a natural counterpart to $\PAT=\PSPACE$
and $\QSAT_{\IF}$ complete for it, where $\PAT_{\IF}$ consists 
of all subsets of $\IF^*$ of the form
\begin{align*}
\big\{ \bar z\in\IF^n \;\big|\; n\in\IN, \;
\exists y_1\in\IF &\;\forall y_2\in\IF \;
\exists y_3\in\IF \;\forall y_4\ldots\; 
\ldots\Quantifier_n y_n\in\IF :
\langle\bar z,\bar y\rangle\in\IV\big\} \qquad\text{ and } \\
\QSAT_{\IF} \;:=\; \big\{ \langle p_1,\ldots,p_k\rangle \;\big|\;&
k,n\in\IN, \; p_j\in\IF[X_1,\ldots,X_n], \; \\[-0.5ex]
\exists y_1\in\IF \; &\forall y_2\in\IF \;
\exists y_3\in\IF \; \forall y_4\ldots
\ldots\Quantifier_n y_n\in\IF :
p_1(\bar y)=\cdots=p_k(\bar y)=0  \big\}
\end{align*}
with $\IV$ running through $\calP_{\IF}$
\mycite{Theorem~4.1}{Cucker};
similarly for $\QSAT^0_{\IF}$
complete for $\calBP(\PAT_{\IF})$.
Natural problems in $\PAR_{\IC}\subset\PAT_{\IC}$ 
and $\PAR_{\IR}\subset\PAT_{\IR}$ traditionally arise
in semi-/algebraic geometry 
\cite{Canny,Giusti,Lecerf,Jeronimo,Scheiblechner1,BssToda,Scheiblechner2}.

\bigskip
In view of our generalization of Boolean satisfiability to
ortholattices $L$ (Definition~\ref{d:QL}f)
this suggests to consider first-order quantified quantum 
(i.e. predicate) logic and to define

\begin{eqnarray*}
\SAT^\ell_L &:=& \big\{\langle t(\bar x^{(1)},\ldots \bar x^{(\ell)})\rangle
\;\big|\;n_1,\ldots n_\ell\in\IN,\;
\exists \bar a^{(1)}\in L^{n_1} \;
\forall \bar a^{(2)}\in L^{n_2} \;
\exists \bar a^{(3)}
\ldots \\[-0.5ex] && \qquad\qquad\qquad \qquad\qquad\qquad \ldots 
\Quantifier_\ell \bar a^{(\ell)}\in L^{n_\ell}:
\; t_L(\bar a^{(1)},\bar a^{(2)},\ldots \bar a^{(\ell)})=\One \big\}  
\enspace, \\ 
\UNSAT^\ell_L &:=& \big\{\langle t(\bar x^{(1)},\ldots \bar x^{(\ell)})\rangle
\;\big|\;n_1,\ldots n_\ell\in\IN,\;
\forall \bar a^{(1)}\in L^{n_1} \;
\exists \bar a^{(2)}\in L^{n_2} \;
\forall \bar a^{(3)}
\ldots \\[-0.5ex] && \qquad\qquad\qquad \qquad\qquad\qquad \ldots 
\Quantifier'_\ell \bar a^{(\ell)}\in L^{n_\ell}:
\; t_L(\bar a^{(1)},\bar a^{(2)},\ldots \bar a^{(\ell)})\neq\One \big\}  
\enspace, \\ 
\sat^\ell_L &:=& \big\{
\langle t(\bar x^{(1)},\ldots \bar x^{(\ell)})\rangle
\;\big|\;n_1,\ldots n_\ell\in\IN,\;
\exists \bar a^{(1)}\in L^{n_1} \;
\forall \bar a^{(2)}\in L^{n_2} \;
\exists \bar a^{(3)}
\ldots \\[-0.5ex] && \qquad\qquad\qquad \qquad\qquad\qquad \ldots 
\Quantifier_\ell \bar a^{(\ell)}\in L^{n_\ell}:
\; t_L(\bar a^{(1)},\bar a^{(2)},\ldots \bar a^{(\ell)})\neq\Zero \big\} 
\enspace, \\ 
\unsat^\ell_L &:=& \big\{\langle t(\bar x^{(1)},\ldots \bar x^{(\ell)})\rangle
\;\big|\;n_1,\ldots n_\ell\in\IN,\;
\forall \bar a^{(1)}\in L^{n_1} \;
\exists \bar a^{(2)}\in L^{n_2} \;
\forall \bar a^{(3)}
\ldots \\[-0.5ex] && \qquad\qquad\qquad \qquad\qquad\qquad \ldots 
\Quantifier'_\ell \bar a^{(\ell)}\in L^{n_\ell}:
\; t_L(\bar a^{(1)},\bar a^{(2)},\ldots \bar a^{(\ell)})=\Zero \big\}  
\enspace,
\\ 
\QSAT_L &:=& \big\{\langle t(x_1,\ldots x_n)\rangle
\;\big|\; 
\exists a_1\in L \; \forall a_2\in L \;
\exists a_3\in L \; \forall a_4 \;\ldots
\Quantifier_n a_n\in L: t_L(\bar a)=\One \big\}
\enspace, \\ 
\qsat_L &:=& \big\{\langle t(x_1,\ldots x_n)\rangle
\;\big|\; 
\exists a_1\in L \; \forall a_2\in L \;
\exists a_3\in L \; \forall a_4 \;\ldots
\Quantifier_n a_n\in L:
\; t_L(\bar a)\neq\Zero \big\}  
\end{eqnarray*}
and analogously for terms \emph{with constants};
cmp. also \cite{Roman}.
Observe that in  all of $\QSAT_{\IF}$, $\QSAT_L$,
and $\qsat_L$  dummy variables are admitted;
in the latter two such  variable $x$  can be camouflaged by
meeting with $(x\vee \neg x)$.

We emphasize that 
the definitions of $\SAT^\ell_L,\sat^\ell_L$ do
\emph{not} depend on the parity of $\ell$.

\begin{theorem} \lab{t:PolyHierarchy}
Fix $\ell\in\IN$.
\begin{enumerate}
\item[a)]
For any infinite two-dimensional MOL $L$,
both $\SAT^\ell_L$ and $\sat^\ell_L$ are
$\PolyS{\ell}$--complete;
and $\UNSAT^\ell_L,\unsat^\ell_L$ are
$\PolyP{\ell}$--complete;
while $\QSAT_L$ and $\qsat_L$ are 
$\PSPACE$--complete.
\item[b)] 
For $\IF\subseteq\IR$ and $d\geq3$,
both $\SAT^\ell_{\Gr(\IF^d)}$ and $\sat^\ell_{\Gr(\IF^d)}$ are
$\calBP(\BssPolyS{\ell,\IF})$--complete;
and $\UNSAT^\ell_{\Gr(\IF^d)}$ and $\unsat^\ell_{\Gr(\IF^d)}$ are
$\calBP(\BssPolyP{\ell,\IF})$--complete;
while $\QSAT_{\Gr(\IF^d)}$ and $\qsat_{\Gr(\IF^d)}$ are
$\calBP(\PAT_\IF)$--complete.
\item[c)] 
For every $\IF$-unitary vector space $\calH$ of
dimension $d\geq3$ admitting an equinormal orthogonal basis,
both $\SAT^\ell_{\Gr(\calH)}$ and $\sat^\ell_{\Gr(\calH)}$ are
$\calBP(\BssPolyS{\ell,\Re\IF})$--complete;
and $\UNSAT^\ell_{\Gr(\calH)}$ and $\unsat^\ell_{\Gr(\calH)}$ are
$\calBP(\BssPolyP{\ell,\Re\IF})$--complete;
while $\QSAT_{\Gr(\calH)}$ and $\qsat_{\Gr(\calH)}$ are
$\calBP(\PAT_{\Re\IF})$--complete.
\item[d)]
Similarly, for $\calH=\IF^d$ with $d\geq3$, 
the above problems \emph{with constants} 
$\SAT^\ell_{\Gr(\calH),\Gr(\calH)}$ and $\sat^\ell_{\Gr(\calH),\Gr(\calH)}$ 
are $\BssPolyS{\ell,\Re\IF}$--complete;
and $\UNSAT^\ell_{\Gr(\calH),\Gr(\calH)}$ and $\unsat^\ell_{\Gr(\calH),\Gr(\calH)}$ 
are $\BssPolyP{\ell,\Re\IF}$--complete;
while $\QSAT_{\Gr(\calH),\Gr(\calH)}$ and $\qsat_{\Gr(\calH),\Gr(\calH)}$ 
are $\PAT_{\Re\IF}$--complete. Here,  
lattice constants are considered
encoded as (ranges of) matrices as in Theorem~\ref{t:IntStarField}c)
\end{enumerate}
\end{theorem}
Since \cite{Poonen} proved $\FEAS^{3}_\IQ$ 
undecidable, we conclude 

\begin{corollary} \lab{c:Poonen}
$\SAT^3_{\Gr(\IQ^3)}$ and $\sat^3_{\Gr(\IQ^3)}$ 
are undecidable 
(to a Turing machine\footnotemark[1]).
\end{corollary}
It should be emphasized that (to the best of our knowledge)
Theorem~\ref{t:PolyHierarchy} cannot be just deduced from 
the satisfiability case $\ell=1$ but requires 
additional considerations such as Lemmas~\ref{l:Boolean2}
and \ref{l:2Dx}.

\begin{proof}[of Theorem~\ref{t:PolyHierarchy}]
We first prove 
$\SAT^\ell_L$  polynomial-time equivalent to $\sat^\ell_L$:
$\exists \bar{a}^{(1)} \forall \bar{a}^{(2)} \ldots : 
t_L(\bar{a}^{(1)},\bar{a}^{(2)},  \ldots )\neq\Zero$
is equivalent within $L$ to
$\exists \big(\bar b,\bar{a}^{(1)}\big),\forall \bar{a}^{(2)} \ldots :
\tilde{h}_d\big(t_L(\bar{a}^{(1)},\bar{a}^{(2)},\ldots  ),\bar b\big)=\One$
according to Lemma~\ref{l:Huhn}b) and Example~\ref{x:huhn},
thus yielding a polynomial-time reduction from
$\sat^\ell_L$ to $\SAT^\ell_L$.
The converse reduction similarly follows from Lemma~\ref{l:Huhn}d).
$\QSAT_L$ and $\qsat_L$ 
are seen polynomial-time equivalent
analogously.
\begin{longenum}
\item[a)]
In view of the above preconsiderations
it suffices to treat the cases $\SAT^\ell_L$
and $\QSAT_L$. Concerning the first, the
case $\ell=1$ amounts to satisfiability and was
already shown complete for $\calNP=\PolyS{1}$ in
Theorem~\ref{t:2DNPc}. To locate
\[ \SAT^2_L \;=\;
\big\{\langle t(x_1,\ldots,x_k,y_1,\ldots,y_m)\rangle
\;\big|\;
\exists \bar a\in L^k \;
\forall \bar b\in L^m : 
\; t_L(\bar a,\bar b)=\One \big\} \]
in $\PolyS{2}$
observe that, for a given $t(\bar x,\bar y)$
in $k+m$ variables and 
according to Lemma~\ref{l:2Dx}a), 
it suffices to consider the case $L=\calMO_{k+m}$.
Now elements $z\in\calMO_{k+m}$
can be encoded as in Proposition~\ref{p:2DinNP}a)
as integers in $\{0,1,\ldots,2k+2m\}$ which, in
unary, correspond to $\bar z\in\{0,1\}^{2(k+m)}$.
Tuples $(a_1,\ldots,a_k)\in\calMO_{k+m}^k$ 
and $(b_1,\ldots,b_m)\in\calMO_{k+m}^m$ 
thus can be encoded
as $(\bar a_1,\ldots,\bar a_k)\in\{0,1\}^{p(k+m)}$
and $(\bar b_1,\ldots,\bar b_m)\in\{0,1\}^{p(k+m)}$
of size $p(n):=4n^2$ polynomial in the 
binary length $n\geq k+m$ of the $(k+m)$-variate input.
This rephrases
$\SAT^2_L$ in the form of Equation~(\ref{e:PolyS})
with 
\begin{multline*}
V \;:=\;\big\{\big\langle t(x_1,\ldots,x_k,y_1,\ldots,y_m),
(\bar a_1,\ldots,\bar a_k),(\bar b_1,\ldots,\bar b_m)\big\rangle\;\big|\; \\
k,m\in\IN, \; t_{\calMO_{k+m}}(a_1,\ldots,a_k,b_1,\ldots,b_m)=\One \big\} \enspace
\end{multline*}
a problem in $\calP$ according to Proposition~\ref{p:2DinNP}a).
More generally regarding $\SAT^\ell_L$,
each $k_j$-tuple $\bar a^{(j)}$
in $\calMO_n$ with $n\geq\sum_{j=1}^\ell k_j$
can be encoded as an $\leq(k_j\cdot n)$-tuple in $\{0,1\}$;
similarly for $\QSAT_L$.

\COMMENTED{ doch keine idee
Concerning $\PolyS{\ell}$--hardness of $\SAT^\ell_L$,
we reduce from $\SAT^\ell$ and 
adapt the proof of Proposition~\ref{p:NPhard}.
Observe that $L$ admits an orthogonal $2$-diamond $\bar c$
(Example~\ref{x:huhn}).
Given such, we require all variables in  $t$ from $\SAT^\ell$,
to commute with $c_i$ (and $c_2$), i.e.  to
take  values only in the four element Boolean MOL
$\{\Zero,\One, c_1,c_2\}$.  Thus, we obtain the sentence
\[  \exists (y_0,y_1,y_2)    \exists \bar{x}^1 \forall \bar{y} ^2 \ldots
 \qf \bar x^\ell.\; h_2(\bar y) \wedge  \bigwedge_i C(x,y_1) 
\wedge t(\bar x^1, \ldots ) =\One \]}

Concerning $\PolyS{\ell}$--hardness of $\SAT^\ell_L$,
we reduce from $\SAT^\ell$ and recall (Fact~\ref{f:Foulis}f)
that pairwise commuting elements `live' in a Boolean 
subalgebra of $L$, i.e. in $\{0,1\}^2$.
Let us therefore abbreviate 
\[ 
\bar C(x_1,\ldots x_n):=\bigwedge\limits_{i<j} C(x_i,x_j) 
\quad\text{ and }\quad
\bar C\big((x_1,\ldots x_n),(y_1,\ldots y_m)\big)
\;:=\;\bigwedge\limits_{{i\leq n,j\leq m}} C(x_i,y_j)
\]
expressing in terms of strong truth
that all $\bar x$ to commute and commute with all $\bar y$,
respectively. Observe that 
$\bar C(\bar x,\bar y\bar z)=\bar C(\bar x,\bar y)\wedge\bar C(\bar x,\bar z)$
and $\bar C(\bar x\bar y)=\bar C(\bar x)\wedge\bar C(\bar y)\wedge\bar C(\bar x,\bar y)$.
We now adapt the proof of Proposition~\ref{p:NPhard}
and want to require all quantifiers to range over
pairwise commuting elements.
In the existential case 
$\exists a_1,\ldots,a_n:t(\bar a)=\One$ this was accomplished by 
proceeding to $t'(\bar a):=t(\bar a)\wedge\bar C(\bar a,\bar a)$.
More generally, $\SAT^\ell$--instance
\[ 
\exists\bar a^{(1)}\; 
\forall\bar a^{(2)}\;
\exists\bar a^{(3)}\;
\ldots
\Quantifier_\ell\bar a^{(\ell)}:
\; t\big(\bar a^{(1)},\bar a^{(2)},\ldots\bar a^{(\ell)}\big)=1 \quad ?
\]
over $\{\Zero,\One\}$ is equivalent to the following formula
with bounded quantifiers over $L$:
\begin{multline} \label{e:PolyHierarchy2}
\big(\exists\bar b^{(1)}.\bar C(\bar b^{(1)})=\One\big)\; 
\big(\forall\bar b^{(2)}.\bar C(\bar b^{(2)},\bar b^{(1)}\bar b^{(2)})=\One\big)\;
\big(\exists\bar b^{(3)}.\bar C(\bar b^{(3)},\bar b^{(1)}\bar b^{(2)}\bar b^{(3)})=\One\big)\;
\ldots \\[-0.5ex]
\ldots \big(\Quantifier_\ell\bar b^{(\ell)}
.\bar C(\bar b^{(\ell)},\bar b^{(1)}\bar b^{(2)}\cdots\bar b^{(\ell)})=\One\big)
:\; t\big(\bar b^{(1)},\bar b^{(2)},\ldots\bar b^{(\ell)}\big)=\One
\end{multline}
Indeed, all quantification in (\ref{e:PolyHierarchy2}) 
is bounded to mutually commuting elements and thus `live'
in $\{\Zero,\One\}^2$; hence Lemma~\ref{l:2Dx}b) applies.
Now according to Lemma~\ref{l:bdquant},
(\ref{e:PolyHierarchy2}) is in turn equivalent over $L$ to the unboundedly quantified
\begin{multline*}
\exists\bar b^{(1)}\;
\forall\bar b^{(2)}\;
\exists\bar b^{(3)}\;
\ldots
\Quantifier_\ell\bar b^{(\ell)}:
\; \beta\Big(t'\big(\bar b^{(1)},\bar b^{(2)},\ldots\bar b^{(\ell)}\big)=\One,
\bar C(\bar b^{(1)})=\One,\bar C(\bar b^{(2)},\bar b^{(1)}\bar b^{(2)})=\One, \\[-0.5ex]
\bar C(\bar b^{(3)},\bar b^{(1)}\bar b^{(2)}\bar b^{(3)})=\One,\ldots
\bar C(\bar b^{(\ell)},\bar b^{(1)}\cdots\bar b^{(\ell)})=\One\Big)
\end{multline*}
with Boolean propositional formula $\beta$ polynomial time computable;
which finally can be converted in similar time 
into the $\SAT^\ell_L$--instance
\[ \exists\bar b^{(1)}\; 
\forall\bar b^{(2)}\;
\exists\bar b^{(3)}\;
\ldots
\Quantifier_\ell \big(\bar b^{(\ell)},\bar x\big):
\; t''\big(\bar b^{(1)},\bar b^{(2)},\ldots,\bar a^{(4)},\bar x\big)=\One \quad ? \]
with $t''$ according to Lemma~\ref{l:Boolean2}.
The reduction from $\QSAT$ to $\QSAT_L$ proceeds analogously.
\item[b)]
Again, we first consider the case $\ell=2$:
Encode $a_j\in\Gr(\IF^d)$ as in Proposition~\ref{p:UpperComplexity}
as any matrix $A_j\in\IF^{d\times d}$ with $a_j=\range(A_j)$.
Then 
\begin{gather*}
\langle t(x_1,\ldots,x_k,y_1,\ldots,y_m)\rangle\in\SAT^2_{\Gr(\IF^d)}
\quad\Leftrightarrow\quad
\exists \bar A\in\IF^{kd^2} \;
\forall \bar B\in\IF^{md^2}: \;
\langle t(\bar x,\bar y),\bar A,\bar B\rangle\in\IV, \\
\IV \;:=\;
\big\{\langle t(x_1,\ldots,x_k,y_1,\ldots,y_m),\bar A,\bar B\rangle \;\big|\;
A_i,B_j\in\IF^{d\times d}, 
t_{\Gr(\IF^d)}(\range\bar A,\range\bar B)=\One
\big\} \end{gather*}
with $\IV\in\calP_{\IF}$ according to Proposition~\ref{p:UpperComplexity}:
hence indeed $\SAT^2_{\Gr(\IF^d)}\in\calBP(\BssPolyS{2,\IF})$.
The case of general $\ell$ proceeds similarly;
and so does $\QSAT_{\Gr(\IF^d)}\in\calBP(\PAT_{\IF})$.

In order to show $\SAT^\ell_{\Gr(\IF^d)}$ hard for
$\calBP(\BssPolyS{\ell,\IF})$ we reduce from $\FEAS_{\IZ,\IF}^{\ell}$.
To this end recall 
the proof of Proposition~\ref{p:IntField},
based on  the isomorphisms
$\Interpret_{\bar A}:\IF\to\IntRing_{\bar A}\subseteq\Gr(\IF^d)$
induced  by any $d$-frame $\bar A$ (and the fact that such exist). 
Recall step  (v), the term  $t(\bar z)$  
such that $t(\bar A)=\One$ if and only if
$\bar A$  constitutes a
$d$-frame, as well 
as
steps i) to iii) 
translating polynomials $p_j\in\IZ[X_1,\ldots,X_n]$
into orthoterms $t_j(\bar A;x_1,\ldots,x_n)$
such that $t_j(\bar A;\Interpret_{\bar A}(r_1),\ldots,
\Interpret_{\bar a}(r_n))=\One$
is equivalent to $p_j(r_1, \ldots,r_n)=0$  for any substitution of
$r_i \in \IF$.
Recall  restriction (iv) for  all  variables $x_i$ to range over
the coordinate ring $\IntRing_{\bar A}=\{U\in\Gr(\IF^d):\rho_{\bar A}(U)\}$.
Now, due to the isomorphisms $\Interpret_{\bar A}$,
a $\FEAS_{\IZ,\IF}^{\ell}$--instance
\[ \exists\bar y^{(1)}\; 
\forall\bar y^{(2)}\;
\exists\bar y^{(3)}\;
\ldots
\Quantifier_\ell\bar y^{(\ell)}:
\; \alpha\big(p_1(\bar y^{(1)},\ldots,\bar y^{(\ell)})=0,
\cdots,p_k(\bar y^{(1)},\ldots,\bar y^{(\ell)})\big) \quad ? \]
with arbitrary Boolean propositional formula $\alpha$
is equivalent to
\begin{multline*}
(\exists\bar z)\;
(\exists\bar y^{(1)}\in\IntRing_{\bar z})\;
(\forall\bar y^{(2)}\in\IntRing_{\bar z})\;(\exists\bar y^{(3)}\in\IntRing_{\bar z})\;
 \ldots
(\Quantifier_\ell\bar y^{(\ell)}\in\IntRing_{\bar z}):  \\[-0.5ex]
t(\bar z)=\One\; \Band\; \alpha\Big(
t_{p_1}\big(\bar z,\bar y^{(1)},\ldots\bar y^{(\ell)}\big)=\One,
\cdots,t_{p_k}\big(\bar z,\bar y^{(1)},\ldots\bar y^{(\ell)}\big)=\One\Big) \quad ? 
\end{multline*}
with 
quantifiers $\qf_i \bar y^{(i)}$  
restricted by   $\bar y^{(i)} \in \IntRing_{\bar z}$, the latter   defined 
as $\Band_h \rho_{\bar z}(y^{(i)}_h)$  
where $\rho_{\bar z}(y)$ is the
formula $(y \wedge z_{22} =\Zero\; \Band\; y \vee
z_{22}= z_{11} \vee z_{22})$ and equivalent to
certain  $t_{\rho}(\bar z,y)=\One$ (by Fact~\ref{f:clear}). 
As in a), these restrictions to the quantifiers 
can be moved, according to Lemma~\ref{l:bdquant},
to the matrix, which thus  becomes  a Boolean combination $\alpha'$ of
$t(\bar z)=\One$, 
$t_{p_j}\big(\bar z,\bar y^{(1)},\ldots\bar y^{(\ell)}\big)=\One$,
and $\bigwedge_h t_{\rho}(\bar z,y_h^{(i)})=\One$;
and, based on the 
$\qf_\ell$-version of Lemma~\ref{l:Boolean2}, the whole sentence can be
transformed 
furtheron to a  form  representing a $\SAT^\ell_{\Gr(\IF^d)}$--instance.
The reduction from $\calBP(\PAT_{\IF})$ to $\QSAT_{\Gr(\IF^d)}$
proceeds analogously.
\item[c)]
Similarly to b), now with 
\[ \IV\;:=\;
\big\{\langle t(x_1,\ldots x_k,y_1,\ldots y_m),
\Re\bar A,\Im\bar A,\Re\bar B,\Im\bar B\rangle 
\;\big|\;
t_{\Gr(\IF^d)}(\range\bar A,\range\bar B)=\One
\big\} \]
$\in\calP_{\Re\IF}$
according to Proposition~\ref{p:UpperComplexity}
to show $\in\calBP(\BssPolyS{2,\Re\IF})$;
and, concerning hardness, with 
modifications  as in the proof of Theorem~\ref{t:IntStarField}a):
conditions added for $\bar A$ to constitute an orthonormal 
$d$-frame and quantifiers to range over the real subring
$\IntRing_0:=\{y\in\IntRing:y=y^*\}$ expressible as
$y \otimes_{\bar z} ((z_{11}\vee a_{22}) \wedge \neg y)=\ominus_{\bar
  z} z_{12}$     by the
 second condition in Lemma~\ref{l:IntStarField}b);
the latter   is translated into certain $t'_{\rho}(\bar z ,y)=\One$  
via relational interpretation (Definition~\ref{d:RelInterp})  of this equation 
and Fact~\ref{f:clear}.

\item[d)]
Similarly to c), now with $\bar A$ fixed to the orthonormal
standard frame and polynomial coefficients $c\in\IF$
encoded as 
lattice constants $C:=\Interpret_{\bar A}(c)$;
conversely, lattice constants are encoded by constant matrices,
i.e. lists of numbers in $\IF$; 
recall the proof of Theorem~\ref{t:IntStarField}c).
\qed\end{longenum}\end{proof}
%

\subsection{Syntactically Restricted Terms} \lab{ss:Syntax}
It is well-known that the Boolean satisfiability problem becomes
no more simple when restricted to terms in conjunctive form
($\bigwedge\bigvee$-terms) and even with at most three literals
per clause 
a problem known a \textsf{3SAT}---whereas \textsf{2SAT}
can be decided in polynomial time;
cmp. e.g. \mycite{\S 9.2}{Papadimitriou}.
\cite{Schaefer1} has succeeded in closely delineating syntactically
the border between $\calP$ and $\calNP$-completeness for
Boolean satisfiability problems; see \cite{Chen} for a modern 
presentation from the general perspective of 
\emph{constraint satisfaction}.
Regarding quantum logic, however,
conjunctive form is semantically a proper restriction:

\begin{myexample} \lab{x:Conjunctive}
$(X\wedge Y)\vee(X\wedge\neg Y)$ 
is a term in disjunctive form
not equivalent over $\Gr(\IF^2)$ 
to any term in conjunctive form.
\end{myexample}
Thus more than two alternations
of $\bigvee$ and $\bigwedge$ 
are required to obtain reasonable
syntactical restrictions of the satisfiability problem;
and Theorem~\ref{t:Syntax} below 
explores the boundary between $\calP$ 
and \BSS-$\calNP$ in terms of the number
of these alternations. 

A different kind of syntactic restriction, 
consider terms without negation: the operation
only a real but no complex \BSS machine can 
compute on $\Gr(\IC^d)$, 
recall Proposition~\ref{p:UpperComplexity}a\,iii).
Now (so-called \emph{lattice})
terms $t$ over $\vee$ and $\wedge$ only are monotone
(Lemma~\ref{l:Conjunctive}e);
hence the question of whether a given lattice term
admits an assignment making it evaluate to $\One$
(or $\Zero$) is easy to decide: by checking the
assignment $(\One,\ldots,\One)$ (or $(\Zero,\ldots,\Zero)$).
The alternation classes $\calQ$ of lattice terms are defined,
inductively: variables are in $\bigwedge$ as well as $\bigvee$;
if $t_1,t_2$  are in $\bigwedge\calQ$, then so is $t_1 \wedge t_2$ 
but $t_1 \vee t_2$ is in $\bigvee \bigwedge \calQ$. 
Dually, for $\bigvee\calQ$.  An ortholattice term is
in the alternation class $\calQ$ if by de Morgan's rules
it reduces to a lattice term  (in variables and negated variables) 
 in the class $\calQ$ (cmp. Lemma~\ref{l:Conjunctive}d).

\begin{theorem} \lab{t:Syntax}
\begin{enumerate}
\item[a)]
Strong satisfiability over $\Gr(\IF^d)$, $d\geq2$,
of $\bigwedge\bigvee$-terms is
(independent of $\IF$ and depending
only on $d$'s parity and in this sense uniformly)
polynomial-time decidable. 
\item[b)]
Strong satisfiability over $L$ of general orthoterms 
is polynomial-time reducible to the strong satisfiability
over $L$ of $\bigwedge\bigvee\bigwedge\bigvee$-terms,
uniformly in $L$. 
\item[c)]
For every $\IF$ and $d\geq3$, the language
\quad $\displaystyle
\big\{\langle t(x_1,\ldots,x_n),s(x_1,\ldots,x_n)\rangle \mid $
\abovedisplayskip=0pt \[
\exists a_1,\ldots,a_n\in\Gr(\IF^d) :
t(a_1,\ldots,a_n)=\One \;\Band\; s(a_1,\ldots,a_n)=\Zero \big\}
\;\;\subseteq\;\;\{0,1\}^* \enspace , \]
that is the  question of whether
two given lattice terms $t$ and $s$ admit a \emph{joint}
assignment over $\Gr(\IF^d)$ making $t$ evaluate 
to $\One$ \emph{and} $s$ to $\Zero$, is complete
for $\calBP(\calNP_\IF)$.
\\
The analogous question for 
terms $(t,s)$ \emph{with} constants is $\calNP^0_\IF$--complete.
\item[d)]
The languages from c) remain $\calBP(\calNP_\IF)$--complete
and $\calNP^0_\IF$--complete, respectively,
when restricting to $t \in \bigwedge \bigvee \bigwedge$ and
 $s \in \bigvee \bigwedge \bigvee$.
\item[e)]
For $\IF$ and $d\geq3$, strong satisfiability over
$\Gr(\IF^d)$ of $\bigwedge\bigvee\bigwedge$-terms
is $\calBP(\calNP_{\Re\IF})$--complete; for $d\leq 2$
it is $\calNP$--complete. 
\end{enumerate}
\end{theorem}
In view of the gap between a), b), and e) we ask
\begin{myquestion} \lab{q:Syntax}
What is the computational complexity of the strong satisfiability
problem for $\bigvee\bigwedge\bigvee$-terms? Does it
 depend on dimension or  the ground field? 
\end{myquestion} 
The proof of Theorem~\ref{t:Syntax} is deferred
until later in this subsection
as it relies on some tools.
Recall that the atoms $a_1,\ldots,a_n$ of $\calMO_n$ satisfy
\begin{equation} \label{e:IndepEvenD}
\One \;=\; a_k\vee a_\ell \;=\; 
a_k\vee\neg a_\ell\;=\;
\neg a_k\vee a_\ell\;=\;\neg a_k\vee\neg a_\ell 
\quad (1\leq k<\ell\leq n) \enspace . 
\end{equation}
Lemma~\ref{l:2D} yields such elements also
in $\Gr(\IF^2)$ and, by Observation~\ref{o:Product}a), 
also in $\Gr(\IF^{2d})$ for every 
$\IF$ and $d$. But Equation~(\ref{e:IndepEvenD}) 
cannot be satisfied in odd dimensions. 
Instead, consider the following 

\begin{myexample} \lab{x:IndepOddD}
For every $d\geq2$ and $\IF$ and $n\in\IN$, 
there exist $a_1,\ldots,a_n\in\Gr(\IF^d)$ with
\[ \One \;=\; a_k\vee\neg a_\ell\;=\;\neg a_k\vee\neg a_\ell
\;=\; a_k\vee a_\ell\vee a_j \quad (k\neq\ell\neq j\neq k) \enspace . \]
Indeed, the case of even $d$ has been treated above.
Whereas the $(3+2d)$-dimensional case follows 
from Observation~\ref{o:Product}a) by combining
a 3D instance $(a_1,\ldots,a_n)$ with an even-dimensional one
$(b_1,\ldots,b_n)$ to $(a_1\obot b_1,\ldots,a_n\obot b_n)$.
In the remaining case $d=3$ observe that any three distinct 
vectors $\vec v_k:=(1,k,k^2)\in\IF^3$
form a Vandermonde matrix and are thus linearly independent;
hence $a_k:=\IF\vec v_k\in\Gr_1(\IF^3)$ satisfy 
$a_k\vee a_\ell\vee a_j=\One$ for every $k\neq\ell\neq j\neq k$.
Moreover $\langle\vec v_k,\vec v_\ell\rangle=1+k\ell+k^2\ell^2>0$
shows $\vec v_k\not\!\,\perp\vec v_\ell$, hence 
$\dim(a_k\vee\neg a_\ell)=\dim(a_k)+3-\dim(a_\ell)-\dim(a_k\wedge\neg a_\ell)
=1+3-1$ by the dimension formula (Fact~\ref{f:mod}e).
Similarly, $\neg a_k\neq\neg a_\ell$ and $\dim(\neg a_k)=2=\dim(\neg a_\ell)$
imply $\dim(\neg a_k\vee\neg a_\ell)\geq 3$.
\end{myexample}

\begin{mylemma} \lab{l:Conjunctive}
\begin{enumerate}
\item[a)] Let $t(x_1,\ldots,x_n)$ denote a 
$\bigwedge\bigvee$-term with at least two different
literals in each clause. Then $t$ is strongly 
satisfiable over $\Gr(\IF^{2d})$ for every $\IF$ 
and every $d$.
\item[b)] Let $t(x_1,\ldots,x_n)$ denote a 
$\bigwedge\bigvee$-term with at least three different
literals in each clause. Then $t$ is strongly satisfiable
over $\Gr(\IF^{2d+1})$ for every $\IF$ and every $d\geq1$. \\
More precisely for $a_1,\ldots,a_n\in\Gr_d(\IF^{2d+1})$
according to Example~\ref{x:IndepOddD},
every choice of $b_j\in\{a_j,\neg a_j\}$ gives rise to
a strongly satisfying assignment $(b_1,\ldots,b_n)$ of $t$.
\item[c)] For a $\bigwedge\bigvee$-term $t(x_1,\ldots,x_n)$
with exactly two different literals in each clause,
the following are equivalent:
\begin{enumerate}
\item[i)] $t$ is strongly satisfiable over $\Gr(\IF^d)$ 
 for some $\IF$ and some odd $d$.
\item[ii)] 
$t$ is strongly satisfiable over $\{\Zero,\One\}$.
\item[iii)]
For $a_1,\ldots,a_n\in\Gr_d(\IF^{2d+1})$ 
from Example~\ref{x:IndepOddD},
there exist $b_j\in\{a_j,\neg a_j\}$ such that $\bar b$
constitutes a strongly satisfying assignment of $t$.
\item[iv)] $t$ is strongly satisfiable over $\Gr(\IF^d)$ for
  every $\IF$ and every odd $d$.
\end{enumerate}
\item[d)]
For any ortholattice term $t(\bar x)$ in $n$ variables,
there exists a lattice term $\hat{t}(\bar x,\bar y)$ in $2n$ variables
such that $L\models \hat{t}(x_1,\ldots,x_n,\neg x_1,\ldots,\neg x_n)=t(x_1,\ldots,x_n)$
holds over any ortholattice $L$.
Moreover, $\hat{t}$ can be calculated from $t$ by a Turing machine in polynomial time.
\item[e)]
For any lattice term $t$ and lattice $L$,
$a_i \leq b_i$ ($1\leq i\leq n)$ implies
$t_L(\bar a) \leq   t_L(\bar b)$.
\end{enumerate}
\end{mylemma}
An orthoterm $\hat{t}(x_1,\ldots,x_n,\neg x_1,\ldots,\neg x_n)$
as in d) is called in \emph{negation normal form},
since it has negations only in front of variables.

\begin{proof}\begin{longenum}
\item[a)]
As argued before Example~\ref{x:IndepOddD}, $\Gr(\IF^{2d})$
contains $a_1,\ldots,a_n$ satisfying Equation~\ref{e:IndepEvenD},
i.e. rendering \true every clause $u\vee v$ with distinct
$u,v\in\{a_1,\neg a_1,\ldots,a_n,\neg a_n\}$.
\item[b)]
Similarly observe that, for any three distinct 
$u,v,w\in\{a_j,\neg a_j,a_k,\neg a_k,a_\ell,\neg a_\ell\}$,
$u\vee v\vee w$ evaluates to $\One$.
\item[c\,i$\Rightarrow$ii)] 
Let $\bar y$ denote a strongly satisfying assignment over $\Gr(\IF^{2d+1})$.
We claim that the derived assignment $b_j:=\One$ for $\dim(y_j)\geq d+1$
and $b_j:=\Zero$ for $\dim(\neg y_j)\geq d+1$ is also a satisfying one.
To this end consider a arbitrary clause $u\vee v$ of $t$ with literals
$u,v\in\{x_1,\neg x_1,\ldots,x_n,\neg x_n\}$.
By hypothesis it evaluates to $\One$ when plugging in $\bar y$ for $\bar x$,
requiring $2d+1\leq\dim(u[\bar y])+\dim(v[\bar y])$
and thus that at least one of $u,v$ had been assigned 
a subspace of dimension $\geq d+1$: which in the
derived assignment becomes $1$ and keeps the clause \true.
\item[c\,ii$\Rightarrow$iii)] 
Let $\bar y$ denote a satisfying assignment over $\{\Zero,\One\}$.
We claim that the derived assignment $b_j:=\neg a_j$ for $y_j=\One$
and $b_j:=a_j$ for $y_j=\Zero$ is a satisfying one.
According to Example~\ref{x:IndepOddD},
this makes all clauses $u\vee v$, $(u\neq v)$, evaluate to \true
for which at least one of the literals $u,v$ are assigned 
to the \emph{complement} of some $a_j$.
On the other hand, since $u[\bar y]\vee v[\bar y]=\One$,
by construction also
at least one of $u[\bar b]$ and $v[\bar b]$ 
is of the form $\neg a_j$.
\item[c\,iii$\Rightarrow$iv)] follows 
from Example~\ref{x:IndepOddD}. And
c\,iv$\Rightarrow$i)~ is a tautology.
\item[d)] 
Recursive application of de Morgan's Laws as in Example~\ref{x:QL}a).
\item[e)]
is straightforward by term induction
since $\vee$ and $\wedge$
 are $\sup$ and $\inf$, respectively.
\qed\end{longenum}\end{proof}
Our next tool is designed to deal with Item~c) of
Theorem~\ref{t:Syntax}.   We consider
 formulae $\phi(\bar x)$ of the form
\[  \exists \bar z.\;\Band\nolimits_{i=1}^k t_i(\bar x, \bar z) =\One  
\; \Band\; \Band\nolimits_{j=1}^\ell  s_j(\bar x, \bar z)=\Zero\]
with lattice terms  $t_i(\bar x, \bar z) \in  \bigwedge \bigvee  $ 
and $s_j(\bar x, \bar z ) \in \bigvee \bigwedge $. 
Let $\Lambda$ denote the set of all such formulae
(that is, certain \emph{pp-formulae} in the
sense of Model Theory).
Recall that a  complemented modular lattice
is a modular lattice $L$ with $\Zero$ and $\One$ such that
for any $a \in L$ there is $b \in L$ such that
$a \wedge b=\Zero$ and $a \vee b=\One$.
Examples are the lattices of all subspaces of 
vector spaces.

\begin{myproposition} \lab{p:collin}
\begin{enumerate}
\item[a)] The Boolean conjunction of two 
formulae  
 in $\Lambda$ is logically equivalent to a formula in
$\Lambda$.
\item[b)] For any  $\phi(\bar x)$ in $\Lambda$ 
 there exist lattice terms $t(\bar x,\bar z)\in 
\bigwedge \bigvee \bigwedge$ and $s(\bar x,\bar z)   
\in \bigvee \bigwedge \bigvee$ 
such that, within any lattice, $\phi(\bar x)$ is 
equivalent to $\exists \bar z.\; t(\bar x,\bar z)=\One  \;\Band\;
s(\bar x,\bar z)=\Zero$.
\item[c)]  Given lattice terms $t(\bar x, u_1, \ldots u_m)$ 
and $t_i(\bar x,u_i)$ ($i=1, \ldots ,m$) such that,
within the lattice $L$, 
the formulae $v=t(\bar x, u_1, \ldots ,u_m)$ 
and $u_i =t_i(\bar x)$ are equivalent to  
 $\phi(v,\bar x, u_1, \ldots ,u_m)$ 
and $\phi_i(u_i,\bar x)$,  respectively, then
\[ \exists u_1 \ldots \exists u_m.\; \phi(v,\bar x,u_1, \ldots u_m)
\;\Band\; \phi_1(u_1,\bar x) \;\Band \ldots \Band \; \phi_m(u_m,\bar
x) \] 
is  equivalent within $L$ to
$v= t\big(\bar x, t_1(\bar x), \ldots ,t_m(\bar x)\big)$.
\item[d)]  To any of the formulae
$x \leq y$, $x=y$, $x \leq y\vee z$, $x \vee y=x \vee z$,
there exists a formula in $\Lambda$ equivalent to 
it within complemented modular lattices.
\item[e)] To every $d\in\IN$ there is 
$\chi_d(\bar y)$  in $\Lambda$  
such that in any complemented modular lattice $L$
one has $L \models \chi_d(\bar a)$ if and only
if $\bar a$ is a $d$-frame of $L$.
\item[f)] To every $d \geq 3$ 
there is $\rho_d(x, \bar y)$  in $\Lambda$
such that for any 
a $d$-frame $\bar A$ in  $\Gr(\calH)$,
$\dim(\cal H)=d$, one has 
$P \in \IntRing_{\bar A}$ if and only $\Gr(\calH)\models
\rho_d(P,\bar A)$.  
\item[g)] 
For any   variable $u$,
  and any given $d \geq 3$ there  are
$\phi_{d,\ominus}(u,x_1,x_2,\bar y)$  and $\phi_{d,\otimes}(u,x_1,x_2,\bar y)$
in $\Lambda$  such that for
any  $d$-frame $\bar A$ in  $\Gr(\calH)$,
$\dim(\cal H)=d\geq 3$,  and $R,P,Q$  in $\IntRing_{\bar A}$
one has $R=P\ominus_{\bar A}Q$   if and only
$\Gr(\calH) \models \phi_{d,\ominus}(R,\bar P,\bar A)$
and  $R=P\otimes_{\bar A}Q$   if and only
$\Gr(\calH) \models \phi_{d,\otimes}(R,\bar P,\bar A)$.
\item[h)] In any of the  above, the claimed formulae
can be obtained in polynomial time from the  data. 
\end{enumerate} 
\end{myproposition} 
\begin{proof} 
Item~a) follows from the basic  rules of logic dealing
with $\exists$ and $\Band$,
In b) let $t=\bigwedge_i t_i$ and $s=\bigvee_j s_j$. 
In c) generalize the reasoning of
Observation~\ref{o:RelInterp}). 
The first two cases in Item~d) are dealt with in Example~\ref{x:QL}f).
Now,  by modularity and existence
of complements,  $x \leq y\vee z$ is equivalent to
$\exists u.\; y \vee z  \vee u =\One \;\Band\;   (x \vee y\vee z) \wedge u  
=\Zero$  and, in any  lattice,  $x\vee y=x\vee z$ to the conjunction of
$y \leq x \vee z$ and $z \leq x \vee y$. 
e) and f) follow immediately --- only  the last condition
in the definition of a frame needs a closer look: 
it is equivalent to 
$a_{ik} \leq  a_{ij} \vee a_{jk}$; indeed, given this
$a_{ik}= (a_i \vee a_k) \wedge  (a_{ij} \vee a_{jk})$ follows by
modularity.

In g),
considering  the subterms of $P \ominus Q$ and $P \otimes Q$ 
and  the calculations in Subsection~\ref{ss:IntField}
(resp. their geometric motivations in Subsection~\ref{ss:Affine})  
we derive the  following,    which proves g) in view of d), c), and a),
\[\begin{array}{lclcl }  
u_1=Q_{13}\quad&\Leftrightarrow & u_1 \leq  Q +A_{13} \;&\Band& 
u_1 \leq A_1 +A_3\\    
u_2=P_{32} &\Leftrightarrow & u_2 \leq P +A_{13} \;&\Band& 
u_2 \leq A_2+ A_3 \\   
u_3=P \otimes Q&\Leftrightarrow &
  u_3 \leq Q_{13}+P_{32}\;&\Band& u_3 \leq A_1+ A_2\\
u_4=S  &\Leftrightarrow &
u_4 \leq Q_{13} +A_2\;&\Band& u_4\leq P+ A_{23}\\
u_5= P\ominus Q&\Leftrightarrow:&u_5 \leq S+A_3 \;&\Band&
u_5\leq A_1+ A_2.\end{array}  \]  
h) is clear from the above proofs.
\end{proof} 

\begin{proof}[of Theorem~\ref{t:Syntax}]
\begin{longenum}
\item[a)]
Given a $\bigwedge\bigvee$-term $t(x_1,\ldots,x_n)$,
first eliminate all clauses with only one literal
by substituting it with $\One$: 
This simplification can obviously be performed in
polynomial time and maintains $t$'s 
$\bigwedge\bigvee$-form as well as strong satisfiability.
If it fails (like for instance in $\neg x\wedge \neg y\wedge(x\vee y)$),
reject $t$. Otherwise, in the even-dimensional case, accept.
In odd dimensions, collect all 
clauses with precisely two (remaining) literals
and report whether this instance of $\textsf{2SAT}$ 
is satisfiable: as mentioned above, in polynomial time.
It remains to assert the correctness of this algorithm.
Regarding the case of even dimensions,
this holds due to Lemma~\ref{l:Conjunctive}a).
Over $\Gr(\IF^{2d+1})$, any satisfying assignment
of $t$ must in particular make all its two-literal 
clauses evaluate to \true; which requires their conjunction 
to be a positive instance for \textsf{2SAT} according to
Lemma~\ref{l:Conjunctive}c\,i$\Rightarrow$ii).
Conversely, if the two-literal clauses are satisfiable
over $\{0,1\}$ then both, they \emph{and} the clauses
with at least three literals, admit a joint satisfying assignment
from $\{a_1,\neg a_1,\ldots,a_n,\neg a_n\}\subseteq\Gr(\IF^{2d+1})$
according to Lemma~\ref{l:Conjunctive}b) and 
c\,ii$\Rightarrow$iii). 
 \item[b)] 
By Lemma~\ref{l:Conjunctive}d) we may w.l.o.g. presume the input
$t(\bar x)$ to be in negation normal form, i.e. equal to
$\hat{t}(\bar x,\neg \bar x)$ with $\hat{t}(\bar x,\bar y)$ over $\vee$ and $\wedge$ only. 
The relational interpretation (Observation~\ref{o:RelInterp})
produces from that in polynomial time a system of basic lattice equations,
i.e. each of one of the forms 
$y_j=y_k\wedge y_\ell$ and $y_j=y_k\vee y_\ell$
in (possibly negated) variables $y_j$
jointly satisfiable iff $t$ admits a strongly satisfying assignment.
Now according to Example~\ref{x:QL}f), these equations are in turn equivalent to
\[ \One\;=\; (y_j\wedge y_k\wedge y_\ell)\vee\big(\neg y_j\wedge(\neg y_k\vee\neg y_\ell)\big)
\;\;\text{ and }\;\;
\One\;=\; 
\big(y_j\wedge(y_k\vee y_\ell)\big)\vee(\neg y_j\wedge\neg y_k\wedge\neg y_\ell) . \]
Their right hand sides are of the form $\bigvee\bigwedge\bigvee$;
hence combining these equations 
according to Fact~\ref{f:clear}a)
yields a single $\bigwedge\bigvee\bigwedge\bigvee$-term $s(\bar y)$
as desired.
\item[c)]
To see the problem in $\calBP(\calNP_\IF)$ 
(and not just in $\calBP(\calNP_{\Re\IF})$),
recall that the evaluation 
in Proposition~\ref{p:UpperComplexity}b) 
needs access to real and imaginary parts only
for terms with negation when invoking
Proposition~\ref{p:UpperComplexity}a\,iii);
whereas for $s$ and $t$ as in this case,
Proposition~\ref{p:UpperComplexity}a\,i)
and a\,ii) suffices.
\\
Regarding completeness, observe that steps
i) to v) in the proof of Proposition~\ref{p:IntField}
produce a finite list of lattice equations: terms $\otimes_{\bar A}$ and $\ominus_{\bar A}$ 
do not invoke negation in ii) and iii), nor does iv)
or v) according to Definition~\ref{d:frame}.
Replacing in Step~vi) Fact~\ref{f:clear}a) with
Fact~\ref{f:clear}b) thus yields two lattice
terms $s$ and $t$ with the claimed properties.
\item[d)]
Regarding completeness as in the proof of c), observe that,
due to Proposition~\ref{p:collin},  Steps~i) to v) in the proof of Proposition~\ref{p:IntField}
produce a finite  list of formulae in $\Lambda$:
 in v) use e) and  in iv) use f). 
In ii) observe that any  polynomial $p$ with integer coefficients
 is equivalent within
(commutative) rings to a term $q(\bar x)$   in constants $0,1$ 
and binary operation symbols $-, \cdot$.  Replace these
according to  Proposition~\ref{p:collin}g)  by $y_{11}$, $y_{12}$, $\phi_{d,\ominus}$,
and $\phi_{d,\otimes}$, respectively,  
  and apply, iteratively, c) and a) of
Proposition~\ref{p:collin}  to obtain
$\phi_{d,p}(u,\bar x,\bar y)$  in $\Lambda$
such that for any $d$-frame $\bar A$ in $\Gr(\calH)$,
$\dim(\calH)=d\geq 3$ and any $R,\bar P$ in $\IntRing_{\bar A}$ one has
$R=p_{\IntRing_{\bar A}}(\bar P)$ if and only if
$\Gr(\calH) \models \phi_{d,p}(R,\Bar P,\bar A)$.   This deals with
Step~iii).    
 Replacing in Step~vi) Fact~\ref{f:clear}b) with
Proposition~\ref{p:collin}b) thus yields two lattice
terms $s$ and $t$ with the claimed properties --- and all this 
in time polynomial in the input length of $p$.
\item[e)] Let $d\geq 3$.
Equation $s=\Zero$ with $s \in \bigvee \bigwedge \bigvee$
is equivalent to $\neg s=\One$ and 
$\neg s $ is equivalent to a term in $ \bigwedge \bigvee \bigwedge$
by de Morgan's Laws.  
In the case $\IF \not\subseteq \IR$ 
we also have to take in account Steps~(v') and (iv') of the proof of 
Theorem~\ref{t:IntStarField}.  
Here, the additional conditions $A_i \leq \neg A_j$ 
and  $X_k \leq \neg  Z_k$  are dealt with by Example~\ref{x:QL}f), again.
In the case $d \leq 2$, for hardness  we refer to the proof of 
Proposition~\ref{p:NPhard}  and the fact that the
commutator is in $\bigvee \bigwedge$;   membership  in $\calNP$ follows
from Proposition~\ref{p:2DinNP}.
\qed\end{longenum}\end{proof}

\subsection{Varying the Ground Field} \lab{ss:Fields}
In dimension two, weak/strong satisfiability
did not depend on the underlying field under 
consideration: recall Lemma~\ref{l:2D}b).
This becomes different starting with dimension three. 

\begin{fact} \lab{f:Fields} 
\begin{enumerate}
\item[ a)] If $\IF$ is a $\ast$-subfield of
$\IE$ then $\Gr(\IF^d)$ embeds into $\Gr(\IE^d)$.
\item[ b)] $\Gr\big(\IF(i)^d\big)$ embeds into $\Gr(\IF^{2d})$.   
\item[ c)] For an irreducible monic polynomial $p(X)\in\IF[X]$ of
  degree $k$ and  $\IF$-vector space 
$V$,  there is an  endomorphism $\phi$ of $V$ such that
$p(\phi)=0$ if and only if $k \leq \dim(V)$. 
\COMMENTED{
the associated \emph{companion matrix} $A_p\in\IF^{d\times d}$
satisfies $p(A_p)=0$.
\item[ d)] Let $V$ denote an $\IF$-vector space,
$\phi:V\to V$ a linear map, and suppose the
non-zero polynomial $p(X) \in \IF[X]$ is irreducible
such that $p(\phi)=0$. Then $\dim(V)\geq\deg(p)$. }
\end{enumerate} 
\end{fact} 
\begin{proof}  
In a)  map $U$ to $\IE U$.
In b) 
consider  $\IF(i)^d$ as an $\IF$-vector space
with scalar product the real part of that in $\IF(i)^d$.
For c)  See \mycite{Ch.III.4}{Jacobson}. 
\COMMENTED{ Choose a zero $\alpha$  of $p(X)=\sum_{k=0}^d a_kX^k$ in $\IC$.
Then $\IF(\alpha)$  has basis $1, \ldots ,\alpha^{d-1}$ 
and $A_p$ is the matrix of
the endomorphism defined by $\phi(r)= \alpha r$.
Now $p(\phi)(r) = \sum_{k=0}^d a_k \phi^k(r) 
=r \sum_{k=0}^d a_k \alpha^k =0$ and it follows $p(A_p)=0$
cmp. \cite[Ch. VII, Thm.4.3]{Hunger}.      
\item[d)]
 Follows
from  the \textsf{Caley-Hamilton Theorem} 
(cmp. \cite[8.33]{Axler}, \cite[Proposition 3.10]{Farenick}): 
Consider a nontrivial cyclic submodule $U$ of the
$\IF[X]$-module defined by $\phi$  on $V$.
Then the minimal polynomial $q$ of $\phi|U$ divides
$p$, and vice versa by irreducibility.
Thus $\deg(q)=\deg(p)$. 
On the other hand, $\dim(U)=\deg(q)$ since
$U$ is cyclic.
\qed\end{longenum}}\end{proof} 
Recall that two algebraic structures are
\textsf{elementarily equivalent} if they satisfy the same first-order
sentences. In particular, all real closed fields $\IF$ 
(such as $\IR$ and $\IA \cap \IR$)   are pairwise  elementarily
equivalent (cf. \mycite{\S 5.5}{Shoenfield});
and so are all the associated  $\ast$-fields  $\IF(i)$. 

\begin{myproposition} \lab{p:Fields}
\begin{enumerate}
\item[a)] For elementarily equivalent  $\IF_1,\IF_2$ 
and $\IF_i$-unitary spaces $\calH_i$  of the same dimension $d$
admitting equinormal orthogonal bases,
every orthologic term $t$ is weakly (resp. strongly) satisfiable
in $\Gr(\calH_1)$ if and only if it is so in $\Gr(\calH_2)$.
\item[b)] A term is weakly (resp. strongly) satisfiable 
either in both or in none of $\Gr(\IC^d)$ and $\Gr(\IA^d)$;
similarly for $\Gr(\IR^d)$ and $\Gr\big((\IR\cap\IA)^d\big)$.
\item[c)]  If $\IF$ is a $\ast$-subfield of $\IE$ then
any term weakly resp. strongly satisfiable
 in $\Gr(\IF^d)$ is so in $\Gr(\IE^d)$.
\item[d)] Any term weakly resp. strongly satisfiable
 in $\Gr\big(\IF(i)^d\big)$ is so in $\Gr(\IF^{2d})$. 
\item[e)] Let polynomial $p(X)\in\IZ[X]$ have a root
in $\IE$ but not in $\IF$ and $d\geq3$.
Then the term $t_{p,d}$ according to Proposition~\ref{p:IntField}c)
is strongly satisfiable in $\Gr(\IE^d)$ 
 but not in $\Gr(\IF^d)$.
\item[f)] For any $d\geq 3$  and any $k\geq 2$ 
there exists a term  $t$ such that 
$t$ is strongly satisfiable in $\Gr(\IR^d)$
and in $\Gr(\IQ^{kd})$ 
but not in any $\Gr(\IQ^n)$  with $n<kd$.
\end{enumerate}
\end{myproposition} 
So Item~e) says that for instance $t_{X^2-2,d}$ 
is strongly satisfiable over $\Gr(\IR^d)$
but not over $\Gr(\IQ^d)$;
and $t_{X^2+1,d}$ 
is strongly satisfiable over $\Gr(\IC^d)$
but not over $\Gr(\IR^d)$.
Proposition~\ref{p:Fields}e) and f)  has
analoga  for
 weak satisfiability according to   Theorem \ref{t:WeakStrong}b).

\begin{proof}[of Proposition~\ref{p:Fields}]
\begin{longenum}
\item[a)] As observed, e.g., in \mycite{Theorem~6}{Hagge1},
weak/strong satisfiability over $\Gr(\IF^d)$ amounts to a
first-order sentence
to be valid in $\IF$. In case $\calH$ admits an equinormal
orthogonal basis, invoke Lemma~\ref{l:Equinormal}a).
\item[b)] follows from a); \quad c) follows from Fact~\ref{f:Fields}a);
\quad d) follows from Fact~\ref{f:Fields}b).
\item[e)] is immediate by Proposition~\ref{p:IntField}c).
\item[f)] Consider a real algebraic integer $\alpha$ of
degree $k$ over $\IQ$ like, for instance, $\alpha:=2^{1/k}$;
cmp. \mycite{Ch. VI Thm.9.1}{Lang}.
Its associated minimal polynomial $p(X)\in\IZ[X]$ is monic
and irreducible --- over $\IQ[X]$ by Gauss' Lemma. 
Fact~\ref{f:MatinGrH}c) below yields a term $t_{p,d}$ 
strongly satisfiable over $\Gr(\IF^n)$ iff 
  there is
an endomorphism $\phi$ of an $m$-dimensional subspace 
$W$ of $\IF^n$ where $n=dm$.  In view of Fact~\ref{f:Fields}c) 
this happens for $\IF=\IQ$ if and only if $k \leq m$.
On the other hand, if $\IA \cap \IR \subseteq \IF$ 
then $p(X)$ factors in $\IF[X]$  into linear factors 
  and  one 
may choose $m=1$. 
\COMMENTED{
a $W\in\Gr_m(\IF^{dm})$ and homomorphism $\phi:W\to W$ with $p(\phi)=0$. 
Now for $\IF=\IA$, such $(W,\phi)$ exists for $m=1$,
namely $W=\IA\times\{0\}^{d-1}$ and $\phi=\alpha\cdot\id_W$. 
For $\IF=\IQ$, such a $(W,\phi)$ exists as well: let $m=k$
and $W=\IQ^k\times\{0\}^{d(k-1)}$ and $\phi$ 
the linear map induced by the companion matrix of $p$
with respect to the standard orthonormal basis
(Fact~\ref{f:Fields}c).
Whereas, according to Fact~\ref{f:Fields}d), no 
$W\in\Gr_m(\IF^{dm})$ for $m=\dim(W)<k$ admits a $\phi\in\Hom(W)$
with $p(\phi)=0$.} 
\qed\end{longenum}\end{proof} 
\cx{
\begin{digression} \lab{d:uclass} 
For $\ast$-subfields  $\IF \subseteq \IE$ of $\IC$,
if $\IE$   is a purely transcendental extension
of $\IF$ then one has an embedding of
$\Gr(\IE^d)$  into some ultrapower of $\Gr(\IF^d)$. 
If $\IE$ is an algebraic extension of $\IF$
then  $\Gr(\IE^d)$  is the direct union of  the
$\Gr(\IE_i^d)$ with $\IE_i$  ranging over the
finite algebraic extensions of $\IF$. 
Thus, the interesting case is
$\IE=\IF(\alpha)$ where $\alpha$ is a root of an
irreducible polynomial $p(X)$  in $\IF[X]$. 

Recall that 
$\Gr(\IE^d)$ embeds into  $\Gr(\IF^n)$ if there is a 
$\ast$-ring embedding  of $\IE^{d \times d}$ into $\IF^{n \times n}$
(whence, in particular, $n=kd$ for some $k$)
and that for  $d\geq 3$ this is also necessary (this  can be derived
from the proofs of the coordinatization results of \cite{Neumann2}).
For $\IE=\IF(i)$   there is a canonical such embedding 
with $k=2$:  consider
$\IF(i)\cong\big\{ \binom{a \; -b}{b \;\; a} \mid a,b \in \IF\big\}$
and convert the entries of a matrix in $\IF(i)^{d \times d} $
into blocks of a matrix over $\IF$, accordingly (cmp. Fact~\ref{f:Fields}b).
In view of Digression~\ref{d:ImagField} 
this leaves to consider the case $\IE \subseteq \IR$.
Here,   existence of a 
$\ast$-ring  embedding amounts  to existence of 
a symmetric $A \in \IF^{k \times k}$ such that $p(A)=0$.
 For $p(X)=X^2-c$ and $k=2$ such $A$  are of those of the form
$A= \binom{a \;\; b}{b \;-a}$ with $c=a^2+b^2$.
Such $A \in \IQ^{2 \times 2}$ exist for  $c=2$ and primes $c \equiv 1\, {\rm
  mod}\, 4$, but not for integers
$c \equiv 3\, {\rm
  mod}\, 4$. 

Now, for a prime $p \equiv 1 \;{\rm mod}\,4$ we have 
integers $a,b$ such that $p =a^2+b^2+2$. Consider
$\IF=\IQ(\sqrt{2},i)$  and the  selfadjoint matrix
  $A=\binom{\!\!\!\!\sqrt{2}  \;\: c}{ c^* \; -\sqrt{2}}$
where $c=a+bi$.
Then $A^2=pI$ and we conclude that $\Gr(\IF(\sqrt{p})^d)$ 
 embeds into $\Gr(\IF^{2d})$ whence in $\Gr(\IQ^8)$.

To us, all other cases are open.  As is the question whether there
is an ortholattice identity valid in all $\Gr(\IQ^n)$ but
not in $\Gr(\IR^3)$. Though, such identity cannot be a mere
lattice identity.
\end{digression} }

\subsection{Dimensions of Satisfiability} \lab{ss:Dimensions}
\begin{observation} \lab{o:Dimensions}
\begin{enumerate}
\item[a)]
A term weakly satisfiable in $\Gr(\calH)$
is also weakly satisfiable in $\Gr(\calH')$
for every subspace $\calH\subseteq\calH'$.
\item[b)]
A term strongly satisfiable in both $\Gr(\calH)$
and in $\Gr(\calH')$
is also strongly satisfiable in $\Gr(\calH\obot\calH')$.
\end{enumerate}
\end{observation}
Put differently, for each term $t$,
the set $\dim_{\IF}(t)\subseteq\IN$ is an ideal in $(\IN,+)$
and $\DIM_{\IF}(t)\subseteq\IN$ a sub-semigroup, where
\begin{equation} \label{e:Dimensions}
\begin{gathered}
\dim_{\IF}(t) \;:=\; \big\{d\in\IN\mid t\text{ is weakly satisfiable over }\Gr(\IF^d)\big\} \\
\DIM_{\IF}(t) \;:=\; \big\{d\in\IN\mid t\text{ is strongly satisfiable over }\Gr(\IF^d)\big\}
\end{gathered}\end{equation}
Indeed, a) follows from Fact~\ref{f:Foulis}b)
and the fact that weak satisfiability of $t$ is the
complement of validity of the equation $t=\Zero$,
the latter being preserved in substructures
and homomorphic images. 
b) follows from Observation~\ref{o:Product}a).

\begin{myexample} \lab{x:Dimensions}
\begin{enumerate}
\item[a)] For any $d$, the term $h_d$ from Lemma~\ref{l:Huhn}a)
is (independent of $\IF$ and)
strongly satisfiable in $\Gr(\calH)$ if and only if 
$d$ divides $\dim(\calH)$.
\item[b)] For any $k$, there exists a term $t_k$ of length
$\calO(k)$ strongly satisfiable over $\Gr(\calH)$ for
$\dim(\calH)=2^k$ but not for $\dim(\calH)<2^k$.
\item[c)] For any $d$, 
the term $g_d$ from Lemma~\ref{l:Huhn}c) 
is (independent of $\IF$ and)
weakly satisfiable in $\Gr(\calH)$ if and only if $d\leq\dim(\calH)$.
\end{enumerate} 
\end{myexample} 
\begin{proof} 
\begin{longenum}
\item[a)] 
Suppose $\dim(\calH)=d\cdot k$. Then 
$\calH=V_1\obot\cdots\obot V_k$ for
appropriate $V_j$ pairwise orthogonal.
By Example~\ref{x:huhn}a),
there exists an orthogonal $d$-diamond
$\bar A_j$ in $\Gr(V_j)$, $1\leq j\leq d$;
and thus one in $\Gr(\calH)$ by Observation~\ref{o:Product}a):
strongly satisfying $h_d$ by Lemma~\ref{l:Huhn}a).
\\
Conversely, any strongly satisfying assignment $\bar A$ of $h_d$
over $\calH$ 
constitutes by Lemma~\ref{l:Huhn}a) an orthogonal $d$-diamond;
hence $d$ divides $\dim(\calH)$ according to Fact~\ref{f:Huhn}c).
\item[b)] Consider the $(2n+1)$-variate term 
\begin{multline*}
x_{n+1}\;\wedge\;\; \bigwedge\nolimits_{i=1}^n \bigg(
(x_i\vee\neg y_i)\;\wedge\;(y_i\vee\neg x_i)\;\wedge\;
(\neg x_i\vee\neg y_i) \\[-2ex] \;\wedge\;
\Big(\big(x_{i+1}\wedge(x_i\vee y_i)\big)\vee\big(\neg x_{i+1}\wedge
\neg(x_i\vee y_i)\big)\Big)\bigg)
\end{multline*}
Note that $x_i\vee\neg y_i=\One$ in $\Gr(\IF^d)$
implies $\dim(x_i)+\dim(\neg y_i)\geq d$;
hence the first two terms in the big conjunction
require $\dim(x_i)=\dim(y_i)$.
The fourth term amounts to condition
$x_{i+1}=x_i\vee y_i$ according to Example~\ref{x:QL}f);
hence $\dim(x_{i+1})=\dim(x_i)+\dim(y_i)$ because of 
$x_i\wedge y_i=\Zero$ (third term). 
Concluding, any satisfying assignment $(a_1,b_1,\ldots,a_n,b_n,a_{n+1})$
has $\dim(a_{i+1})=2\cdot\dim(a_i)$.
Therefore $2^n\cdot\dim(a_1)=\dim(a_{n+1})=d$ 
by the very first term. Conversely, 
the following is easily verified
to constitute a satisfying assignment:
\begin{multline*}
a_1:=\IF\times\{0\}, \;\; b_1:=\{(x,x):x\in\IF\}, \quad\; 
a_2:=\IF^2\times\{0\}^2, \;\; b_2:=\{(\vec x,\vec x):\vec x\in\IF^2\}, \\
\quad \ldots, \quad
a_{i+1}:=\IF^{2^i}\times\{0\}^{2^i}, \;\; b_{i+1}:=\{(\vec x,\vec x):\vec x\in\IF^{2^i}\}
\quad \ldots, \quad a_{n+1}=\IF^{2^n} 
\end{multline*}
(all understood embedded into $\IF^{2^n}$ by appending zeros).
\item[c)]
As pointed out in Corollary~\ref{c:Hagge},
$g_d$ is not weakly satisfiable over $\Gr(V)$ for $\dim(V)<d$
but weakly satisfiable for $\dim(V)=d$.
In case $\dim(\calH)>d$, there exists $V\in\Gr_d(\calH)$;
hence by Observation~\ref{o:Dimensions}a),
$g_d$ is weakly satisfiable also over $\Gr(\calH)$.
\qed\end{longenum}\end{proof} 
 \cx{
\begin{digression}
An alternative term 
(and historically the perhaps first,
although of length cubic in $d$)
weakly satisfiable in $\Gr(\calH)$ if and only if $d\leq\dim(\calH)$
is also due to \person{Huhn} \cite{Huhn}: He has shown 
that the ``\emph{$(d-1)$-distributive law}''
\begin{equation} \label{e:d-distributive}
 x \vee \bigwedge\nolimits_{j=1}^d y_j = \bigwedge\nolimits_{j=1}^d \bigl(x \vee
\bigwedge\nolimits_{j\neq i=1}^d y_i \bigr) 
\end{equation}
fails in a modular lattice $L$ if and
only if $L$ contains $v<u$ and a
$d$-diamond in $[v,u]$ --- where the condition
$a_i \leq \neg a_j$ for $i \neq j$ is replaced
by the independence requirement
$a_i \wedge \bigvee_{j=1,\,j\neq i}^d =\Zero$ for $i\geq 1$.
As in Fact~\ref{f:Huhn} this implies $\dim[v,u]\geq d$.
Now consider $r \wedge \neg \ell$
and observe that $\ell\leq r$
holds in any lattice, where $\ell$ and $r$
denote left and right hand side of Equation~(\ref{e:d-distributive}), 
respectively.
\end{digression}}
A partial converse to Example~\ref{x:Dimensions}c)
appears in Item~b) of the next

\begin{mylemma} \lab{l:Dimensions}
Recall from Lemma~\ref{l:Conjunctive}d)+e) that 
a lattice term is one over $\vee$ and $\wedge$ only,
i.e. contains no negations.
\begin{enumerate}
\item[ a)]  For any lattice term $t$  and  $\vec v \in t_{\Gr(\calH)}(U_1, \ldots ,U_n)$ 
there are $V_i \subseteq U_i$  in $\Gr(\calH)$ 
such that $\sum_{i=1}^n \dim V_i \leq|t|$   and
$\vec v \in   t_{\Gr(\calH)}(V_1, \ldots ,V_n)$. 
\item[b)]
A term $t$ is weakly satisfiable in $\Gr(\calH)$ 
if and only if there exists $W\in \Gr(\calH)$ with
$\dim W\leq|t|$ such that $t$ is weakly satisfiable in $\Gr(W)$.
\end{enumerate} 
\end{mylemma}
\begin{proof} 
\begin{longenum}
\item[a)] by term induction:
Let $t=t_1 \vee t_2$  and  $m_j:=|t_j|$ whence
$|t|\geq m_1+m_2$.    Then
$\vec v \in t_{1 \Gr(\calH)}(U_1, \ldots ,U_n) 
+t_{2 \Gr(\calH)}(U_1, \ldots ,U_n)$, i.e.
there are $\vec v_j \in   t_{j,\Gr(\calH)}(U_1, \ldots ,U_n)$   
such that $\vec v =\vec v_1 +\vec v_2$
and by inductive hypothesis there are $V_{ji} \subseteq U_i$ 
for $j=1,2$ such that $\sum_{i=1}^n \dim V_{ji}\leq m_j$ 
and $\vec v_j  \in    t_{j,\Gr(\calH)}(V_1, \ldots ,V_n)$.
So put $V_i:=V_{1i}+V_{2i}$. 
\\
For $t= t_1 \wedge t_2$
the same arguing applies with $\vec v =\vec v_1=\vec v_2$.
\item[b)]
In view of Lemma~\ref{l:Conjunctive}d) we may assume $t(\bar x)=\hat t(\bar x,\neg\bar x)$
with a lattice term $\hat t(\bar x,\bar y)$.
Now, assume $\vec v \in t_{\Gr(\calH)}(U_1, \ldots ,U_n)$.
Then $\vec 0\neq  \vec v \in \tilde{t}_{\Gr(\calH)}(U_1, \ldots ,U_n,U_1^\perp,
\ldots ,U_n^\perp)$.
According to a) 
there are  $V_i \subseteq U_i$ and $V'_i \subseteq U_i^\perp$ 
such that $\sum_{i=1}^n \dim V_i +\dim V_i' \leq m =|t|$ 
and $\vec v \in  
\tilde{t}_{\Gr(\calH)}(V_1, \ldots ,V_n,V_1',
\ldots ,V_n')$.  Put $W:=\sum_{i=1}^n V_i+V_i'$
and $W_i :=U_i \cap W$. 
Then $\dim W \leq m$, $V_i \subseteq W_i \subseteq U_i$, and
$V_i' \subseteq  W\cap W_i^\bot$ since $V_i' \subseteq 
U_i^\bot \subseteq W_i^\bot$. By Lemma~\ref{l:Conjunctive}e) 
it follows $\vec v \in 
 \tilde{t}_{\Gr(\calH)}(W_1, \ldots ,W_n, W \cap W_1^\perp,
\ldots , W \cap  W_n^\perp)
=t_{\Gr(W)}(W_1,\ldots, W_n)$.
\qed\end{longenum}\end{proof}

\section{Satisfiability in Indefinite (yet Finite) Dimension} \lab{s:Indef}
In an aim to approach the infinite-dimensional case,
we now consider satisfiability questions quantifying existentially
not just over assignments but also over the (finite) dimension
the assignment lives in. 

\begin{mydefinition}  \lab{d:Indef}
\begin{enumerate}
\item[a)]
Call a term weakly respectively
strongly satisfiable in $\Gr(\IF^*)$
if and only if it is so in $\Gr(\IF^d)$ for some $d\in\IN$.
\item[b)]
The corresponding decision problems are abbreviated as
$\sat_{\Gr(\IF^*)}:=\sat_{\{\Gr(\IF^d):d\in\IN\}}=\bigcup_d\sat_{\Gr(\IF^d)}$
and
$\SAT_{\Gr(\IF^*)}:=\SAT_{\{\Gr(\IF^d):d\in\IN\}}=\bigcup_d\SAT_{\Gr(\IF^d)}$,
respectively.
\end{enumerate}
\end{mydefinition}
\begin{theorem} \lab{t:WeakIndef}
\begin{enumerate}
\item[a)]
For Pythagorean $\IF$ it holds 
$\sat_{\Gr(\IF^*)}\in\calBP(\calNP_{\Re\IF})$.
In particular, $\sat_{\Gr(\IR^*)},\sat_{\Gr(\IC^*)}\in\calBP(\calNP_{\IR})$.
\item[b)]
A term is weakly (resp. strongly) satisfiable 
either in both or in none of $\Gr(\IF^*)$ and $\Gr(\Re\IF^*)$,
i.e. it holds $\SAT_{\Gr(\IF*)}=\SAT_{\Gr(\Re\IF^*)}$
and $\sat_{\Gr(\IF*)}=\sat_{\Gr(\Re\IF^*)}$.
\end{enumerate}
\end{theorem}
That is, weak satisfiability over $\Gr(\IF^*)$
is decidable by a nondeterministic polynomial-time 
\BSS-machine over $\Re\IF$---and thus no more hard
than in the fixed-dimensional case.
A qualitative precursor to this result 
is the decidability of the equational
theory of the class $\{\Gr(\IC^d)\mid d<\omega\}$ in \cite{hn}. 

\begin{proof}
\begin{longenum}
\item[a)]
According to Lemma~\ref{l:Dimensions}d), $t$ is weakly
satisfiable over $\Gr(\IF^*)$ iff it is so over 
$\Gr(W)$ for some $d$-dimensional $\IF$-unitary space $W$
for $d:=|t|$; iff it is over $\Gr(\IF^d)$, according to
Fact~\ref{f:GramSchmidt}. Now recall 
(Proposition~\ref{p:UpperComplexity}c)
that satisfiability over $\Gr(\IF^d)$
can be decided by a 
nondeterministic \BSS-machine over $\Re\IF$
in time polynomial in $|t|$ \emph{and $d$}.
\item[b)]
follows from Proposition~\ref{p:Fields}c+d)
in view of $\IF=\big(\Re\IF\big)(i)$ for $\ast$-fields $\IF\not\subseteq\IR$.
\qed\end{longenum}\end{proof}
\begin{myquestion} \lab{q:WeakIndef}
Is $\sat_{\Gr(\IR^*)}$ hard for $\calNP$
or even for $\calBP(\calNP^0_\IR)$?
\end{myquestion}
The rest of this section explores a strong counterpart to
Theorem~\ref{t:WeakIndef}a) and Question~\ref{q:WeakIndef}, 
namely the complexity (and computability) of $\SAT_{\Gr(\IF^*)}$.
To this end, the next subsections extend 
Subsections~\ref{ss:vonNeumann}, \ref{ss:IntField}, and \ref{ss:OrthonFrame}
from interpreting into quantum logic not just the ring $\IF$ of scalars 
but the ring $\IF^{m\times m}$ of matrices, uniformly in $m$,
and similarly for $\ast$-rings.
We first discuss the feasibility problems for those.

\subsection{Feasibility in Matrix Rings} \lab{ss:intermat}
\COMMENTED{Matrices form a non-commutative ring
naturally equipped with transposition
(and possibly complex conjugation).
Recall that (cmp. \mycite{\S 2.13}{Rowen}, \mycite{\S 5.1}{Farenick})
a $\ast$-\textsf{ring} is a ring $R$ equipped
with an \textsf{involution} $r \mapsto  r^\ast$,
i.e. a map  such that  $(r^\ast)^\ast=r$,
$(r+s)^\ast=r^\ast+s^\ast$, and 
$(rs)^\ast =s^\ast r^\ast$.
Thus, the  ring $\IF^{d \times d}$ of $d$-by-$d$-matrices over 
the $\ast$-subfield $\IF$ of $\IC$  
becomes a $\ast$-ring with the 
involution $A \mapsto B= A^\ast$, the conjugate transpose
 of $A$, i.e. $b_{ij}= \bar a_{ji}$.}

\COMMENTED{
For (commutative) 
involutive rings we have already naturally been led
to consider $\ast$-polynomials in commuting variables.
}
Dealing with non-commutative rings,
the classical concept corresponding to the notion
of a term  amounts to that of polynomials
$p\in\IZ\langle X_1,\ldots,X_n\rangle$
with integer coefficients in non-commuting variables $\bar X$.
Similarly for non-commutative $\ast$-rings, a term
$t(X_1, \ldots,X_n)$   is
(equivalent to) a polynomial $p\in\IZ\langle X_1,\ldots,X_n,
X_1^\dagger,\ldots,X_n^\dagger\rangle$ --
again with the convention that $X^\dagg$ has to be interpreted
as $A^*=A^\dagg$ if $X$ is interpreted as $A$.  
And the ring  anti-automorphism  interchanging $X_i$ and 
$X_i^\dagg$ yields a polynomial equivalent to the term
 $t^\ast$ within $\ast$-rings.
These naturally suggest  to generalize the family of problems
$\FEAS_{\IZ,\Ring}$ from Example~\ref{x:BP} and the feasibility
of commutative $\ast$-polynomials from Theorem~\ref{t:IntStarField}.

\begin{mydefinition} \lab{d:Feas2}
\begin{enumerate}
\item[a)]
For a not necessarily commutative 
ring $\Ring\supseteq\IZ$ 
let $\FEAS_{\IZ,\Ring}$ denote the following decision problem:
\begin{quote}\it
Given polynomials
$p_1,\ldots,p_k\in\IZ\langle X_1,\ldots,X_n\rangle$,
do they admit a common root in $\Ring$, i.e. some 
assignment $\bar x\in\Ring^n$ such that
$p_1(\bar x)=\cdots=p_k(\bar x)=0$ ?
\end{quote} 
\item[b)]
$\QUAD{\Ring}$ is the similar problem but
restricted to quadratic polynomial systems
with coefficients from $\{0,\pm1\}$.
\item[c)]
For a $\ast$-ring $\Ring$, $\FEAD_{\IZ,\Ring}$ is the 
question of whether non-commutative $\ast$-polynomials
$p_j\in\IZ\langle X_1,X_1^\dagger,\ldots,X_n,X_n^\dagger\rangle$
admit a common root in $\Ring$.
\item[d)]
$\QUART{\Ring}$ is the similar problem but
restricted to one single quartic $\ast$-polynomial
$p(X_1,X_1^\dagger,\ldots,X_n,X_n^\dagger)$ 
with coefficients from $\{0,\pm1,\pm2,\ldots,\pm2n\}$.
\item[e)]
For $\calR$ a class of rings, 
$\FEAS_{\IZ,\calR}:=\bigcup_{R\in\calR} \FEAS_{\IZ,\calR}$ 
asks for a common root to $p_1,\ldots,p_k$ over
\emph{some} $R\in\calR$. Similarly for $\QUAD{\calR}$.
\item[f)]
For a class $\calR$ of $\ast$-rings, write
$\FEAD_{\IZ,\calR}:=\bigcup_{R\in\calR} \FEAD_{\IZ,\calR}$
and similarly for $\QUART{\calR}$.
\end{enumerate}
\end{mydefinition}
We point out that (the complements of) these computational 
problems arise in practice \cite{Benanti},
hence their complexity is worthwhile investigating.
Fact~\ref{f:RingInterp} carries over to the case
of noncommutative ($\ast$-) rings:

\begin{mylemma} \lab{l:RingInterp2}
\begin{enumerate}
\item[a)]
For $\Ring$ a (not necessarily commutative) ring,
$\FEAS_{\IZ,\Ring}$ and $\QUAD{\Ring}$
are polynomial-time equivalent, independently of $\Ring$.
\item[b)]
The matrix $\ast$-ring $\IF^{d\times d}$ is \emph{formally real}
in the sense that $\sum_{j=1}^J B_j^\adjoint\cdot B_j=0$
implies $B_j=0$ ($1\leq j\leq J)$.
\item[c)]
For $\Ring$ a formally real $\ast$-ring,
$\FEAD_{\IZ,\Ring}$ and $\QUART{\Ring}$
are polynomial-time equivalent, independently of $\Ring$.
\end{enumerate}
\end{mylemma}
\begin{proof}
\begin{longenum}
\item[a)] follows straightforwardly using the relational interpretation
Example~\ref{x:RingInterp} in connection with
Observation~\ref{o:MultChain}.
\item[b)]
Suppose there is some $k$ and some $\vec x\in\IF^d$
such that $B_k\vec x\neq0$. Then
$0=\big\langle\vec x\big|\sum\nolimits_j B_j^\adjoint\cdot B_j\vec x\big\rangle
= \sum\nolimits_j {\langle B_j\vec x|B_j\vec x\rangle}
\geq \|B_j\vec x\|^2 > 0$: contradiction.
\item[c)]
Similarly to a), Observation~\ref{o:RelInterp} yields in polynomial time
a system of quadratic $\ast$-polynomials $q_n$ ($1\leq n\leq N)$
whose joint feasibility
is equivalent to that of the original ones. From these, proceed to the
single $\ast$-polynomial $\sum_{n=1}^N q_n^*\cdot q_n$.
Similarly to Fact~\ref{f:RingInterp}, the latter can be achieved
to have coefficients bounded by the number of `literals'
$X_j$ and $X_j^\dagg$.
\qed\end{longenum}\end{proof}

\begin{myproposition}  \lab{p:feas} 
Abbreviate $\FEAS_{\IZ,\IF^*}:=\FEAS_{\IZ,\{\IF^{d\times d}:d\in\IN\}}$,
$\QUAD{\IF^*}:=\QUAD{\{\IF^{d\times d}:d\in\IN\}}$,
$\FEAD_{\IZ,\IF^*}:=\FEAD_{\IZ,\{\IF^{d\times d}:d\in\IN\}}$,
and $\QUART{\IF^*}:=\QUART{\{\IF^{d\times d}:d\in\IN\}}$.
\begin{enumerate}
\item[a)] For any fixed field $\IF$ and $d \geq 1$,
the three problems $\FEAS_{\IZ,\IF}$  and $\FEAS_{\IZ,\IF^{d \times d}}$  
and $\QUAD{\IF^{d\times d}}$ are mutually polynomial-time equivalent.
\item[b)]
$\FEAS_{\IZ,\IF^*}$ is polynomial-time equivalent to $\QUAD{\IF^*}$.
\item[c)] 
For algebraically closed
$\IF$ (say, $\IF=\IC$),   $\FEAS_{\IZ,\IF}$ reduces  polynomially 
 to $\FEAS_{\IZ,\IF^*}$.
\item[d)]
For an $\ast$-subfield $\IF$ of $\IC$ and $d \geq 1$,
the following are polynomial time equivalent: \\
$\FEAD_{\IZ,\IF}$, \; $\FEAS_{\IZ,\Re\IF}$, \;
$\FEAD_{\IZ,\IF^{d \times d}}$, \; $\QUART{\IF^{d\times d}}$.
\item[e)]
For an $\ast$-field $\IF$,
$\FEAD_{\IZ,\IF^*}$ is polynomial-time equivalent to $\QUART{\IF^*}$.
\item[f)] 
If $\IF$ is real or algebraically closed
(e.g. $\IF=\IR,\IC$) then there is 
a polynomial time reduction from $\FEAD_{\IZ,\IF}$ to
  $\FEAD_{\IZ,\IF^*}$. 
\end{enumerate} 
\end{myproposition} 
\begin{proof}  
\begin{longenum}
\item[a)] 
In view of Lemma~\ref{l:RingInterp2}a), it remains to
show the equivalence of $\FEAS_{\IZ,\IF}$ and $\FEAS_{\IZ,\IF^{d\times d}}$.
Concerning one direction, 
a polynomial-time nondeterministic \BSS machine over $\IF$ can
guess the entries of a matrix assignment over $\IF^{d\times d}$
and verify them by evaluating the polynomials using operations
from $\IF$. This demonstrates 
$\FEAS_{\IZ,\IF^{d\times d}}\in\calBP(\calNP_\IF)$.
Since $\FEAS_{\IZ,\IF}$ is complete for $\calBP(\calNP_\IF)$,
the first reduction follows. \\
More directly, the
relational interpretation (Observation~\ref{o:RelInterp}) 
expresses matrix calculations as a system of equations over $\IF$.
\\
For the converse reduction,
apply Fact~\ref{f:matunit} 
requiring, in addition to the given equations
the existence of a $d$-system of matrix units commuting with
all variables. 
\item[b)] Apply Lemma~\ref{l:RingInterp2}a) 
to $\Ring=\IF^{d\times d}$, uniformly in $d$.
\item[c)]
Given equations $p_1(\bar X)=\ldots=p_k(\bar X)=0$ 
add all commutativity conditions  $X_iX_j-X_jX_i$.
Clearly, if $a_1, \ldots a_n$ is a solution in $\IF$ 
for the old system, then $a_1I, \ldots a_nI$ 
is a solution in $\IF^{d \times d}$ for any $d$.
Conversely, if $A_1, \ldots ,A_n$ is a solution
of the extended system in  some $\IF^{d \times d}$,
then the $A_1, \ldots ,A_n$ commute with each other
whence have a common eigenvector $\vec v$ 
(cmp. \mycite{p.109}{Gelfand}, \mycite{Ch.IV.9}{Jacobson})
to respective eigenvalues $\lambda_1, \ldots ,\lambda_n\in\IF$;
which constitute a solution in $\IF$.
\item[d)]
To reduce $\FEAD_{\IZ,\IF}$ to $\FEAS_{\IZ,\Re\IF}$,
use the relational interpretation and  the description
of operations in $\IF= (\Re \IF)(i)$ in  terms of the Cartesian
representation $x=x_1+ix_2$  e.g. replacing
$X=YZ$ by     $X_1=Y_1Z_1-Z_2Z_2=0$ and $X_2= Y_1Y_2 +Y_2Z_1=0$.   
For the converse, add the condition $X_j-X_j^\dagg=0$ for 
each  variable.   
The remaining equivalences follows as in a)
using $\ast$-matrix units.  
\item[e)] Apply Lemma~\ref{l:RingInterp2}b+c) 
to the  $\ast$-ring $\IF^{d\times d}$, uniformly in $d$.
\item[f)] 
Add the commutativity conditions 
and also conditions  $X_j=X_j^\dagg$,
that is, require the $X_j$ to be selfadjoint.
Having a solution $A_1, \ldots A_n$, it holds also
in the real case that 
each invariant subspace $U$ of  $A_i$  contains
an eigenvector of $A_i$ (since the endomorphism
defined by $A_i$ is selfadjoint, also when  restricted to $U$).
Thus, by the same proof as for $\IF=\IC$ in c), the $A_i$ have some 
common eigenvector and the proof proceeds as before
(cmp. \mycite{Ch.VI. Thm.7}{Jacobson}.
\qed\end{longenum}\end{proof}

\subsection{Interpreting the Matrix Rings $\IF^{m \times m}$ into $\Gr(\IF^{3m})$, Uniformly in $m$}  
\lab{ss:MatinGrH}
The next step is to recover
the ring $\IF^{n \times m}$  in  $\Gr(\IF^{3m})$,
a subject dealt with exhaustively  by \person{von Neumann},
who also derived  existence and uniqueness
of coordinatizing $\ast$-rings in much more
general context. Again, we first 
use general frames in order to 
interpret the matrix ring with\emph{out} involution;
then (Subsection~\ref{ss:adjoint}) take also the
latter into account based on orthonormal frames.
The intuitions and motivations from 
Section~\ref{s:Arith2QL} still apply,  though under a more abstract view.

In the first step, we  consider 
   the $\Hom(U,V)$ 
   as $\IF$-vector spaces  
and may allow  infinite dimension.
For subspaces $U,V$ of a vector space $\calH$ define
\[
\IntRing_{U,V} := \{X \in \Gr(\calH) \mid X \cap V=\Zero,\;X+V=U+V\}
\quad\text{ and }\quad
\Gamma_{U,V}(\phi):=\{\vec u - \phi(\vec u) \mid \vec u \in U\}
\enspace , \]
for maps $\phi: U \rightarrow V$, that is, the graph of $-\phi$.

\begin{fact} \lab{endoring} 
Consider  $U \cap V=\Zero$ in $\Gr(\calH)$.
\begin{enumerate}
\item[ a)] 
$\Gamma_{U,V}$ is a bijection of $\Hom(U,V)$ onto 
 $\IntRing_{U,V}$.
\item[ b)] $\phi\in \Hom(U,V)$  is an isomorphism
 if and only if $\Gamma_{U,V}(\phi)   \in \IntRing_{V,U}$.
 Moreover, $\Gamma_{V,U}(\phi^{-1})
=\Gamma_{U,V}(\phi)$ in this case. 
\item[ c)]  For a direct sum $U \oplus V \oplus W$ of 
subspaces of $\calH$ and $\phi\in \Hom(U,V)$ and
$\psi \in \Hom(V,W)$
\[ \Gamma_{U,W}(\psi \circ \phi) = (U+W)\cap \bigl(
\Gamma_{U,V}(\phi) +\Gamma_{V,W}(\psi) 
\bigr) \]  
\item[ d)]
 $\bar A=(A_{ij}\mid 1\leq i,j \leq d)$ is a 
$d$-frame  of $\Gr(\calH)$  if and
only if there is a $d$-system $\bar \vep$ of matrix units such that
$A_i=A_{ii}=\range(\vep_{ii})$ 
and   for all pairwise distinct $i,j$
\[A_{ij}=\{ x -\vep_{ij}(x) \mid x \in A_i\}.\] 
In particular,   one has the maps
$\eta_{ij}$ of Fact~\ref{f:matunit}b). 
\end{enumerate} 
\end{fact} 
\begin{proof} a)-c) are well known  exercises in Linear Algebra.
  By hypothesis, $U+V=U\oplus V$ is a direct
sum, i.e. the $\vec x \in U+V$ have
unique representation $\vec x= \vec u +\vec v$ with
$\vec u \in U$ and $\vec v \in V$. 
\begin{longenum}
\item[a)]
Consider $\phi\in \Hom(U,V)$
and $X=\Gamma_{U,V}(\phi)$.
By linearity, $X \in \Gr(\calH)$. 
Moreover, if $\vec u -\phi(\vec u) \in V$
then $ \vec u =\vec 0$ whence $\phi(\vec u) =\vec 0$ 
by linearity. This proves $X\cap V=0$.
On the other hand, $\vec u = \vec x +\vec v$ 
with $\vec x = \vec u -\phi(\vec u)  \in X$ and
$\vec v =\phi(\vec u)$ whence $U \subseteq X+V$ 
and so $U+V=X+V$. 
\\
Conversely, consider $X \in \IntRing_{U,V}$.
Then $X+V=X\oplus V$ is also a direct sum.
Since $U \subseteq X+V$,
for each $\vec u \in U$ 
there are unique $\vec x \in X$ and
$\vec v \in V$ such hat $\vec u =\vec x +\vec v$.
Define  $\phi(\vec u)=\vec v$.  
Then $X=\Gamma_{U,V}(\phi)$ by definition.
Linearity of $\phi$ follows from the assumption
that $X$ is a linear subspace of $\calH$.
\item[b)]
If the inverse $\phi^{-1} \in \Hom(V,U)$ 
exists, apply a) to it and observe that 
$\vec u -\phi(\vec u) = -(\vec v -\phi^{-1}(\vec v))$
where $\vec v =\phi(\vec u)$ resp.
$\vec u = \phi^{-1}(\vec v)$.
Assuming $X=\Gamma_{U,V}(\phi) \in \IntRing_{V,U}$  
apply a) to conclude that $X=\Gamma_{V,U}(\psi)$ for some
$\psi\in \Hom(V,U)$  and use the preceding argument
to derive $\psi =\phi^{-1}$.
\item[c)] 
\begin{minipage}[c]{0.9\textwidth}
\[\begin{array}{lclllllll} 
\Gamma_{U,V}(\phi)+\Gamma_{V,W}(\psi)&=&\{ \vec u&-\phi(\vec u)+\vec v
&-\psi(\vec v)& \mid \vec u \in U,&\vec v \in V&&\}\\
U+W&=&\{\vec x&&-\vec w&\mid \vec x \in U,&&\vec w \in W&\}
\end{array}\]
\end{minipage} \\
whence, by uniqueness of summands
$\vec u \in U$, $\vec v \in V$, and
$\vec w \in W$ 
in $\vec x= \vec u+\vec v +\vec w$, 
 the intersection is obtained as
$\{ \vec u
-\psi(\vec v) \mid \vec u \in U,\,\vec v \in V,\; \vec v= \phi(\vec
u)\}
= \Gamma_{U,W}(\psi \circ \phi)$.
\item[d)] That $\calH$ is the  direct sum of the
$A_i$ is encoded into the first two conditions
in Definition~\ref{d:frame} .
The remaining ones and a)-c)  give the $\eta_{ij}$ as 
 in Fact~\ref{f:matunit}b).  
\qed\end{longenum}\end{proof}    
Guided by the development in Subsection~\ref{ss:IntField}
and following von  Neumann \cite{Neumann2}, one is led to

\begin{mydefinition} \lab{d:MatinGrH}
For any $d \geq 3$ and variables $x,y$, and 
$\bar y =(y_{ij}\mid 1\leq i,j \leq d$), define
the following lattice formulae and terms. 
\[\begin{array}{lcl} 
\rho_{\bar y }(x)& :\Leftrightarrow & x \wedge y_2=\Zero\; \Band \;x \vee
y_2= y_1 \vee y_2\\
\pi^{ij}_{{\bar y} ik}(x)&:=&(y_i\vee y_k)\wedge (x\vee y_{jk})\\    
\pi^{ij}_{{\bar y} kj}(x)&:=&(y_k\vee y_j)\wedge (x\vee y_{ik})\\
x\ominus_{\bar y} z&:=& 
\bigl( [(\pi^{12}_{{\bar y}13}(z) \vee y_2) \wedge  (x \vee y_{23})] \vee y_3\bigr) \wedge
(y_1 \vee y_2)\\  
\ominus_{\bar y} z&:=& y_1 \ominus_{\bar y} z\\
x\oplus_{\bar y}  z&:=& x\ominus_{\bar y} (\ominus_{\bar y} z)\\ 
x \otimes_{\bar y} z&:=&(\pi^{12}_{{\bar y}13}(z) \vee  \pi^{12}_{{\bar y}32}(x)) \wedge   (y_1\vee y_2)
\end{array} \]
\end{mydefinition} 
For $d$-frame $\bar A$
of $\Gr(\calH)$  and $\phi: A_i \rightarrow A_j$ abbreviate
$\IntRing_{\bar A ij}:= \IntRing_{A_i,A_j}$, $\Gamma_{\bar Aij}(\phi):=\Gamma_{A_i,A_j}(\phi)$.

\begin{fact} \lab{f:MatinGrH}
\begin{enumerate}
\item[a)]
Fix a  $d$-frame $\bar A$ in $\Gr(\calH)$ 
with  $\eta_{ij}$ as in Fact~\ref{endoring}d).
For any pairwise distinct $i,j,k \leq 3$, the polynomials
$\pi^{ij}_{\bar A ik}$  and  
$\pi^{ij}_{\bar A kj}$  of Definition~\ref{d:MatinGrH} provide
bijections from $\IntRing_{\bar A ij}$
onto $\cal R_{\bar A ik}$ and $\IntRing_{\bar A kj}$, respectively.
Moreover, for  
$\phi\in \Hom(A_i,A_j)$,
\[   \pi^{ij}_{\bar A ik}\big(\Gamma_{\bar A ij}(\phi)\big)= \Gamma_{\bar A ik}(\eta_{jk} \circ
\phi),\quad
\pi^{ij}_{\bar A kj} \big(\Gamma_{\bar A ij}(\phi)\big)= \Gamma_{\bar A kj}( \phi \circ\eta_{ki})\]
\item[ b)] 
The map $\Interpret_{\bar A}$, $\phi\mapsto
 \Gamma_{\bar A 12}(\eta_{12} \circ \phi)$,
is an isomorphism 
from the ring 
$\End(A_1)$  with operations $+$, $-$, $0$, $\circ$, $\id$
onto $\IntRing_{\bar A 12}$ with operations
$\oplus_{\bar A}$, 
$\ominus_{\bar A}$, $A_1$, $\otimes_{\bar A}$, $A_{12}$.
\item[c)] 
To every $d\geq3$ and every finite family $\bar p$ of
integer polynomials $p_j\in\IZ\langle X_1,\ldots,X_n\rangle$ 
($1\leq j\leq J)$
in $n$ non-commuting variables, there exists a 
term $t_{\bar p,d}$ in the language of ortholattices such that
the following holds: \\
For every field $\IF$ and every $\IF$-unitary vector space $\calH$,
$t_{\bar p,d}$ is strongly satisfiable over $\Gr(\calH)$
~iff~ there exists $m\in\IN$ and $W\in\Gr_m(\calH)$
and $\phi_1,\ldots,\phi_n\in\End(W)$ such that 
$\dim(\calH)=md$ and $p_j(\phi_1,\ldots,\phi_n)=0$ for every $1\leq j\leq J$.
\end{enumerate}
\end{fact}
Recall that Claim~c) has already been employed in
the proof of Proposition \ref{p:Fields}f). 
\begin{proof}  
\begin{longenum}
\item[a)]
and bijectivity of $\Interpret_{\bar A}$  follow immediately from Fact~\ref{endoring}.
\item[b)]
Consider
$\phi,\psi \in \End(A_1)$. By a) we have
\[ \pi^{12}_{\bar A 32}(\Interpret_{\bar A}(\phi))= \Gamma_{\bar A 32}(\eta_{12} \circ \phi 
\circ \eta_{31}),\quad
  \pi^{12}_{\bar A 13}(\Interpret_{\bar A}(\psi))= \Gamma_{\bar A 12}( \eta_{23}
  \circ\eta_{12} \circ \psi) =  \Gamma_{\bar A 12}( \eta_{13}  \circ \psi) \]  
and with c) of Fact~\ref{endoring} it follows
$\Interpret_{\bar A}(\phi \circ \psi) =\Interpret_{\bar A}(\phi) \otimes_{\bar A} \Interpret_{\bar A}(\psi)$.   
Moreover
\[\begin{array}{lcllclcllll}  \Interpret_{\bar A}(\phi) + A_{23}&=&
\{\vec x &+&\vec y -\eta_{12}(\phi(\vec x)) &-& \eta_{23}(\vec y)
&\mid \vec  x \in A_1,&\vec y \in A_2&\} 
\\
\Interpret_{\bar A}(\psi) + A_{2}&=&
\{\vec u &+&\vec v  &-&\eta_{13}(\psi(\vec u)) 
&\mid \vec  u \in A_1,&\vec v \in A_2&\}
\end{array}\] 
whence for $S:=  (\Interpret_{\bar A}(\psi) + A_2)\cap (\Interpret_{\bar A}(\phi) + A_{23})$
one obtains
\[ S= \{\vec x +\vec y -\eta_{12}(\phi(\vec x)) - \eta_{23}(\vec y)
\mid \vec  x \in A_1,\;\vec y \in A_2,\;\eta_{23}(\vec y)=
\eta_{13}(\psi(\vec x)) \} \]
Now,  applying $\eta_{32}$ we see that 
$\eta_{23}(\vec y)= \eta_{13}(\psi(\vec x))$ 
is equivalent to $\vec y= \eta_{12}(\psi(\vec x))$ and 
\[ \begin{array}{lclclcll}
 S&=& \{\vec x &+& \eta_{12}(\psi(\vec x)) -\eta_{12}(\phi(\vec x)) &-& \eta_{23}(\vec y)
&\mid \vec  x \in A_1,\;\vec y \in A_2 \}\\
&=& \{\vec x &-& \eta_{12}((\phi -\psi)(\vec x))  &-& \eta_{23}(\vec y)
&\mid \vec  x \in A_1,\;\vec y \in A_2 \}
\end{array} \]
whence
$\Interpret_{\bar A}(\phi)\ominus_{\bar A}
\Interpret_{\bar A}(\psi)= (S+A_3)\cap (A_1+A_2)= \Interpret_{\bar A}(\phi -\psi)$.
\item[c)]
Similar to the proof of Proposition~\ref{p:IntField}c),
let $t_{\bar p,d}(\bar y,\bar x)$ be the conjunction of 
conditions requiring $\bar y$ to evaluate to a $d$-frame $\bar A$
(Definition~\ref{d:frame}), all variables $x_j$ to 
belong to the induced ring $\IntRing_{\bar A 12}$, 
integer constants replaced by terms (Observation~\ref{o:MultChain}),
and operations $+,-,0,\circ,\id$ by 
$\oplus_{\bar A},\ominus_{\bar A},A_1,\otimes_{\bar A},A_{12}$,
respectively.
Note that, according to Fact~\ref{endoring}, 
a $d$-frame exists precisely in dimensions an integer multiple of $d$.
\qed\end{longenum}\end{proof} 
\begin{scholium} \lab{s:MatinGrH}
The term $t_{\bar p,d}$ from Fact~\ref{f:MatinGrH}c)
can be computed from $\bar p$
and $d$ by a Turing machine in time polynomial in $d$
and in the binary encoding length of $\bar p$. 
\end{scholium}

\subsection{Uniformly Interpreting the $\ast$-Rings 
$\IF^{m\times m}$ into $\Gr(\IF^{3m})$}\lab{ss:adjoint} 
In order to include  adjunction into the above result,
we have to relate adjoint pairs of linear maps between
subspaces of $\calH$ to the ortholattice $\calH$.

\begin{fact} \lab{f:adjoint} 
Consider 
$U,V \in \Gr(\calH)$, $\calH$ a finite-dimensional $\IF$-unitary space.
\begin{enumerate} 
\item[a)]
If $U \perp V$ then
$\phi$ and $\psi$ are adjoint to each other
if and only if $\Gamma_{U,V}(\phi) \perp \Gamma_{V,U}(-\psi)$.
\COMMENTED{
\item[h)] Fact~\ref{f:matunit} 
is valid mutatis mutandis, i.e. $\IF$-vector spaces to be $\IF$-unitary
idempotents are required to be projections
 direct sums to be orthogonal, bases to be
equinormal orthogonal, and linear   isomorphisms  to be isometries;
in particular, $\eta_{ij}^*= \eta_{ij}^{-1} =\eta_{ji}$.  } 
\item[b)] 
 $\bar A=(A_{ij}\mid 1\leq i,j \leq d)$ is an orthonormal
$d$-frame  of $\Gr(\calH)$  if and
if and only if  the $\eta_{ij}$ in Fact~\ref{endoring}d)
are isometries, in other words, $\eta_{ij}^{-1}=
\eta_{ji}=\eta_{ij}^\ast$.
\item[c)] If $d$ divides $\dim(\calH)$ and $\IF$ is Pythagorean,
there is an orthonormal $d$-frame of
$\Gr(\calH)$.
\end{enumerate} 
\end{fact} 
Put differently, b)  establishes a $1$-$1$-correspondence
between orthonormal $d$-frames of $\Gr(\calH)$  and
$d$-systems of matrix units of the $\ast$-ring $\End(\calH)$.

\begin{proof} 
Concerning a), $U \perp V$ implies
\[ \langle \vec u - \phi(\vec u) \mid 
\vec v  +\psi(\vec v) \rangle 
=
\langle \vec u \mid \psi(\vec v) \rangle -
\langle \phi(\vec u)\mid \vec v\rangle 
\;\;\mbox{ for all } \vec u \in U,\,\vec v \in V\]
proving the claim.
For b)   observe that by Fact~\ref{endoring}e)
$A_{12}= \Gamma_{\bar A 12} (\id)$ and 
$\ominus_{\bar A12} A_{12} =\Gamma_{\bar A 12}(- \id)$. These being orthogonal
to each other means  $
 \langle \vec x \mid \vec y \rangle 
= \langle \eta_{12}(\vec x) \mid \eta_{12}(\vec y) \rangle$
since, in view of  $A_1 \perp A_2$,  
 for $\eta=\eta_{12}$ and
all $\vec x,\vec y \in A_1$  
\[ \langle \vec x - \eta(\vec x)   \mid \vec y +\eta(\vec y)\rangle       
=  \langle \vec x \mid \vec y  \rangle
+\langle \vec x  \mid \eta(\vec y)  \rangle
- \langle \eta(\vec x) \mid \vec y \rangle
-\langle \vec x \mid \eta(\vec y) \rangle
= \langle \vec x \mid \vec y \rangle 
- \langle \eta(\vec x) \mid \eta(\vec y) \rangle.\] 
Symmetry applies to the $\eta_{1k}$.
The claim for the $\eta_{ij}$ follows from 
$\eta_{ij}= \eta_{1j} \circ \eta_{1i}^{-1}$.  
c) ) By virtue of Fact~\ref{f:GramSchmidt}b)
choose an orthonormal basis $\vec v_{ik}$, $i=1, \ldots ,d$,
$k=1, \ldots, n$ and put 
 $A_i=\sum_{k=1}^n \IF\vec v_{ik}$, $A_{ij}=   \sum_{k=1}^n \IF(\vec
 v_{ik}
-\vec v_{jk})$.
\end{proof}
Now, suppose that an orthonormal $d$-frame $\bar A$ of 
$\Gr(\calH)$ is given
as in Fact~\ref{f:adjoint}b) 
for finite dimensional $\calH$.   For $\phi \in \End(A_1)$ we want
to recover $\phi^*$ in terms of   $\Gr(\calH)$.
Write $\Gamma_{\bar A ij}=\Gamma_{A_1,A_j}$  
and consider $\psi \in \End(A_1)$.  We claim that
\[ \psi = \phi^* 
\;\mbox{ if and only if } \;  
\Gamma_{\bar A21}( \psi \circ \eta_{21}) = (A_1 + A_2) \cap  (\Gamma_{\bar A12}
( \eta_{12} \circ \phi))^\bot.\] 
Indeed, by  Fact~\ref{f:adjointa}b+f) 
it holds that
\[ \psi=\phi^* ~\mbox{ if and only if }~
\psi \circ \eta_{21} =  \phi^* \circ \eta_{21}
= \phi^* \circ \eta_{12}^* = (\eta_{12} \circ \phi)^*.\]
By Fact~\ref{f:adjoint}b)  the latter   amounts to  
$\Gamma_{\bar A21}( \psi \circ \eta_{21}) \subseteq  (A_1 + A_2) \cap  \big(\Gamma_{\bar A12}
(-\eta_{12} \circ \phi)\big)^\bot$.
Here, equality holds since both sides are
complements of $A_1$ in $[\Zero,A_1+A_2]$ ---
with $X= \Gamma_{\bar A12}
( \eta_{12} \circ -\phi)$ one derives
$
 A_1 + X^\bot = A_1^{\bot\bot} +X^\bot  
=(A_1^\bot \cap X)^\bot =( \sum_{i=2}^dA_i \cap X)^\bot
=\Zero^\bot  =\One$ whence $A_1+\bigl( (A_1+A_2)\cap X^\bot\bigr) 
= A_1+A_2$. 

Thus, it remains to recover  $\Gamma_{\bar A12}(\eta_{12}\circ \psi)$ 
from $\Gamma_{\bar A21}(\psi \circ \eta_{21})$ by means
of lattice operations.   This is easily done
using Fact~\ref{f:MatinGrH}a):
\[ \eta_{12} \circ     \psi =  \eta_{12} \circ \psi \circ \eta_{11}
= \eta_{12} \circ \psi \circ \eta_{21} \circ \eta_{32} \circ
\eta_{13}
= \bigl( \eta_{12} \circ ((\psi \circ \eta_{21}) \circ \eta_{32}\bigr)
\bigr) \circ
\eta_{13} \]
implies 
$\Gamma_{\bar A12}( \eta_{12} \circ     \psi)
= t_{\bar A}\big( \Gamma_{\bar A21}( \psi \circ \eta_{21})\big)$
where  
 $t_{\bar y}(x)$ is defined as 
\begin{equation}
\label{e:Dagger}
t_{\bar y}(x) \;:= \; \pi^{32}_{\bar y12}\big(\pi^{31}_{\bar y
  32}[\pi^{21}_{\bar y31}(x)]\big)
\qquad\text{and}\qquad
 x^{\dagger_{\bar y}} \;:=\;
t_{\bar y}\big((y_1 \vee y_2) \wedge \neg(\ominus_{\bar y}  x )\big)
\end{equation}
captures adjunction.
So, together with Fact~\ref{f:MatinGrH}b) 
and Scholium~\ref{s:MatinGrH} we have proved

\begin{mylemma} \lab{l:MatinGrH} 
\begin{enumerate}
\item[a)]
Let $\IF$ be Pythagorean, $\calH$ a
finite dimensional $\IF$-unitary space, and
$d =\dim \calH \geq 3$.
 Given  an orthonormal $d$-frame $\bar A$ of $\Gr(\calH)$,
there 
is an isomorphism 
from the $\ast$-ring 
$\End(A_1)$  with operations $+$, $-$, $0$, $\circ$, $\id$, $^\ast$ 
onto $\IntRing_{\bar A 12  }=\{U \in \Gr(\calH)\mid \rho_{\bar A}(U)\}$ with
operations
defined by the terms
$x\oplus_{\bar A}z$, 
$\ominus_{\bar A}z$, $A_1$, $x\otimes_{\bar A}z$, $A_{12}$,
$x^{\dagger_{\bar A}}$
 with constants 
$\bar A =(A_{ij}\mid 1\leq i,j \leq 3)$ 
from $\Gr(\calH)$,  
 as given by Definition~\ref{d:MatinGrH} and Equation~(\ref{e:Dagger}).
\item[b)]
Given $d\geq3$ and a finite family $\bar p$ of
$\ast$-polynomials $p_j\in\IZ\langle X_1,X_1^\dagger,\ldots,X_n,X_n^\dagger\rangle$ 
($1\leq j\leq J)$
in $n$ noncommuting variables, 
a polynomial-time Turing machine can produce a 
term $t_{\bar p,d}$ in the language of ortholattices such that
the following holds: \\
For every Pythagorean $\ast$-field $\IF$, and every $\IF$-unitary vector space $\calH$,
$t_{\bar p,d}$ is strongly satisfiable over $\Gr(\calH)$
~iff~ $\dim(\calH)=dm$ and 
there exist $\phi_1,\ldots,\phi_n\in\End(\IF^m)$ such that 
$\dim(\calH)=md$ and $p_j(\phi_1,\phi_1^\ast,\ldots,\phi_n,\phi_n^\ast)=0$ 
for every $1\leq j\leq J$.
\end{enumerate}
\end{mylemma} 
Item b) thus generalizes Theorem~\ref{t:IntStarField}a)  to the noncommutative
case, uniformly in $m$.

\subsection{Strong Satisfiability in Indefinite Finite Dimension}
\lab{ss:StrongIndef}
Picking up on Theorem~\ref{t:WeakIndef}a),
and since Theorem~\ref{t:WeakStrong}
depends on the dimension being fixed, we now 
consider the complexity and computability of 
\emph{strong} satisfiability over indefinite finite dimensions
--- and approach the boundary between complexity and mere computability.

\begin{myproposition} \lab{p:StrongIndef} For  a fixed
  $\ast$-subfield $\IF$ of  of $\IC$ there are
the following  polynomial time reductions:
\begin{enumerate}
\item[a)] $\SAT$ to  $\SAT_{\Gr(\IF^*)}$ 
to $\FEAD_{\IZ,\IF^*}$.
\item[b)]
\COMMENTED{
$\sat_{\Gr(\IF^*)}$ to  any of $\SAT_{\Gr(\IF^*)}$,  $\FEAD_{\IZ,\IF^*}$, and
 $\FEAD_{\IZ,\IF}$. 
\item[c)] 
Moreover, }
For Pythagorean $\IF$,  $\FEAD_{\IZ,\IF^*}$  to $\SAT_{\Gr(\IF^*)}$.
\item[c)] 
For real resp. algebraically closed $\IF$, 
$\FEAD_{\IZ,\IF}$  to $\SAT_{\Gr(\IF^*)}$.
\end{enumerate} 
\end{myproposition} 

\begin{proof}
a) The first reduction is Proposition~\ref{p:NPhard}.
According to Proposition~\ref{p:UpperComplexity}c),
$\SAT_{\Gr(\IF^d)}$ reduces in polynomial time 
to $\FEAS_{\IZ,\Re\IF}$ and, by Fact~\ref{f:BSSBP}a),
furtheron to any problem in $\calBP(\calNP_{\Re\IF})$ --
such as $\FEAD_{\IZ,\IF^{d\times d}}$. However both
reductions depend on $d$ while  
  Fact~\ref{f:lat}  provides
a  short-cut uniformly in $d$.   
\COMMENTED{
\item[b)] 
For any  ortholattice term $t$
let $\chi(t) =\phi_{|t|}(t)$  
according to Theorem~\ref{t:WeakStrong}a),
translate this into a  solvability condition $\Sigma_t$ for
$\ast$-rings as in a).
In that,   replace each (matrix)  variable $X$ by a 
$|t|\times|t|$ matrix $(x_{ij})$ of (coefficient)   variables
and express matrix operations via  these coefficients
to obtain a solvability  condition $\Sigma_t^\#$.
Observe that one has 
to consider the dimensions $d$ given in unary
representation  in order to obtain uniformity.
Here, we have $d=|t|$, so this representation is inherent in
the term  $t$.  
Then, 
the following are equivalent: 
  Weak satisfiability  of $t$ in $\Gr(\IF^*)$;
 Weak satisfiability  of $t$ in $\Gr(\IF^{|t|})$ 
(by Lemma~\ref{l:Dimensions}b);
 Strong  satisfiability  of $\chi(t)$ in $\Gr(\IF^{|t|})$
(by Theorem~\ref{t:WeakStrong}a);
Feasibility of $\Sigma_t$ in $\IF^{|t|\times|t|}$
(by a);
Feasibility of  $\Sigma_t^{\#}$  in $\IF$ (by the proof of Proposition~\ref{p:feas}d).   
This proves the reductions claimed in b).
}
b) 
follows  from Lemma~\ref{l:MatinGrH}b).
For c), combine b) with Proposition~\ref{p:feas}f).
\qed\end{proof} 
\begin{theorem} \lab{t:StrongIndef}
$\SAT_{\Gr(\IR^*)}=\SAT_{\Gr(\IC^*)}$ is $\calBP(\calNP_\IR)$--hard.
Moreover, the following are equivalent:
\begin{enumerate}
\item[i)] $\SAT_{\Gr(\IC^*)}$ is decidable
\item[ii)] $\FEAD_{\IZ,\IC^*}$ is decidable
\item[iii)] $\QUART{\IR^*}$ is decidable 
\item[iv)] There exists a total recursive function
$\delta:\IN\to\IN$ such that the following holds:
\begin{equation} \label{e:StrongIndef1}
\begin{minipage}{0.7\textwidth}
Whenever a term $t$ is strongly satisfiable over $\Gr(\IC^*)$, \\
it is so over $\Gr(\IC^d)$ for some $d\leq\delta(|t|)$.
\end{minipage}
\end{equation}
\item[v)] There exists a total recursive function
$\delta':\IN\to\IN$ such that the following holds:
\begin{equation} \label{e:StrongIndef2}
\begin{minipage}{0.87\textwidth}
Whenever a quartic 
$p\in\{0,\pm1,\pm2,\ldots,\pm2n\}\langle X_1,X_1^\dagger,\ldots,X_n,X_n^\dagger\rangle$
admits \\ a root over $\IR^{d\times d}$ for some $d$,
it also does so over $\IR^{d'\times d'}$ for some $d'\leq\delta'(n)$.
\end{minipage}
\end{equation}
\end{enumerate}
\end{theorem}
Note that the following functions $\delta,\delta'$ do
satisfy Equations~(\ref{e:StrongIndef1}) and (\ref{e:StrongIndef2}), 
respectively: 
\begin{align}
\delta:n\mapsto\max\big\{d\;\big|\;&\forall 
\text{ ortholattice terms } t: \\[-0.5ex]
& \big(|t|=n \;\wedge\;
\exists d'\in\IN \; \exists A'_1,\ldots,A'_n\in\Gr(\IR^{d'}): \nonumber
t_{\Gr(\IR^d)}(A'_1,\ldots,A'_n)=\One \big) \\[-0.5ex]
&\quad\qquad\;\Rightarrow\;  \label{e:StrongIndef3}
\exists d'\leq d \; \exists A_1,\ldots,A_n\in\Gr(\IR^{d'}):
t_{\Gr(\IR^{d'})}(A_1,\ldots,A_n)=\One \big\} \nonumber
\\
\delta':n\mapsto\max\big\{d'\;\big|\;&\forall 
p\in\{0,\pm1,\pm2,\ldots,\pm2n\}\langle X_1,X_1^\dagger,\ldots,X_n,X_n^\dagger\rangle: \\[-0.5ex]
& \big(\exists d\in\IN \; \exists A_1,\ldots,A_n\in\IR^{d\times d}: \label{e:StrongIndef4}
p(A_1,A_1^\adjoint,\ldots,A_n,A_n^\adjoint)=0\big) \nonumber \\[-0.5ex]
&\qquad\;\Rightarrow\; \nonumber
\exists d\leq d' \; \exists A'_1,\ldots,A'_n\in\IR^{d\times d}: 
p(A'_1,A_1'^\adjoint,\ldots,A'_n,A_n'^\adjoint)=0 \big\} \enspace .
\end{align}
Indeed, the number of terms of given length is finite;
hence the maximum exists. Similarly, a quartic $n$-variate
polynomial consists of $\calO(n^4)$ monomials;
hence there can be no more than $|K|^{\calO(n^4)}$ of
them with coefficients in $K:=\{0,\pm1,\pm2,\ldots,\pm n\}$.

\begin{myquestion} \lab{q:StrongIndef}
Is the function $\delta'$ 
(unlike the \textsf{busy beaver})
recursively bounded?
\end{myquestion}
According to Example~\ref{x:Dimensions}b),
such a bound has to be at least exponential.

\begin{proof}[of Theorem~\ref{t:StrongIndef}]
$\SAT_{\Gr(\IC^*)}=\SAT_{\Gr(\IR^*)}$ holds due to
Fact~\ref{f:Fields}a+b); for hardness invoke
Proposition~\ref{p:StrongIndef}c).
\begin{longenum}
\item[i)$\Leftrightarrow$ii)] Proposition~\ref{p:StrongIndef}a+b).
\item[i)$\Leftrightarrow$iii)] 
$\SAT_{\Gr(\IC^*)}=\SAT_{\Gr(\IR^*)}$ is equivalent
$\FEAS_{\IZ,\IR^*}$ according to Proposition~\ref{p:StrongIndef}a+b);
which by Proposition~\ref{p:feas}e)
is in turn equivalent to $\QUART{\IR^*}$.
\item[i)$\Rightarrow$iv)] 
Based on an algorithm deciding $\SAT_{\Gr(\IC^*)}$,
the function $\delta$ from Equation~(\ref{e:StrongIndef3})
can be computed as follows: Given $n$ enumerate all 
(the finitely many, up to renaming variables)
terms $t$ of length $n$ strongly satisfiable over
$\Gr(\IC^*)$; for each one 
search (Fact~\ref{f:Tarski})
for the first dimension $t$ is strongly satisfiable in
and return the maximum.
\item[iv)$\Rightarrow$i)] 
Given $t$, calculate $d:=\delta(|t|)$ and 
decide satisfiability of $t$ over
$\Gr(\IC^1)$, $\Gr(\IC^2)$, \ldots, $\Gr(\IC^d)$:
if none succeeds then $t$ is not satisfiable
over $\Gr(\IC^*)$ either.
\item[ii)$\Rightarrow$v)]
Similarly to i)$\Rightarrow$iv),
enumerate all (the finitely many) $n$-variate quartic 
$\ast$-polynomials over $\{0,\pm1,\pm2,\ldots,\pm2n\}$ 
in $\QUART{\IR^*}$; for each search for the first
$d$ such that it admits a root in $\IR^{d\times d}$
and return the maximum.
\item[v)$\Rightarrow$ii)]
Similarly to iv)$\Rightarrow$i),
given $p\in\{0,\pm1,\pm2,\ldots,\pm2n\}\langle X_1,\ldots,X_n,X_1^\dagg,\ldots,X_n^\dagg\rangle$
calculate $d':=\delta'(n)$ and invoke
Proposition~\ref{p:feas}d) to decide whether
$p$ admits a root in $\IR^{1\times 1}$, $\IR^{2\times 2}$,
\ldots, $\IR^{d'\times d'}$:
if none succeeds then $p\not\in\FEAD_{\IZ,\IR^*}$.
\qed\end{longenum}\end{proof}
\cx{
\begin{digression}
For Pythagorean $\IF$,  
the following problem allows a polynomial time  reduction to strong
satisfiability in $\Gr(\IF^*)$.  Namely,
to decide for any finite  semigroup 
resp. group presentation (the latter are considered  as 
semigroup presentations with 
a generator symbol $g'$ and relations $gg'=e=g'g$ 
added for each group generator symbol $g$)
whether it admits a   finite
dimensional $\IF$-algebra model $\calA$.
The proof is by interpreting the relations of the
presentation as in Lemma~\ref{l:MatinGrH} (with $d=3$), 
Then, the ortholattice relations are
combined into strong satisfiability of $t$ 
via Fact~\ref{f:clear}a). If there is a 
model $\calA$  as required, then $t$  will be strongly  satisfied in $\Gr(\calH)$
where $\calH$ is the
$3$-fold orthogonal sum of $\IF[S]$
considered as a right vector space with canonical scalar product
and where  any generator symbol $g$  is interpreted as negative graph
of the linear map $x \mapsto gx$.
Conversely, if $t$ is strongly satisfied, 
then the endomorphisms of $U_1$ corresponding
to generator symbols  
generate a  finite dimensional model  $\calA$. 
\end{digression} }

\section{Towards the Infinite-Dimensional Case} \lab{s:Neumann}

Recall that, for an infinite dimensional Hilbert space $\calH$,
$\Gr(\calH)$ consists of all \emph{closed} linear
subspaces of $\calH$. We also define
$\Grl(\calH)$ as the modular lattice of all
linear subspaces of  the $\IF$-vector space
$\calH$; and 
$\Grf(\calH)$ as the subset of all 
$U$ and $U^\bot$ where $U \in \Lat(\calH)$ is finite
dimensional (and in particular closed).

\begin{fact}
 $\Grf(\calH)$ is a sub-ortholattice of 
$\Gr(\calH)$  and a sublattice of $\Grl(\calH)$,
in particular a modular ortholattice. All claims from
Lemma~\ref{l:Dimensions} remain valid
when applied to  $\Grf(\calH)$. 
The analogous  results hold  for abstract atomic MOLs cf. \cite{hn}. 
\end{fact} 
%
Though, one has to be aware of the fact that
$\Grl(\calH)$ is  much more special than $\Gr(\calH)$ even
 if $\calH$ is a  separable Hilbert space.
In particular, nothing is known about
satisfiability in  $\Gr(\calH)$. Neither 
is it known whether the equational theory of
 $\Gr(\calH)$ is decidable.

Recall that a von Neumann algebra ${\bf A}$
(of bounded operators on a Hilbert  space)  
is \textsf{finite} 
if $AA^\ast=I$ implies $A^\ast A=I$ for all $A \in {\bf A}$.
It is a \textsf{factor} if its center is $\IC$.
Its \textsf{projections} $P=P^\ast=P^2$  are
partially ordered by $P\leq Q \Leftrightarrow QP=P$.  
The following is due to \person{Murray} and \person{von Neumann}
 cmp. \cite{Takesaki}.  

\begin{fact} The projections of a finite
von Neumann algebra ${\bf A}$ form a  
 continuous modular ortholattice $\Neu({\bf A})$.
If ${\bf A}$ is a factor then either
$\Neu({\bf A}) \cong \Gr(\IC^d)$ for some $d <\infty$ 
or $\Neu({\bf A})$ is simple with no atoms;
 the latter are  called of  \emph{type} $\textrm{II}_1$
and contain  all $\Gr(\IC^d)$ as sub-ortholattices.
\end{fact}
\begin{proof} 
Most of this is in 
\cite{mn}. The key to modularity is that
${\bf A}$ extends to a $\ast$-regular ring having the same
projections (\S 14.1). In   
 \mycite{Thm. XIII}{mn4} it is shown that a type $\textrm{II}_1$-algebra
has a subalgebra $\IC^{d \times d}$ for any $d$. 
From this it follows that $\Gr(\IC^{d \times d})$
embeds into $\Neu({\bf A})$ -- cmp. \mycite{Cor.3.2}{hn} 
 \end{proof}  

The $\Neu({\bf A})$ with $\bf A$ of type II$_1$ 
are examples of \emph{continuous geometries}, i.e.
continuous modular ortholattices $L$,
admitting a  dimension function
mapping $L$ onto the unit interval.
In \cite{Neumann3} \person{von Neumann}
introduced additional structure and axioms
on the latter, motivated by his  view on what
a `quantum logic' should be, to  recover these
abstract continuous  geometries as
projection MOLs of algebras of
operators,  cf. \cite{Redei}.

\begin{theorem}[\cite{hn}]
 Given a finite von Neumann algebra factor ${\bf A}$
of type II$_1$, an ortholattice identity
is valid in $\Neu({\bf A})$ if and only if it is valid
in  $\Gr(\IC^d)$  for all  resp. infinitely many $d$.
\end{theorem}

\begin{corollary} \begin{enumerate} 
\item[ a)] The (indefinite)  weak satisfiability problem  for 
any class of  of projection ortholattices of
 finite von Neumann algebra factors  is in $\calBP(\calNP_\IR)$.
\item[ b)] The decision problem  for 
the equational  theory  of  any class of  projection ortholattices
of  
 finite von Neumann algebra factors  is in $\calBP(\co\calNP_\IR)$.
\end{enumerate} 
\end{corollary}

\cx{
\begin{digression}
The concept of
 $\ast$-\emph{regular rings},  $\ast$-rings such that for any $a$
there is $x$ such that $a=axa$ and $aa^*=0$ only for $a=0$,
has been introduced by   \person{von Neumann}
as an abstract approach to  certain algebras of unbounded operators.
He showed that $\ast$-regular rings  have MOLs of projections and 
 that  
any MOL  admitting an orthogonal $d$-frame $\bar a$, $d\geq 4$,
is isomorphic to the lattice of projections
of some $\ast$-regular ring (cf. 
\cite{Neumann2}), actually  the matrix ring $\IntRing_{\bar a}^{d \times d}$
with involution mimicking Fact~{f:adjoint}c). 
This has been
extended to `large partial orthogonal $d$-frames'
by \person{J\'{o}nsson} \cite{jon1}, including the case $d=3$
under the stronger Arguesian Law (cf. \cite{coord})
which is valid in all MOL of projections. In particular, this
applies to all simple Arguesian  MOLs of height $\geq 3$.
\end{digression}} 

\end{document}